\documentclass[10pt,reqno]{amsart}
\usepackage{amssymb,amscd}
\usepackage{color}

\newtheorem{thm}{Theorem}[section]
\newtheorem{prop}[thm]{Proposition}
\newtheorem{lemma}[thm]{Lemma}
\newtheorem{cor}[thm]{Corollary}

\theoremstyle{remark}

\newtheorem{remark}[thm]{Remark}

\theoremstyle{definition}
\newtheorem{definition}[thm]{Definition}

\def\bN{\mathbb{N}}
\def\U{\mathrm{U}}
\def\tr{\mathrm{tr}}
\def\orb{\mathrm{orb}}
\def\1{\mathbf{1}}
\def\eps{\varepsilon}
\def\cX{\mathcal{X}}
\def\bR{\mathbb{R}}
\def\bK{\mathbb{K}}
\def\bP{\mathbb{P}}
\def\cP{\mathcal{P}}
\def\cN{\mathcal{N}}
\def\bC{\mathbb{C}}
\def\bT{\mathbb{T}}

\def\cA{\mathcal{A}}
\def\sym{\mathrm{sym}}
\def\SU{\mathrm{SU}}
\def\HS{\mathrm{HS}}
\def\bX{\mathbf{X}}
\def\bY{\mathbf{Y}}

\def\bv{\mathbf{v}}
\def\bU{\mathbf{U}}
\def\bV{\mathbf{V}}
\def\bA{\mathbf{A}}

\def\<{\langle}
\def\>{\rangle}

\begin{document}

\title[Orbital approach to free entropy]{Orbital approach to microstate free entropy}
\author[F. Hiai]{Fumio Hiai$\,^{1}$}
\address{Graduate School of Information Sciences,
Tohoku University, Aoba-ku, Sendai 980-8579, Japan}
\author[T. Miyamoto]{Takuho Miyamoto}
\address{Graduate School of Information Sciences,
Tohoku University, Aoba-ku, Sendai 980-8579, Japan}
\author[Y. Ueda]{Yoshimichi Ueda$\,^{2}$}
\address{Graduate School of Mathematics,
Kyushu University, Fukuoka 810-8560, Japan}

\thanks{$^1\,$Supported in part by Grant-in-Aid for Scientific Research
(B)17340043.}
\thanks{$^2\,$Supported in part by Grant-in-Aid for Young Scientists
(B)17740096.}

\maketitle

\begin{abstract}
Motivated by Voiculescu's liberation theory, we introduce the orbital free entropy
$\chi_\orb$ for non-commutative self-adjoint random variables (also for
``hyperfinite random multi-variables"). Besides its basic properties the relation of
$\chi_\orb$ with the usual free entropy $\chi$ is shown. Moreover, the dimension
counterpart $\delta_{0,\orb}$ of $\chi_\orb$ is discussed, and we obtain
the relation of $\delta_{0,\orb}$ with the original free entropy dimension $\delta_0$
with applications to $\delta_0$ itself. 
\end{abstract}

\allowdisplaybreaks
\section*{Introduction}
\setcounter{equation}{0}

We propose a somewhat new approach to Voiculescu's theory of free entropy (see
e.g., \cite{V-Survey} for a survey), and introduce the orbital free entropy
$\chi_\orb(X_1,\dots,X_n)$. This quantity is an extension of the projection free
entropy $\chi_{\mathrm{proj}}(p_1,\dots,p_n)$ studied in \cite{HU3} following
Voiculescu's proposal in \cite{V6}. Our essential idea is to restrict microstates for
$(X_1,\dots,X_n)$ to unitary-orbital ones, that is, to use only the unitary parts of
microstates with disregarding the diagonal parts under their diagonalization. We prove
the exact relation (Theorem \ref{T-2.6})
\begin{equation}\label{F-0.1}
\chi(X_1,\dots,X_n)=\chi_\orb(X_1,\dots,X_n)+\sum_{i=1}^n\chi(X_i)
\end{equation}
between $\chi_\orb$ and the usual free entropy $\chi$ as naturally
expected from the definition. The particular formula
$$
-\chi_\orb(X,Y)=-\chi(X,Y)+\chi(X)+\chi(Y)
$$
resembles the expression $I(X;Y)=-H(X,Y)+H(X)+H(Y)$ of the classical mutual
information in terms of the Boltzmann-Gibbs entropy $H(\cdot)$. It should be
emphasized here that this expression of $I(X;Y)$ motivated Voiculescu to develop his
liberation theory and introduce the mutual free information $i^*$ in \cite{V6}.
In this way, we may regard $-\chi_\orb(X,Y)$ as a kind of free analog of the mutual
information and also as one possible  microstate version of the mutual free information
$i^*(W^*(X);W^*(Y))$. Among other properties, we prove (Theorem \ref{T-3.1}) that
$\chi_\orb(X_1,\dots,X_n)=0$ if and only if $X_1,\dots,X_n$ are freely independent
without any extra assumption. Together with the relation \eqref{F-0.1} this directly
implies the characterization \cite{V2,V4} of freeness by the additivity of $\chi$.
The proof of the theorem is based on a certain transportation cost inequality as in
the projection case in \cite{HU3}. An advantage of our orbital approach is that one of
equivalent definitions of $\chi_\orb$ as well as its all properties is valid even for
non-commutative random multi-variables $\bX_1,\dots,\bX_n$ each of which generates a
hyperfinite von Neumann algebra $W^*(\bX_i)$. But the present definition obeys the
essential restriction of hypefiniteness due to Jung's result \cite{Ju4} (or Lemma
\ref{L-1.2} below).

Furthermore, we study the dimension counterpart $\delta_{0,\orb}$ of $\chi_\orb$ for
hyperfinite random multi-variables $\bX_1,\dots,\bX_n$. It is defined similarly to the
modified free entropy dimension $\delta_0$ with replacing the semicircular deformation
by the liberation process (\cite{V6}). The $\delta_{0,\orb}$ enjoys properties similar
to $\delta_0$; for example, $\delta_{0,\orb}(\bX_1,\dots,\bX_n)=0$ if
$\chi_\orb(\bX_1,\dots,\bX_n)>-\infty$ (in particular, this is the case if
$\bX_1,\dots,\bX_n$ are free). Moreover, we prove the covering/packing formula of
$\delta_{0,\orb}$ based on Jung's approach \cite{Ju1,Ju2} to $\delta_0$, and
furthermore we prove the following general formula: 
\begin{equation*}\label{F-0.2}
\delta_0(\bX_1\sqcup\cdots\sqcup\bX_n) = \delta_{0,\orb}(\bX_1,\dots,\bX_n)
+ \sum_{i=1}^n\delta_0(\bX_i).
\end{equation*}
The orbital theory developed in this paper has several applications to the original
free entropy dimension $\delta_0$ itself. Among others, the most important one is
the following lower semicontinuity result for $\delta_0$: Let $\bX_1,\dots,\bX_n$ be
hyperfinite random multi-variables, and assume, for each $1\le i \le n$, that we have
a sequence of hyperfinite random multi-variables $\bX_i^{(k)}$ converging to $\bX_i$
in moments. In this setup, we will see that, if $\bX_i^{(k)} \subset W^*(\bX_i)$ is
further assumed for every $k \in \bN$ and $1 \le i \le n$, then 
\begin{equation*} 
\liminf_{k\rightarrow\infty}\delta_0(\bX^{(k)}_1\sqcup\cdots\sqcup\bX^{(k)}_n) \ge 
\delta_0(\bX_1\sqcup\cdots\sqcup\bX_n). 
\end{equation*}
Note here that the lower semicontinuity of $\delta_0$ was shown by Voiculescu in the
single variable case, see \cite{V-Survey} for the history on the semicontinuity problem
of $\delta_0$ until that time. Also, a certain related result was obtained by
Jung \cite[Lemma 7.3]{Ju1} based on his result on the $\delta_0$ of hyperfinite
algebras. However, Shlyakhtenko \cite{Sh2} pointed out that $\delta_0$ is never
lower semicontinuous in general. Our result is probably the first affirmative
semicontinuity result for $\delta_0$ of non-commutative nature.

\medskip\noindent
{\it Acknowledgment.} This work was completed during the third-named
author's stay in the Fields institute in 2007. He would like to thank the institute
for hospitality and also acknowledge Professor George Elliott for inviting him to the
operator algebra program. 

\section{Preliminaries}
\setcounter{equation}{0}

\subsection{Notations} We will use the standard notations; $M_N(\mathbb{C})$ denotes
the $N\times N$ matrix algebra, $\mathrm{Tr}_N$ stands for the usual (non-normalized)
trace and $\mathrm{tr}_N :=N^{-1}\mathrm{Tr}_N$ for its normalization. The operator
norm of a bounded operator $a$ is denoted by $\Vert a \Vert_\infty$. For a (tracial)
state $\phi$ on a $C^*$-algebra $A$ we write
$\Vert a \Vert_{p,\phi} := \phi(|a|^p)^{1/p}$ for $a \in A$ and $1 \leq p <\infty$.
Let $M_N^{sa}$ be the set of all $N\times N$ self-adjoint matrices. Let
$\mathrm{U}(N)$ be the $N\times N$ unitary group and $\mathrm{T}(N)$ be its (standard)
maximal torus consisting of all diagonal matrices in $\mathrm{U}(N)$ (isomorphic to
the $N$-dimensional torus $\mathbb{T}^N$ in the obvious way). We also use the
$N\times N$ special unitary group $\mathrm{SU}(N)$ in some places. The canonical
quotient map from $\mathrm{U}(N)$ onto the left coset space
$\mathrm{U}(N)/\mathrm{T}(N)$ is denoted by $q_N^\mathrm{U}$. We denote by
$\gamma_{\mathrm{U}(N)}$ the Haar probability measures on $\mathrm{U}(N)$ and by
$\gamma_{\mathrm{U}(N)/\mathrm{T}(N)}$ the probability measure on
$\mathrm{U}(N)/\mathrm{T}(N)$ induced from $\gamma_{\U(N)}$, that is, 
$\gamma_{\mathrm{U}(N)/\mathrm{T}(N)}=\gamma_{\mathrm{U}(N)}\circ(q_N^\U)^{-1}$.
The Haar probability measure on $\mathrm{SU}(N)$ is similarly denoted by
$\gamma_{\mathrm{SU}(N)}$.

\subsection{Matricial microstate spaces}  The following measure-space isomorphism is
well-known: 
\begin{equation}\label{F-1.1}
\Phi_N : 
\big((\mathrm{U}(N)/\mathrm{T}(N))\times\mathbb{R}^N_\geq,
\gamma_{\mathrm{U}(N)/\mathrm{T}(N)}\otimes\mu_N\big) \cong
\big(M_N^{sa}, \Lambda_N\big), 
\end{equation} 
where  
\begin{itemize} 
\item $\mathbb{R}^N_\geq$ denotes all $(x_1,\dots,x_N) \in \mathbb{R}^N$ with
$x_1 \geq x_2 \geq \cdots \geq x_N$,  
\item $\displaystyle\mu_N := C_N \prod_{i<j}(x_i-x_j)^2\prod_{i=1}^ndx_i$
with $\displaystyle C_N :={(2\pi)^{N(N-1)/2}\over\prod_{k=1}^{N-1}k!}$,
\item $\Lambda_N$ is the usual Lebesgue measure on $M_N^{sa}\cong\bR^{N^2}$.
\end{itemize} 
The map $\Phi_N$ is given by the continuous surjection
$$
([U],D) \in (\mathrm{U}(N)/\mathrm{T}(N))\times\mathbb{R}^N_\geq \mapsto
UDU^* \in M_N^{sa}
$$
with identifying $(x_1,\dots,x_N) \in \mathbb{R}_\geq^N$ with the diagonal matrices
$D$ whose diagonal entries are $x_1,\dots,x_N$ in decreasing order from
the upper diagonal corner. In what follows, we will identify an element in
$\mathbb{R}^N$ with a diagonal matrix in this canonical way. The above description
of $M_N^{sa}$ will be a key of our approach. 

\subsection{Technical lemmas} Here we recall three well-known lemmas, which
will be main technical ingredients in our discussions. The first lemma is just a
reformulation of Voiculescu's lemma \cite[Lemma 4.3]{V2} (also \cite[Lemma 4.3.4]{HP}).

\begin{lemma}\label{L-1.1} Let $1\le p<\infty$, $R > 0$ and $\eps >0$ be all
arbitrary. Then there exist an $m \in \mathbb{N}$ and a $\delta > 0$ so that for
every $N \in \mathbb{N}$ and every pair of diagonal matrices
$D_1, D_2 \in \mathbb{R}_\geq^N$ with $\Vert D_2 \Vert_\infty \leq R$, the condition 
\begin{equation*} 
\left|\mathrm{tr}_N(D_1^k) - \mathrm{tr}_N(D_2^k)\right|
< \delta \quad \text{for all $1 \leq k \leq m$}
\end{equation*} 
ensures that $\Vert D_1 - D_2 \Vert_{p,\mathrm{tr}_N} < \eps$. 
\end{lemma}    

Among some generalizations of Voiculescu's lemma above, an ultimate result due to Jung
\cite{Ju4} is the following:

\begin{lemma}\label{L-1.2}{\rm(Jung \cite{Ju4})}\enspace
Let $M$ be a von Neumann algebra with a faithful normal tracial state $\tau$, and
assume that $M$ is embeddable into the ultraproduct $R^\omega$ of the hyperfinite
II$_1$ factor and has a finite number of self-adjoint generators, say
$X_1,\dots,X_n \in M$. Let $1\le p<\infty$ be given. The following properties are
equivalent:
\begin{itemize}
\item[(1)] $M$ is hyperfinite.
\item[(2)] Any two embeddings of $M$ into $R^\omega$ are unitarily equivalent in
$R^\omega$.
\item[(3)] For each $\eps > 0$ there exist an $m \in \mathbb{N}$ and a $\delta > 0$
so that for every $N \in \mathbb{N}$, if two $n$-tuples
$(A_1,\dots,A_n)$, $(B_1,\dots,B_n) \in (M_N^{sa})^n$ satisfy 
$$
\qquad\left|\mathrm{tr}_N(A_{i_1}\cdots A_{i_k})
- \tau(X_{i_1}\cdots X_{i_k})\right| < \delta, \quad  
\left|\mathrm{tr}_N(B_{i_1}\cdots B_{i_k})
- \tau(X_{i_1}\cdots X_{i_k})\right| < \delta  
$$
for all $1 \leq i_1,\dots,i_k \leq n$ and $1 \leq k \leq m$, then there is a single
unitary $U \in \mathrm{U}(N)$ such that
$\Vert UA_i U^* - B_i \Vert_{p,\mathrm{tr}_N} < \eps$ for $1\le i\le n$.
\end{itemize}
\end{lemma}

Remark that Jung dealt with only the $2$-norm $\Vert\cdot\Vert_{2,\mathrm{tr}_N}$
but his argument clearly works for any $p$-norm.  

\medskip
Let $(A,\phi)$ be a non-commutative probability space, and $(\Omega_i)_{i\in I}$ be a
family subsets of $A$. Let  $(A^{\star I},\phi^{\star I})$ be the free product
of copies of $(A,\phi)$ indexed by $I$, and denote by $\iota_i$ the canonical
embedding of $A$ onto the $i$th copy of $A$ in $A^{\star I}$. For each $\eps > 0$ and
$m \in \mathbb{N}$ we will say that $(\Omega_i)_{i\in I}$ are {\it $(m,\eps)$-free}
in $(A,\phi)$ if 
\begin{equation*} 
\left|\phi(a_1\cdots a_k)
- \phi^{\star I}(\iota_{i_1}(a_1)\cdots \iota_{i_k}(a_k))\right| < \eps 
\end{equation*} 
for all $a_j \in \Omega_{i_j}$, $i_j \in I$ with $1 \leq j \leq k$ and
$1 \leq k \leq m$. The next result due to Voiculescu \cite[Corollary 2.13]{V-IMRN} is
a key ingredient when dealing with the freely independent situation. 

\begin{lemma}\label{L-1.3}{\rm(Voiculescu \cite{V-IMRN})}\enspace
Let $R>0$, $\eps > 0$, $\theta > 0$ and $m \in \mathbb{N}$ be given. Then there exists
an $N_0 \in \mathbb{N}$ such that 
\begin{align*} 
\gamma_{\mathrm{U}(N)}^{\otimes n}\big(\big\{ (U_1,\dots,U_p) \in \mathrm{U}(N)^p:
&\ \{T_1^{(0)},\dots,T_{q_0}^{(0)}\},
\{U_1 T_1^{(1)}U_1^*,\dots,U_1 T_{q_1}^{(1)}U_1^*\}, \\ 
&\dots,\{U_p T_1^{(p)}U_p^*,\dots,U_p T_{q_p}^{(p)}U_p^*\}
\ \text{are $(m,\eps)$-free} \big\}\big) > 1-\theta
\end{align*} 
whenever $N \geq N_0$ and $T_j^{(i)} \in M_N(\mathbb{C})$ with
$\Vert T_j^{(i)}\Vert_\infty \leq R$ for $1 \leq j \leq q_i$, $1 \leq q_i \leq m$,
$0 \leq i \leq p$ and $1 \leq p \leq m$. 
\end{lemma}

\section{Orbital free entropy $\chi_{\mathrm{orb}}$ and its basic properties}
\setcounter{equation}{0}

Throughout this section, let $(M,\tau)$ be a tracial $W^*$-probability space and
$(X_1,\dots,X_n)$ be an $n$-tuple of self-adjoint random variables in $(M,\tau)$. We
will use the standard notations such as the microstate set
$\Gamma_R(X_1,\dots,X_n;N,m,\delta)$ appearing in the course of defining the
{\it microstate free entropy} $\chi(X_1,\dots,X_n)$ (see \cite{V2}). We define a
free entropy-like quantity as follows.

\begin{definition}\label{D-2.1}{\rm
For each $\delta > 0$, $m,N \in \mathbb{N}$, $R > 0$ and $1 \leq i \leq n$, we denote
by $\Delta_R(X_i;N,m,\delta)$ the set of all $N\times N$ diagonal matrices
$D \in \mathbb{R}^N_\geq$ satisfying $\Vert D \Vert_{\infty} \leq R$ and
$\big|\mathrm{tr}_N(D^k) - \tau(X_i^k)\big| < \delta$ for all $1 \leq k \leq m$,
and by $\Gamma_{\mathrm{orb},R}(X_1,\dots,X_n; N,m,\delta)$ the set of all $n$-tuples
$(U_1,\dots,U_n)$ of $N\times N$ unitary matrices such that there exists an $n$-tuple
$(D_1,\dots,D_n)$ in $\prod_{i=1}^n \Delta_R(X_i;N,m,\delta)$ satisfying
$$
\big|\mathrm{tr}_N(U_{i_1}D_{i_1}U_{i_1}^*\cdots U_{i_k}D_{i_k}U_{i_k}^*)
- \tau(X_{i_1}\cdots X_{i_k})\big| < \delta
$$ 
for all $1\leq i_1,\dots,i_k \leq n$ with $1 \leq k \leq m$. We define 
\begin{align*} 
\chi_{\mathrm{orb},R}(X_1,\dots,X_n)
:= \lim_{m\rightarrow\infty, \delta \searrow 0} \limsup_{N\rightarrow\infty}    
\frac{1}{N^2}\log\gamma_{\mathrm{U}(N)}^{\otimes n}
(\Gamma_{\mathrm{orb},R}(X_1,\dots,X_n;N,m,\delta)),
\end{align*}
and 
\begin{equation*}
\chi_{\mathrm{orb}}(X_1,\dots,X_n) := \sup_{R>0}\chi_{\mathrm{orb},R}(X_1,\dots,X_n).
\end{equation*}
}\end{definition}

The above definition of $\Gamma_{\mathrm{orb},R}(X_1,\dots,X_n;N,m,\delta)$
clearly contains a superfluous condition. In fact, it can be rephrased more simply as
the set of all $(U_1,\dots,U_n) \in \mathrm{U}(N)^n$ such that
\begin{equation}\label{F-2.1}
(U_1 D_1 U_1^*,\dots,U_n D_n U_n^*) \in \Gamma_R(X_1,\dots,X_n;N,m,\delta)
\end{equation}
for some diagonal matrices $D_1,\dots,D_n \in \mathbb{R}^N_\geq$. The map
$\Phi_N : ([U],D) \mapsto UDU^*$ in \eqref{F-1.1} gives rise to the continuous
surjection $\Phi_N^n : (\mathrm{U}(N)/\mathrm{T}(N))^n \times (\mathbb{R}^N_\geq)^n
\rightarrow (M_N^{sa})^n$, which provides a measure-space isomorphism between those
measure spaces. Denote by $\mathrm{pr}_N^{\mathrm{U}}$ the projection map from
$(\mathrm{U}(N)/\mathrm{T}(N))^n \times(\mathbb{R}_\geq^N)^n$ onto the first $n$
factors $(\mathrm{U}(N)/\mathrm{T}(N))^n$. It is obvious that 
$$
(q_N^\mathrm{U})^n(\Gamma_{\mathrm{orb},R}(X_1,\dots,X_n; N,m,\delta))
= \mathrm{pr}_N^{\mathrm{U}}\bigl((\Phi_N^n)^{-1}
(\Gamma_R(X_1,\dots,X_n;N,m,\delta))\bigr), 
$$
so that $\Gamma_{\mathrm{orb},R}(X_1,\dots,X_n; N,m,\delta)$ is essentially the
projection of $\Gamma_R(X_1,\dots,X_n; N,m,\delta)$ to the unitary part via
matrix diagonalization.

In the following let us introduce two more definitions of $\chi_\orb$. The first one
is a slight modification of $\chi_\mathrm{orb}(X_1,\dots,X_n)$, where the (operator
norm) cut-off procedure is removed.   

\begin{definition}\label{D-2.2}{\rm
We define $\Gamma_\mathrm{orb}(X_1,\dots,X_n;N,m,\delta)$ to be the set of all
$(U_1,\dots,U_n)\in\U(N)^n$ satisfying \eqref{F-2.1} for some
$D_1,\dots,D_n\in\bR_\ge^N$ with the microstate set
$\Gamma(X_1,\dots,X_n;\allowbreak N,m,\delta)$ without cut-off by parameter $R$. Define
\begin{align*} 
\chi_{\mathrm{orb},\infty}(X_1,\dots,X_n)
:= \lim_{m\rightarrow\infty, \delta \searrow 0}
\limsup_{N\rightarrow\infty}&\frac{1}{N^2}\log
\gamma_{\mathrm{U}(N)}^{\otimes n}(\Gamma_{\mathrm{orb}}(X_1,\dots,X_n;N,m,\delta)). 
\end{align*} 
}\end{definition}

The next definition is a natural generalization of the projection free
entropy $\chi_{\mathrm{proj}}$ introduced and studied in \cite{HU3} following
Voiculescu's proposal in \cite[14.2]{V6}.   

\begin{definition}\label{D-2.3}{\rm
For each $1\le i\le n$ let us first choose and fix an $n$-tuple $(\xi_1,\dots,\xi_n)$
of sequences $\xi_i=\{\xi_i(N)\}$ of $\xi_i(N)\in M_N^{sa}$, $N\in\bN$, such that
$\xi_i(N)$ converges to $X_i$ in moments as $N \rightarrow \infty$ for $1\le i\le n$.
{\rm(}Of course such sequences always exist.{\rm)} We define
$\Gamma_{\mathrm{orb}}(X_1,\dots,X_n:\xi_1(N),\dots,\xi_n(N); N,m,\delta)$ to be the
set of all $(U_1,\dots,U_n)\in\U(N)^n$ such that 
$$
\big|\mathrm{tr}_N(U_{i_1}\xi_{i_1}(N)U_{i_1}^*\cdots
U_{i_k}\xi_{i_k}(N)U_{i_k}^*) - \tau(X_{i_1}\cdots X_{i_k})\big| < \delta
$$ 
for all $1\leq i_1,\dots,i_k \leq n$ with $1 \leq k \leq m$, that is, the set of all
$(U_1,\dots,U_n)\in\U(N)^n$ such that $(U_i\xi_i(N)U_i^*)_{i=1}^n$ is in
$\Gamma(X_1,\dots,X_n;N,m,\delta)$. Define
\begin{align*} 
&\chi_\mathrm{orb}(X_1,\dots,X_n:\xi_1,\dots,\xi_n) \\
&\quad:= \lim_{m\rightarrow\infty, \delta \searrow 0}
\limsup_{N\rightarrow\infty}\frac{1}{N^2}\log
\gamma_{\mathrm{U}(N)}^{\otimes n}(
\Gamma_{\mathrm{orb}}(X_1,\dots,X_n:\xi_1(N),\dots,\xi_n(N);N,m,\delta)).
\end{align*}
}\end{definition}

The next lemma asserts that all the three definitions in Definitions
\ref{D-2.1}--\ref{D-2.3} are equivalent. Thus, all the quantities will be denoted by
the same symbol $\chi_\mathrm{orb}$, and we call $\chi_\orb(X_1,\dots,X_n)$ the
{\it orbital free entropy} of $(X_1,\dots,X_n)$ since the definition is based on
``unitary-orbital microstates."

\begin{lemma}\label{L-2.4}
For any choice of $R\ge\max_{1\le i\le n}\|X_i\|_\infty$ and for any choice of an
approximating $n$-tuple $(\xi_1,\dots,\xi_n)$ one has  
\begin{align*}
\chi_{\mathrm{orb},\infty}(X_1,\dots,X_n)
&=\chi_{\mathrm{orb}}(X_1,\dots,X_n) \\
&=\chi_{\mathrm{orb},R}(X_1,\dots,X_n) \\
&=\chi_{\mathrm{orb}}(X_1,\dots,X_n:\xi_1,\dots,\xi_n).
\end{align*}
\end{lemma}

\begin{proof}
First, due to the invariance of $\gamma_{\U(N)}$ under unitary conjugation, we may 
and do assume that $(\xi_1,\dots,\xi_n)$ is  an $n$-tuple of sequences
$\xi_i=\{D_i(N)\}$  of diagonal matrices in $\mathbb{R}^N_\ge$.  Then it is obvious
that $\Gamma_{\mathrm{orb}}(X_1,\dots,X_n:D_1(N),\dots,D_n(N);N,m,\delta) \subset
\Gamma_{\mathrm{orb}}(X_1,\dots,X_n;N,m,\delta)$, which implies
$\chi_\orb(X_1,\dots,X_n:\xi_1,\dots,\xi_n)\le\chi_{\orb,\infty}(X_1,\dots,X_n)$.
Moreover, one can choose $D_i(N)$ so that $\|D_i(N)\|_\infty\le\|X_i\|_\infty$ for
all $N\in\bN$. In this case, whenever $R\ge R_0:=\max_{1\le i\le n}\|X_i\|_\infty$,
one has
$$
\Gamma_{\mathrm{orb}}(X_1,\dots,X_n:D_1(N),\dots,D_n(N);N,m,\delta)
\subset\Gamma_{\mathrm{orb},R}(X_1,\dots,X_n;N,m,\delta)
$$
and hence
$\chi_\orb(X_1,\dots,X_n:\xi_1,\dots,\xi_n) \le \chi_{\orb,R}(X_1,\dots,X_n)$
($\le \chi_{\orb,\infty}(X_1,\dots,X_n)$). Thus, it suffices to prove that for any
approximating sequences $\xi_i=\{D_i(N)\}$ and for every
$m \in \mathbb{N}$ and $\delta > 0$, there are an $m'\in\mathbb{N}$, a $\delta'>0$
and an $N_0\in\bN$ so that
\begin{equation}\label{F-2.2} 
\Gamma_{\mathrm{orb}}(X_1,\dots,X_n;N,m',\delta') \subset
\Gamma_{\mathrm{orb}}(X_1,\dots,X_n:D_1(N),\dots,D_n(N);N,m,\delta) 
\end{equation} 
for all $N \ge N_0$. Choose a $\rho\in(0,1)$ with $m(R_0+1)^{m-1}\rho<\delta/2$.
By Lemma \ref{L-1.1} one can find an $m'\in\bN$ with $m'\ge 2m$, a $\delta'>0$ with
$\delta'\le\min\{1,\delta/2\}$ and an $N_0\in\bN$ such that for every $1\le i\le n$
and every $D_i\in\bR_\ge^N$ with $N\ge N_0$, if
$|\mathrm{tr}_N(D_i^k) - \tau(X_i^k)|<\delta'$ for all $1\le k\le m'$, then
$\|D_i-D_i(N)\|_{m,\tr_N}<\rho$. Suppose $N\ge N_0$ and $(U_1,\dots,U_n)$ is in
the left-hand side of \eqref{F-2.2} so that
$(U_iD_iU_i^*)_{i=1}^n\in\Gamma(X_1,\dots,X_n;N,m',\delta')$ for some
$D_1,\dots,D_n\in\bR_\ge^N$. Since $\|D_i-D_i(N)\|_{m,\tr_N}<\rho$ and
\begin{align*}
\|D_i\|_{m,\tr_N}&\le\tr_N(D_i^{2m})^{1/2m}<(\tau(X_i^{2m})+\delta')^{1/2m} \\
&\le(R_0^{2m}+1)^{1/2m}\le R_0+1,
\end{align*}
we get
\begin{align*} 
&\left|\mathrm{tr}_N(U_{i_1}D_{i_1}(N)U_{i_1}^*\cdots U_{i_k}D_{i_k}(N)U_{i_k}^*)
- \tau(X_{i_1}\cdots X_{i_k})\right| \\
&\quad\le
\left|\mathrm{tr}_N(U_{i_1}D_{i_1}(N)U_{i_1}^*\cdots U_{i_k}D_{i_k}(N)U_{i_k}^*)
- \mathrm{tr}_N(U_{i_1}D_{i_1}U_{i_1}^*\cdots U_{i_k}D_{i_k}U_{i_k}^*)\right| \\
&\qquad
+\left|\mathrm{tr}_N(U_{i_1}D_{i_1}U_{i_1}^*\cdots U_{i_k}D_{i_k}U_{i_k}^*)
- \tau(X_{i_1}\cdots X_{i_k})\right| \\
&\quad\le m(R_0+1)^{m-1}\rho+\delta'<\delta
\end{align*}   
for all $1\leq i_1,\dots,i_k \leq n$ with $1 \leq k \leq m$. The above latter
inequality is seen by the H\"older inequality. This implies that $(U_1,\dots,U_n)$ is
in the right-hand side of \eqref{F-2.2}.
\end{proof}

Some basic properties of $\chi_\mathrm{orb}$ are summarized in the next proposition.
The properties (1)--(3) are obvious, and (4) is seen by using Definition \ref{D-2.2}
due to Lemma \ref{L-2.4}. 

\begin{prop}\label{P-2.5}
$\chi_\mathrm{orb}$ enjoys the following properties: 
\begin{itemize}
\item[(1)] $\chi_\mathrm{orb}(X) = 0$ for any single random variable. 
\item[(2)] $\chi_\mathrm{orb}(X_1,\dots,X_n) \leq 0$.
\item[(3)] $\chi_\mathrm{orb}(X_1,\dots,X_n) \leq
\chi_\mathrm{orb}(X_1,\dots,X_k) + \chi_\mathrm{orb}(X_{k+1},\dots,X_n)$ for every
$1\le k<n$.
\item[(4)] If $(X_1^{(k)},\dots,X_n^{(k)})$, $k\in\bN$, are $n$-tuples of
self-adjoint random variables converging to $(X_1,\dots,X_n)$ in the distribution
sense as $k \rightarrow \infty$, then 
\begin{equation*} 
\chi_\mathrm{orb}(X_1,\dots,X_n) \geq 
\limsup_{k\rightarrow\infty} \chi_\mathrm{orb}(X_1^{(k)},\dots,X_n^{(k)}). 
\end{equation*}  
\end{itemize} 
\end{prop}

The following exact relation between $\chi_\orb$ and the usual $\chi$ is
the main result of this section.

\begin{thm}\label{T-2.6}
$$
\chi(X_1,\dots,X_n)=\chi_{\mathrm{orb}}(X_1,\dots,X_n) + \sum_{i=1}^n \chi(X_i).
$$  
\end{thm} 

\begin{proof}
Let $R\ge\max_{1\le i\le n}\|X_i\|_\infty$. Since
\begin{align*}
&(\Phi_N^n)^{-1}(\Gamma_R(X_1,\dots,X_n;N,m,\delta)) \\
&\qquad \subset
(q_N^\mathrm{U})^n(\Gamma_{\mathrm{orb},R}(X_1,\dots,X_n;N,m,\delta))
\times \prod_{i=1}^n \Delta_R(X_i;N,m,\delta) 
\end{align*}
and $\Gamma_{\mathrm{orb},R}(X_1,\dots,X_n;N,m,\delta)$ is invariant under the right
multiplication by elements of $\mathrm{T}(N)^n$, we get by \eqref{F-1.1}
\allowdisplaybreaks{    
\begin{align*} 
&\log\Lambda_N^{\otimes n}(\Gamma_R(X_1,\dots,X_n;N,m,\delta)) \\
&\quad\leq 
\log\gamma_{\mathrm{U}(N)}^{\otimes n}
(\Gamma_{\mathrm{orb},R}(X_1,\dots,X_n;N,m,\delta))
+ \sum_{i=1}^n \log\mu_N(\Delta_R(X_i;N,m,\delta)) \\
&\quad= 
\log\gamma_{\mathrm{U}(N)}^{\otimes n}
(\Gamma_{\mathrm{orb},R}(X_1,\dots,X_n;N,m,\delta))
+ \sum_{i=1}^n \log\Lambda_N(\Gamma_R(X_i;N,m,\delta)).
\end{align*}  
}By Lemma \ref{L-2.4} this immediately implies the inequality $\le$ for the required
equality.

For the reverse inequality we show that for each $m \in \mathbb{N}$ and $\delta>0$
there are an $m' \in \mathbb{N}$, a $\delta' >0$ and an $N_0 \in\mathbb{N}$ so that 
\begin{align}\label{F-2.3} 
&(q_N^\mathrm{U})^n(\Gamma_{\mathrm{orb},R}(X_1,\dots,X_n;N,m',\delta'))
\times \prod_{i=1}^n \Delta_R(X_i;N,m',\delta') \notag \\
&\qquad\subset
(\Phi_N^n)^{-1}(\Gamma_R(X_1,\dots,X_n;N,m,\delta)) 
\end{align} 
for all $N \geq N_0$. By Lemma \ref{L-1.1} one can find an $m'\in\bN$, a
$\delta'\in(0,\delta/2)$ and an $N_0\in\bN$ such that for every $1\le i\le n$ and
every $N\ge N_0$, if $D_i,D'_i\in\Delta_R(X_i;N,m',\delta')$ then
\begin{equation*}
\Vert D_i - D'_i\Vert_{1,\mathrm{tr}_N} < \frac{\delta}{2m(R+1)^{m-1}}.
\end{equation*}
Now suppose $N\ge N_0$ and $([U_1],\dots,[U_n],D_1,\dots,D_n)$ is in the left-hand
side of \eqref{F-2.3}. Then we have
$(U_iD'_iU_i^*)_{i=1}^n \in \Gamma_R(X_1,\dots,X_n;N,m',\delta')$ for some
$(D'_i)_{i=1}^n \in \prod_{i=1}^n \Delta_R(X_i;\allowbreak N,m',\delta')$. Since
\begin{align*} 
&\left|\mathrm{tr}_N(U_{i_1}D_{i_1}U_{i_1}^*\cdots U_{i_k}D_{i_k}U_{i_k}^*)
- \tau(X_{i_1}\cdots X_{i_k})\right| \\
&\quad\leq 
\left|\mathrm{tr}_N(U_{i_1}D_{i_1}U_{i_1}^*\cdots U_{i_k}D_{i_k}U_{i_k}^*)
- \mathrm{tr}_N(U_{i_1}D'_{i_1}U_{i_1}^*\cdots U_{i_k}D'_{i_k}U_{i_k}^*)\right| \\ 
&\quad\quad+ 
\left|\mathrm{tr}_N(U_{i_1}D'_{i_1}U_{i_1}^*\cdots U_{i_k}D'_{i_k}U_{i_k}^*)
- \tau(X_{i_1}\cdots X_{i_k})\right| \\
&\quad\le
m(R+1)^{m-1} \max_{1\le i\le n}\Vert D_i - D'_i\Vert_{1,\mathrm{tr}_N} +\delta'<\delta
\end{align*}   
for all $1\leq i_1,\dots,i_k \leq n$ and $1 \leq k \leq m$, it follows that
$(U_iD_iU_i^*)_{i=1}^n \in \Gamma_R(X_1,\dots,X_n;\allowbreak N,m,\delta)$, proving
\eqref{F-2.3}. By Lemma \ref{L-2.4} we thus obtain
\allowdisplaybreaks{
\begin{align*} 
&\chi_{\mathrm{orb}}(X_1,\dots,X_n) + \sum_{i=1}^n \chi(X_i) \\
&\quad\leq \limsup_{N\rightarrow\infty}\frac{1}{N^2}\log
\gamma_{\mathrm{U}(N)}^{\otimes n}(
\Gamma_{\mathrm{orb},R}(X_1,\dots,X_n;N,m',\delta')) \\
&\quad\quad+
\sum_{i=1}^n \biggl(\lim_{N\rightarrow\infty}\frac{1}{N^2}\log
\Lambda_N\big(\Gamma_R(X_i;N,m',\delta')\big) + \frac{1}{2}\log N\biggr) \\
&\quad\leq \limsup_{N\rightarrow\infty}\biggl(
\frac{1}{N^2}\log\Lambda_N^{\otimes n}(\Gamma_R(X_1,\dots,X_n;N,m,\delta)
+ \frac{n}{2}\log N\biggr)
\end{align*}
}for every $m \in \mathbb{N}$ and $\delta>0$. This implies inequality $\ge$ for the
desired equality. (A point in the above proof is that $\limsup$ can be replaced by
$\lim$ in the definition of $\chi(X)$ in the single variable case, see
\cite[5.6.2]{HP}.)    
\end{proof}

Theorem \ref{T-2.6} in particular gives
\begin{equation*}
-\chi_{\mathrm{orb}}(X,Y) = -\chi(X,Y) + \chi(X) + \chi(Y)
\end{equation*}
for two (non-commutative) self-adjoint random variables $X,Y$ in $(M,\tau)$ with
$\chi(X),\chi(Y)>-\infty$. The above expression suggests that $-\chi_\orb$
is a kind of free probability counterpart of the so-called mutual information $I(X;Y)$
for two real random variables $X,Y$. In fact, recall the expression
\begin{equation}\label{F-2.4}
I(X;Y)=-H(X,Y)+H(X)+H(Y)
\end{equation}
in terms of the Boltzmann-Gibbs entropy $H(\cdot)$, which holds as long as $H(X)$ and
$H(Y)$ are finite. The following remark is another justification
for the analogy between $-\chi_\orb$ and the classical mutual information.

\begin{remark}\label{R-2.7}{\rm
Let $\gamma_{\mathfrak{S}_N}$ denote the uniform probability measure on the symmetric
group $\mathfrak{S}_N$. Let $X_1,\dots,X_n$ be bounded real random variables on a
classical probability space. For $N,m\in\bN$ and $\delta>0$ define
$\Delta(X_1,\dots,X_n;N,m,\delta)$ to be the set of all $n$-tuples $(x_1,\dots,x_n)$
of vectors $x_i=(x_{i1},\dots,x_{iN})$ in $\bR^N$ such that
$$
\left|{1\over N}\sum_{j=1}^Nx_{i_1j}\cdots x_{i_kj}
-\mathbb{E}(X_{i_1}\cdots X_{i_k})\right|<\delta
$$
for all $1\le i_1,\dots,i_k\le n$ and $1\le k\le m$, where $\mathbb{E}(\cdot)$ is the
expectation. Moreover, define $\Delta_\sym(X_1,\dots,X_n;N,m,\delta)$ to be the set
of all $n$-tuples $(\sigma_1,\dots,\sigma_n)$ of $\sigma_i\in\mathfrak{S}_N$ such that
$(\sigma_1(x_1),\dots,\sigma_n(x_n))\in\Delta(X_1,\dots,X_n;N,m,\delta)$ for some
$x_1,\dots,x_n\in\bR_\ge^N$, where
$\sigma_i(x_i):=(x_{i\sigma_i(1)},\dots,x_{i\sigma(N)})$. We then define
$$
H_\sym(X_1,\dots,X_n):=\lim_{m\rightarrow\infty,\delta\searrow 0}
\limsup_{N\to\infty}{1\over N}\log\gamma_{S_N}^{\otimes n}
\bigl(\Delta_\sym(X_1,\dots,X_n;N,m,\delta)\bigr).
$$
One can show that
$$
H(X_1,\dots,X_n)=H_\sym(X_1,\dots,X_n)+\sum_{i=1}^nH(X_i).
$$
In particular, when $X$ and $Y$ are real bounded random variables with
$H(X),H(Y)>-\infty$, we have $I(X;Y)=-H_\sym(X,Y)$. In this way, the
``classical analog" of $-\chi_\orb(X,Y)$ provides a new definition (a kind of
``discretization") of the classical mutual information $I(X;Y)$. More on
this idea are examined in \cite{HP2}.
}\end{remark}

It seems that the expression \eqref{F-2.4} was one of the motivations of Voiculescu
to introduce the mutual free information $i^*(A_1;\dots;A_n)$ for subalgebras
$A_1,\dots,A_n$ in \cite{V6} (in particular, see Introduction there). For any $n$-tuple
of projections $(p_1,\dots,p_n)$ in a $W^*$-probability space, from the definition in
\cite{HU3} and Lemma \ref{L-2.4} we notice that
\begin{equation*} 
\chi_\mathrm{proj}(p_1,\dots,p_n) = \chi_\mathrm{orb}(p_1,\dots,p_n).  
\end{equation*} 
In \cite{HU2} we conjectured that $-\chi_\mathrm{proj}(p,q)$ coincides with the mutual
free information $i^*(\mathbb{C}p+\mathbb{C}(1-p);\mathbb{C}q+\mathbb{C}(1-q))$
for two projections $p,q$, and gave a heuristic computation supporting it. It would be
further conjectured that $-\chi_\mathrm{orb}(X,Y) = i^*(W^*(X) ; W^*(Y))$ holds for
any $X, Y$; however this is out of scope of this paper. Here note that this is true
when $X,Y$ are freely independent (see Proposition \ref{P-2.9} below). From the above
point of view we are tempted to write 
$i(W^*(X_1);\dots;W^*(X_n)):= -\chi_{\mathrm{orb}}(X_1,\dots,X_n)$ and use the term
``microstate mutual free information." However we leave the symbol $i$ to further
progress on the subject.

In view of the analogy between $-\chi_{\mathrm{orb}}$ and $I(X;Y)$ the following
proposition is strongly expected.     

\begin{prop}\label{P-2.8}
$\chi_{\rm orb}(X_1,\dots,X_n)$ depends only upon $W^*(X_1),\dots,W^*(X_n)$, where
$W^*(X_i)$ means the von Neumann subalgebra of $M$ generated by $X_i$ {\rm(}and the
unit $\1 \in M${\rm)}.
\end{prop}

\begin{proof}
Let $(X'_1,\dots,X'_n)$ be another $n$-tuple of self-adjoints in $M$ with
$W^*(X'_i)=W^*(X_i)$ for $1\le i\le n$. By symmetry and Lemma \ref{L-2.4} it
suffices to prove that for each $m\in\bN$ and $\delta>0$ there are an $m'\in\bN$ and
a $\delta'>0$ such that
\begin{align}\label{F-2.5}
&\gamma_{\U(N)}^{\otimes n}(
\Gamma_\orb(X_1,\dots,X_n:\xi_1(N),\dots,\xi_n(N);N,m',\delta')) \nonumber\\
&\qquad\qquad\qquad
\le\gamma_{\U(N)}^{\otimes n}(\Gamma_\orb(X'_1,\dots,X'_n;N,m,\delta))
\end{align}
for all $N\in\bN$, where $\{\xi_i(N)\}$ is an approximating sequence for $X_i$
as in Definition \ref{D-2.3} for $1\le i\le n$. Let
$R:=\max_{1\le i\le n}\|X'_i\|_\infty$. For each $1\le i\le n$, since
$X'_i\in W^*(X_i)$, one can choose, by the Kaplansky theorem, a real polynomial
$P_i(t)$ such that $\|X'_i-P_i(X_i)\|_{1,\tau}<\allowbreak\delta/2m(R+1)^{m-1}$ and
$\|P_i(X_i)\|_\infty\le\|X'_i\|_\infty$. For each $m\in\bN$ and $\delta>0$ one can
choose an $m'\in\bN$ and a $\delta'>0$ (depending on $P_1,\dots,P_n$ as well)
such that $(A_i)_{i=1}^n\in\Gamma(X_1,\dots,X_n;N,m',\delta')$ implies
$(P_i(A_i))_{i=1}^n\in\Gamma(P_1(X_1),\dots,P_n(X_n);N,m,\delta/2)$ for every
$N\in\bN$. If
$(U_1,\dots,U_n)\in\Gamma_\orb(X_1,\dots,X_n:\xi_1(N),\dots,\xi_n(N);N,m',\delta')$,
then we get
\begin{align*}
&|\tr_N(U_{i_1}P_{i_1}(\xi_{i_1}(N))U_{i_1}^*\cdots
U_{i_k}P_{i_k}(\xi_{i_k}(N))U_{i_k}^*)
-\tau(X'_{i_1}\cdots X'_{i_k})| \\
&\quad\le|\tr_N(P_{i_1}(U_{i_1}\xi_{i_1}(N)U_{i_1}^*)\cdots
P_{i_k}(U_{i_k}\xi_{i_k}(N)U_{i_k}^*))
-\tau(P_{i_1}(X_{i_1})\cdots P_{i_k}(X_{i_k}))| \\
&\quad\quad+|\tau(P_{i_1}(X_{i_1})\cdots P_{i_k}(X_{i_k}))
-\tau(X'_{i_1}\cdots X'_{i_k})| \\
&\quad\le{\delta\over2}+m(R+1)^{m-1}\max_{1\le i\le n}\|X'_i-P_i(X_i)\|_{1,\tau}
<\delta
\end{align*}
for all $1\le i_1,\dots,i_k\le n$ and $1\le k\le m$. Now, for each $N\in\bN$ and
$1\le i\le n$ write $P_i(\xi_i(N))=V_i(N)D_i(N)V_i(N)^*$ with $D_i(N)\in\bR_\ge^N$
and $V_i(N)\in\U(N)$. Then we have
\begin{align*}
&\Gamma_\orb(X_1,\dots,X_n:\xi_1(N),\dots,\xi_n(N);N,m',\delta'))
\cdot(V_1(N),\dots,V_n(N)) \\
&\qquad\qquad\qquad\qquad
\subset\Gamma_\orb(X'_1,\dots,X'_n;N,m,\delta),
\end{align*}
and \eqref{F-2.5} follows from the right invariance of the Haar measure
$\gamma_{\U(N)}$.
\end{proof} 

If $\chi(X_i) > -\infty$ for all $1 \leq i \leq n$ and $X_1,\dots,X_n$ are freely
independent, then the additivity theorem \cite[Proposition 5.4]{V2} and Theorem
\ref{T-2.6} show that $\chi_{\mathrm{orb}}(X_1,\dots,X_n)=0$ (or the additivity of
$\chi_\orb$ in view of Proposition \ref{P-2.5}\,(1)). The next proposition shows that
this is still true even when the finiteness assumption of the $\chi(X_i)$'s is
dropped.

\begin{prop}\label{P-2.9}
If $X_1$ is freely independent of $X_2,\dots,X_n$, then
$$
\chi_\orb(X_1,X_2,\dots,X_n)=\chi_\orb(X_2,\dots,X_n).
$$
Consequently, $\chi_\orb(X_1,\dots,X_n)=0$ if $X_1,\dots,X_n$ are freely independent.
\end{prop}

\begin{proof}
The proof is based on the method due to Voiculescu \cite{V-IMRN} (or Lemma
\ref{L-1.3}) while it is easier than that for the additivity of $\chi$. 
By (1) and (3) of Proposition \ref{P-2.5} we may prove that
\begin{equation}\label{F-2.6}
\chi_\orb(X_1,X_2,\dots,X_n)\ge\chi_\orb(X_2,\dots,X_n).
\end{equation}
under the assumption $\chi_\orb(X_2,\dots,X_n)>-\infty$. Choose an approximating
sequence $\{\xi_i(N)\}$ for $X_i$ with $\|\xi_i(N)\|_\infty\le\|X_i\|_\infty$ for
$1\le i\le n$. For $N,m\in\bN$ and $\delta,\delta'>0$ set
\begin{align*}
\Psi(N,m,\delta)&:=\Gamma_\orb(X_1,\dots,X_n:\xi_1(N),\dots,\xi_n(N);N,m,\delta), \\
\Phi(N,m,\delta')&:=\Gamma_\orb(X_2,\dots,X_n:\xi_2(N),\dots,\xi_n(N);N,m,\delta'),
\end{align*}
and moreover
\begin{align*}
\Omega(N,m,\delta')&:=\bigl\{(U_1,U_2,\dots,U_n)\in\U(N)^n:
\{U_1\xi_1(N)U_1^*\}\ \mbox{and}\\
&\hskip2cm\{U_2\xi_2(N)U_2^*,\dots,U_n\xi_n(N)U_n^*\}
\ \mbox{are $(m,\delta')$-free}\bigr\}.
\end{align*}
For every $m\in\bN$ and $\delta>0$ one can find a $\delta'>0$ such that if
$U_1\in\Gamma_\orb(X_1:\xi_1(N);N,m,\delta')$ and
$(U_2,\dots,U_n)\in\Phi(N,m,\delta')$ and if $(U_1,U_2,\dots,U_n)$ is in
$\Omega(N,m,\delta')$, then $(U_1,\dots,U_n)$ is in $\Psi(N,m,\delta)$. Note that
$\Gamma_\orb(X_1:\xi_1(N);N,m,\delta')$ is the whole $\U(N)$ for sufficiently large
$N$. Hence, by Lemma 1.3 there is an $N_0\in\bN$ such that
$\Gamma_\orb(X_1:\xi_1(N);N,m,\delta')=\U(N)$ and
$$
\gamma_{\U(N)}\bigl(\{U_1\in\U(N):(U_1,U_2,\dots,U_n)
\in\Omega(N,m,\delta')\}\bigr)\ge{1\over2}
$$
for all $N\ge N_0$ and every $(U_2,\dots,U_n)\in\U(N)^{n-1}$. From the assumption
$\chi_\orb(X_2,\dots,X_n)\allowbreak>-\infty$, we may assume that
$\gamma_{\U(N)}^{\otimes n-1}(\Phi(N,m,\delta'))>0$ for all $N\ge N_0$. Hence, with
the measure
$$
\mu_N:={1\over\gamma_{\U(N)}^{\otimes n-1}(\Phi(N,m,\delta'))}
\,\gamma_{\U(N)}^{\otimes n-1}\big|_{\Phi(N,m,\delta')},
$$
we get for every $N\ge N_0$
\begin{align*}
&{\gamma_{\U(N)}^{\otimes n}(\Psi(N,m,\delta))
\over\gamma_{\U(N)}^{\otimes n-1}(\Phi(N,m,\delta'))} \\
&\quad\ge\int_{\Phi(N,m,\delta')}\Biggl(\int_{\U(N)}
\1_{\Omega(N,m,\delta')}(U_1,U_2,\dots,U_n)
\,d\gamma_{\U(N)}(U_1)\Biggr)\,d\mu_N(U_2,\dots,U_n)\ge{1\over2}.
\end{align*}
Therefore,
$$
\limsup_{N\to\infty}{1\over N^2}\log
\gamma_{\U(N)}^{\otimes n}(\Psi(N,m,\delta))
\ge\limsup_{N\to\infty}{1\over N^2}\log
\gamma_{\U(N)}^{\otimes n-1}(\Phi(N,m,\delta')),
$$
which implies \eqref{F-2.6} thanks to Lemma \ref{L-2.4}.
\end{proof}

\section{Characterization of freeness by $\chi_\orb=0$}
\setcounter{equation}{0}

Let $(X_1,\dots,X_n)$ be an $n$-tuple of self-adjoint random variables as in the
preceding section. This section is devoted to proving the converse implication of
the second assertion of Proposition \ref{P-2.9}; consequently we have the
following:

\begin{thm}\label{T-3.1}
$\chi_\orb(X_1,\dots,X_n)=0$ if and only if $X_1,\dots,X_n$ are freely independent.
\end{thm}

To prove the theorem, we will provide a certain transportation cost inequality
similarly to the projection case in \cite[\S5]{HU3}. In what follows we adopt the
description of $\chi_\orb$ as
$\chi_{\mathrm{orb}}(X_1,\dots,X_n:\xi_1,\dots,\xi_n)$ due to Lemma \ref{L-2.4}.
For $1\le i\le n$ let us choose and fix a sequence $\{\xi_i(N)\}$ of
$\xi_i(N)\in M_N^{sa}$ such that $\|\xi_i(N)\|_\infty\le\|X_i\|_\infty$ and
$\xi_i(N)\to X_i$ in moments as $N\to\infty$.

With $R :=\max_{1\le i\le n}\Vert X_i\Vert_\infty$, $C[-R,R]$ is the $C^*$-algebra
of continuous functions on $[-R,R]$. Let $\mathcal{A}_R := C[-R,R]^{\star n}$ be
the universal free product $C^*$-algebra of $n$-copies of $C[-R,R]$ with canonical
self-adjoint free generators $Z_1,\dots,Z_n$, i.e., $Z_i(t)=t$ in the $i$th $C[-R,R]$.
We denote by $TS(\mathcal{A}_R)$ the set of all tracial states on $\mathcal{A}_R$ and
by $\mathcal{P}(\mathrm{SU}(N)^n)$ the set of all probability measures on the $n$-fold
product $\mathrm{SU}(N)^n$. For each $\lambda \in \mathcal{P}(\mathrm{SU}(N)^n)$ we
associate a unique $\widehat{\lambda} \in TS(\mathcal{A}_R)$ as follows:
\begin{equation}\label{F-3.1}
\widehat{\lambda}(h) := \int_{\mathrm{SU}(N)^n}
\mathrm{tr}_N(h(U_1\xi_1(N)U_1^*,\dots,U_n\xi_n(N)U_n^*))\,d\lambda
\quad \text{for $h \in \mathcal{A}_R$},
\end{equation}
where $h(U_1\xi_1(N)U_1^*,\dots,U_n\xi_n(N) U_n^*)$ is the image of
$h\in\mathcal{A}_R$ by the $*$-homomorphism from $\mathcal{A}_R$ to $M_N(\mathbb{C})$
sending each $Z_i$ to $U_i \xi_i(N) U_i^*$. Similarly,
$\tau_{(X_1,\dots,X_n)} \in TS(\mathcal{A}_R)$ is defined by
\begin{equation*}
\tau_{(X_1,\dots,X_n)}(h) := \tau(h(X_1,\dots,X_n))
\quad \text{for $h \in \mathcal{A}_R$}.
\end{equation*}
From the trivial fact that the image of $\mathrm{SU}(N)$ by the quotient map
$q^\mathrm{U}_N$ is exactly $\mathrm{U}(N)/\mathrm{T}(N)$, it is clear that no
difference occurs when $\mathrm{SU}(N)$ is used in place of $\mathrm{U}(N)$ in the
definition of $\chi_\mathrm{orb}$ (Definitions \ref{D-2.1}--\ref{D-2.3}). Letting 
\begin{equation*} 
\Gamma(N,m,\delta)
:= \Gamma_\mathrm{orb}(X_1,\dots,X_n:\xi_1(N),\dots,\xi_n(N);N,m,\delta)
\cap \mathrm{SU}(N)^n,
\end{equation*} 
we thus have 
\begin{equation*} 
\chi_\mathrm{orb}(X_1,\dots,X_n)
= \lim_{m\rightarrow\infty, \delta \searrow 0}\limsup_{N\rightarrow\infty}
\frac{1}{N^2}\log\gamma_{\mathrm{SU}(N)}^{\otimes n}(\Gamma(N,m,\delta)).
\end{equation*} 
Now, we can choose a subsequence $N_1 < N_2 < \cdots$ in such a way that 
\begin{equation}\label{F-3.2}
\chi_\mathrm{orb}(X_1,\dots,X_n)
= \lim_{m\rightarrow\infty}\frac{1}{N_m^2}\log
\gamma_{\mathrm{SU}(N_m)}^{\otimes n}(\Gamma(N_m,m,1/m)),  
\end{equation}
and define 
\begin{equation*} 
\lambda_m := \frac{1}{\gamma_{\mathrm{SU}(N_m)}^{\otimes n}(\Gamma_m)}
\,\gamma_{\mathrm{SU}(N_m)}^{\otimes n}\big|_{\Gamma_m}
\in \mathcal{P}(\mathrm{SU}(N)^n)  
\end{equation*}
with $\Gamma_m := \Gamma(N_m,m,1/m)$. Then the next lemma can be proven in the same
way as in the proof of \cite[(2.5) in p.\,401]{HU1}.   
    
\begin{lemma}\label{L-3.2}
$\lim_{m\rightarrow\infty}\widehat{\lambda}_m = \tau_{(X_1,\dots,X_n)}$ in the
weak$^*$ topology. 
\end{lemma} 

The following is essentially a kind of reformulation of Voiculescu's asymptotic
freeness result \cite{V0,V-IMRN} (also \cite[\S4.3]{HP}) for unitary random
matrices (related to Lemma \ref{L-1.3}). A simple proof based on Lemma \ref{L-1.3} is
given for completeness.

\begin{lemma}\label{L-3.3}
$\lim_{N\rightarrow\infty}\widehat{\gamma_{\mathrm{SU}(N)}^{\otimes n}}
= \tau_{(X_1,\dots,X_n)}^\mathrm{free}$ in the weak$^*$ topology,
where\ $\tau_{(X_1,\dots,X_n)}^\mathrm{free} := \star_{i=1}^n\tau_{X_i}$ is the
free product of the states $\tau_{X_i}$ on $C[-R,R]$ induced from the distribution
measure of $X_i$.
\end{lemma}

\begin{proof}
For each $m\in\bN$ and $\delta,\theta>0$, Lemma \ref{L-1.3} implies that
$\gamma_{\U(N)}^{\otimes n}(\Omega(N,m,\delta))>1-\theta$ for all
sufficiently large $N$, where
$$
{\Omega}(N,m,\delta):=\{(U_1,\dots,U_n)\in\U(N)^n:
\{U_1\xi_1(N)U_1^*\},\dots,\{U_n\xi_n(N)U_n^*\}\ \mbox{are $(m,\delta)$-free}\}.
$$
For any $1\le i_1,\dots,i_k\le n$ with $1\le k\le m$, notice that
\begin{align*}
\widehat{\gamma_{\SU(N)}^{\otimes n}}(Z_{i_1}\cdots Z_{i_k})
&=\int_{\SU(N)^n}\tr_N(U_{i_1}\xi_{i_1}(N)U_{i_1}^*\cdots
U_{i_k}\xi_{i_k}(N)U_{i_k}^*)\,d\gamma_{\SU(N)}^{\otimes n} \\
&=\int_{\U(N)^n}\tr_N(U_{i_1}\xi_{i_1}(N)U_{i_1}^*\cdots
U_{i_k}\xi_{i_k}(N)U_{i_k}^*)\,d\gamma_{\U(N)}^{\otimes n}
\end{align*}
and
$$
\tr_N^{\star n}(\iota_{i_1}(\xi_{i_1}(N))\cdots\iota_{i_k}(\xi_{i_k}(N))
=\tr_N^{\star n}(\iota_{i_1}(U_{i_1}\xi_{i_1}(N)U_{i_1}^*)\cdots
\iota_{i_k}(U_{i_k}\xi_{i_k}(N)U_{i_k}^*))
$$
for all $(U_1,\dots,U_n)\in\U(N)^n$ (with the notations given before Lemma
\ref{L-1.3}). Hence one can immediately estimate
\begin{align*} 
&\left|\widehat{\gamma_{\SU(N)}^{\otimes n}}(Z_{i_1}\cdots Z_{i_k})
-\tr_N^{\star n}(\iota_{i_1}(\xi_{i_1}(N))\cdots\iota_{i_k}(\xi_{i_k}(N)))\right| \\
&\quad\le\Biggl(\int_{\Omega(N,m,\delta)}
+\int_{\U(N)^n\setminus\Omega(N,m,\delta)}\Biggr)
\big|\tr_N(U_{i_1}\xi_{i_1}(N)U_{i_1}^*\cdots U_{i_k}\xi_{i_k}(N)U_{i_k}^*) \\
&\hskip4cm -\tr_N^{\star n}(\iota_{i_1}(U_{i_1}\xi_{i_1}(N)U_{i_1}^*)\cdots
\iota_{i_k}(U_{i_k}\xi_{i_k}(N)U_{i_k}^*))\big|\,d\gamma_{\U(N)}^{\otimes n} \\
&\quad<\delta+2(R+1)^m\theta
\end{align*}
for every $N\ge N_0$. Since
$$
\lim_{N\to\infty}
\tr_N^{\star n}(\iota_{i_1}(\xi_{i_1}(N))\cdots\iota_{i_k}\xi_{i_k}(N)))
=\tau_{(X_1,\dots,X_n)}^\mathrm{free}(Z_{i_1}\cdots Z_{i_k}),
$$
the desired assertion follows.
\end{proof}

Let $W_{2,\mathrm{free}}(\tau_1,\tau_2)$ denote the free probabilistic $2$-Wasserstein
distance between $\tau_1, \tau_2 \in TS(\mathcal{A}_R)$ introduced by Biane and
Voiculescu \cite{BV} (see \cite[\S1.3]{HU1} for a brief summary fit to our arguments).
We need the next lemma comparing the free $2$-Wasserstein distance with the original
one (for measures) under the transformation
$\lambda \in \mathcal{P}(\mathrm{SU}(N)^n) \mapsto
\widehat{\lambda} \in TS(\mathcal{A}_R)$ defined in \eqref{F-3.1}.

\begin{lemma}\label{L-3.4}
For any $\lambda_1, \lambda_2 \in \mathcal{P}(\mathrm{SU}(N)^n)$ one has
\begin{equation*}
W_{2,\mathrm{free}}(\widehat{\lambda}_1,\widehat{\lambda}_2) \leq
\frac{2R}{\sqrt{N}} W_{2,\Vert\cdot\Vert_{\HS}}(\lambda_1,\lambda_2) \leq
\frac{2R}{\sqrt{N}} W_{2,\mathrm{geod}}(\lambda_1,\lambda_2),  
\end{equation*}
where $W_{2,\Vert\cdot\Vert_{\HS}}$ and $W_{2,\mathrm{geod}}$ are the $2$-Wasserstein
distances for measures with respect to the Hilbert-Schmidt norm
$\Vert\cdot\Vert_{\HS}$ and the geodesic distance, respectively.
\end{lemma}

\begin{proof}
The proof goes along the same line as that of \cite[Lemma 1.3]{HU1} with slight
modifications in the following two points. First, let $\Pi(\lambda_1,\lambda_2)$
denote the set of all probability measures on $\mathrm{SU}(N)^n\times\mathrm{SU}(N)^n$
whose left and right marginal measures are $\lambda_1$ and $\lambda_2$, respectively.
For each $\pi \in \Pi(\lambda_1,\lambda_2)$ we associate the state
$\widehat{\pi} \in TS(\mathcal{A}_R\star\mathcal{A}_R)$ via (the free product of two
copies of) the $*$-homomorphism sending each $Z_i$ to $U_i \xi_i(N)U_i^*$ as above.
Then one can easily observe that 
\begin{equation*} 
W_{2,\mathrm{free}}(\widehat{\lambda}_1,\widehat{\lambda}_2) \leq
\sqrt{\int_{\mathrm{SU}(N)^n}\int_{\mathrm{SU}(N)^n}
\sum_{i=1}^n \Vert U_i \xi_i(N) U_i^* - V_i^* \xi_i(N) V_i^*\Vert_{\HS}^2\,d\pi}
\end{equation*}   
for any $\pi \in \Pi(\lambda_1,\lambda_2)$, where the first integration is for
$(U_1,\dots,U_n)$ and the second for $(V_1,\dots,V_n)$. 
Secondly, we need the following elementary estimate:
\begin{equation*} 
\Vert U_i\xi_i(N)U_i^* - V_i\xi_i(N) V_i^*\Vert_{\HS} \leq
2 \Vert\xi_i(N)\Vert_\infty \Vert U_i - V_i\Vert_{\HS} \leq
2R \Vert U_i - V_i\Vert_{\HS},  
\end{equation*}
which is the reason why $R$ appears in the desired inequality. Finally, the latter
inequality is trivial because the geodesic distance majorizes the Hilbert Schmidt norm
distance.
\end{proof}   

We are now in a position to show the following transportation cost
inequality. Since $W_{2,\mathrm{free}}$ is indeed a metric, this yields the
implication from $\chi_\mathrm{orb}(X_1,\dots,X_n) = 0$ to the freeness of
$X_1,\dots,X_n$, thus proving the theorem.

\begin{prop}\label{P-3.5}
\begin{equation*} 
W_{2,\mathrm{free}}\big(\tau_{(X_1,\dots,X_n)},
\tau_{(X_1,\dots,X_n)}^\mathrm{free}\big) \leq
4\biggl(\max_{1\le i\le n}\|X_i\|_\infty\biggr)
\sqrt{-\chi_\mathrm{orb}(X_1,\dots,X_n)}.
\end{equation*}
\end{prop}

\begin{proof}
The proof is also same as that of \cite[Theorem 2.2]{HU1}, and thus we only give an
outline. Since the Ricci curvature of $\mathrm{SU}(N)^n$ (with respect to the inner
product induced from $\mathrm{Re}\,\mathrm{Tr}_N$) is known to be the constant $N/2$,
the transportation cost inequality 
\begin{equation*} 
W_{2,\mathrm{geod}}(\lambda_m,\gamma_{\mathrm{SU}(N_m)}^{\otimes n}) \leq
\sqrt{\frac{4}{N_m} S(\lambda_m,\gamma_{\mathrm{SU}(N_m)}^{\otimes n})}   
\end{equation*}
holds due to \cite{OV}, where $S(\cdot\,,\cdot)$ is the relative entropy. Since
$$
S(\lambda_m,\gamma_{\mathrm{SU}(N_m)}^{\otimes n})
= -\log\gamma_{\mathrm{SU}(N_m)}^{\otimes n}(\Gamma_m),
$$
we have by Lemma \ref{L-3.4}
\begin{equation*}
W_{2,\mathrm{free}}(\widehat{\lambda}_m,
\widehat{\gamma_{\mathrm{SU}(N_m)}^{\otimes n}}) \leq
4R \sqrt{-\frac{1}{N_m}\log\gamma_{\mathrm{SU}(N_m)}^{\otimes n}(\Gamma_m)}.
\end{equation*}
The desired inequality follows as $m\rightarrow\infty$ thanks to \eqref{F-3.2},
Lemmas \ref{L-3.2} and \ref{L-3.3} together with the joint lower
semicontinuity of $W_{2,\mathrm{free}}$.
\end{proof} 

\section{Generalization of $\chi_{\mathrm{orb}}$ to hyperfinite random multi-variables}
\setcounter{equation}{0}

For $1\le i\le n$ let $\bX_i = (X_{i1},\dots,X_{ir(i)})$ be a
non-commutative self-adjoint random multi-variable (called a {\it random multi-variable}
for short), which means a tuple consisting of self-adjoint random variables in
$(M,\tau)$. What we want here is to
generalize the orbital free entropy $\chi_{\mathrm{orb}}$ for random variables to
that for those multi-variables $\bX_1,\dots,\bX_n$. But there is a serious difficulty
in so doing in the general setting because we have no right counterpart of the map
$\Phi_N$ in \eqref{F-1.1} for the $n$-tuple space $(M_N^{sa})^n$. However, the
description of $\chi_\mathrm{orb}$ as
$\chi_{\mathrm{orb}}(X_1,\dots,X_n:\xi_1,\dots,\xi_n)$ (see Lemma \ref{L-2.4}) and
Jung's lemma (Lemma \ref{L-1.2}) allow us to define
$\chi_{\mathrm{orb}}(\bX_1,\dots,\bX_n)$ only when all
$W^*(\bX_i) := W^*(X_{i1},\dots,X_{i r(i)})$'s are hyperfinite. Throughout this
section we assume that $\bX_1,\dots,\bX_n$ are all {\it hyperfinite} in this sense.
Now, the definition of the orbital free entropy $\chi_\orb(\bX_1,\dots,\bX_n)$ is
similar to Definition \ref{D-2.3} as follows.

\begin{definition}\label{D-4.1}{\rm
For each $1\le i\le n$ let us choose a sequence $\{\Xi_i(N)\}$ consisting of
$r(i)$-tuples $\Xi_i(N) = (\xi_{i1}(N),\dots,\xi_{ir(i)}(N))$ of
$\xi_{ij}(N)\in M_N^{sa}$, $N\in\bN$, such that $\Xi_i(N)$ converges to $\bX_i$ in the
distribution sense (or in mixed moments) as $N \rightarrow \infty$. (Such a sequence
always exists due to the hyperfiniteness for $\bX_i$.) Define
$\Gamma_\mathrm{orb}(\bX_1,\dots,\bX_n:\Xi_1(N),\dots,\Xi_n(N);N,m,\delta)$ to be the
set of all $n$-tuples $(U_1,\dots,U_n) \in \mathrm{U}(N)^n$ such that 
\begin{equation*}
\left|\mathrm{tr}_N(U_{i_1}\xi_{i_1 j_1}(N)U_{i_1}^*\cdots
U_{i_k}\xi_{i_k j_k}(N) U_{i_k}^*)
- \tau(X_{i_1 j_1}\cdots X_{i_k j_k}) \right| < \delta
\end{equation*}
for all $1 \leq i_t \leq n$, $1 \leq j_t \leq r(i_t)$, $1\le t\le k$ and
$1 \leq k \leq m$, that is, $(U_i\Xi_i(N)U_i^*)_{i=1}^n\in
\Gamma(\bX_1\sqcup\cdots\sqcup\bX_n;\allowbreak N,m,\delta)$, where
$U_i\Xi_i(N)U_i^*$ means $(U_i\xi_{i1}(N)U_i^*,\dots,U_i\xi_{ir(i)}U_i^*)$. Then
we define
\begin{align*} 
&\chi_\mathrm{orb}(\bX_1,\dots,\bX_n) \\
&\quad:= \lim_{m\rightarrow\infty, \delta \searrow 0}\limsup_{N\rightarrow\infty}
\frac{1}{N^2}\log\gamma_{\mathrm{U}(N)}^{\otimes n}(
\Gamma_\mathrm{orb}(\bX_1,\dots,\bX_n:\Xi_1(N),\dots,\Xi_n(N);N,m,\delta)).
\end{align*}
}\end{definition}

If each $\bX_i$ consists of a single random variable, then the above
$\chi_\mathrm{orb}(\bX_1,\dots,\bX_n)$ clearly coincides with the $\chi_\orb$ in \S2
by definition. Moreover, the above definition is satisfactory as shown in
the next lemma. The proof is similar to that of Lemma \ref{L-2.4}.

\begin{lemma}\label{L-4.2}
$\chi_\mathrm{orb}(\bX_1,\dots,\bX_n)$ is independent of the choice of
$(\{\Xi_1(N)\},\dots,\{\Xi_n(N)\})$.
\end{lemma}

\begin{proof}
Let $(\{\Xi'_1(N)\},\dots,\{\Xi'_n(N)\})$ be another approximating $n$-tuple. By
symmetry it suffices to show that for each $\delta>0$ and $m \in \mathbb{N}$,
\begin{align*} 
&\gamma_{\mathrm{U}(N)}^{\otimes n}(
\Gamma_\mathrm{orb}(\bX_1,\dots,\bX_n:\Xi'_1(N),\dots,\Xi'_n(N);N,m,\delta/2)) \\
&\qquad\leq \gamma_{\mathrm{U}(N)}^{\otimes n}(
\Gamma_\mathrm{orb}(\bX_1,\dots,\bX_n:\Xi_1(N),\dots,\Xi_n(N);N,m,\delta)) 
\end{align*} 
for all sufficiently large $N$. Since $\Xi_i(N)$ and $\Xi'_i(N)$ converge to the same
$\bX_i$ in distribution, by Lemma \ref{L-1.2} one can choose an $N_0 \in \mathbb{N}$
so that for every $N \geq N_0$ there is an $n$-tuple
$(V_1(N),\dots,V_n(N)) \in \mathrm{U}(N)^n$ satisfying 
\begin{equation*} 
\Vert V_i(N)\xi_{ij}(N)V_i(N)^* - \xi'_{ij}(N)\Vert_{m,\mathrm{tr}_N}
< \frac{\delta}{2m(R+1)^{m-1}} 
\end{equation*} 
for all $1 \leq j \leq r(i)$ and $1\le i\le n$, where
$$
R:=\sup\{\|\xi_{ij}(N)\|_{m,\tr_N},\|\xi'_{ij}(N)\|_{m,\tr_N}:
1\le j\le r(i),\,1\le i\le n,\,N\in\bN\}\ (<+\infty).
$$
If $(U_1,\dots,U_n) \in \Gamma_\mathrm{orb}(\bX_1,\dots,\bX_n:
\Xi'_1(N),\dots,\Xi'_n(N);N,m,\delta/2)$ with $N\ge N_0$, then as in the proof of
Lemma \ref{L-2.4} we get
\begin{align*} 
&\big|\mathrm{tr}_N(U_{i_1}(V_{i_1}(N)\xi_{i_1 j_1}(N)V_{i_1}(N)^*)U_{i_1}^*
\cdots U_{i_k}(U_{i_k}(N)\xi_{i_k j_k}(N)U_{i_k}(N)^*)U_{i_k}^*) \\
&\qquad - \tau(X_{i_1 j_1}\cdots X_{i_k j_k})\big| \\
&\quad\leq 
m(R+1)^{m-1}\max_{i,j}\Vert V_i(N)\xi_{ij}(N)V_i(N)^*
-\xi'_{ij}(N)\Vert_{m,\mathrm{tr}_N} + \frac{\delta}{2} \\
&\quad< \delta  
\end{align*} 
for all $1 \leq i_t \leq n$, $1 \leq j_t \leq r(i_t)$, $1 \leq t \leq k$ and
$1 \leq k \leq m$. This means that 
\begin{align*} 
&\Gamma_\mathrm{orb}(\bX_1,\dots,\bX_n:\Xi'_1(N),\dots,\Xi'_n(N);N,m,\delta/2) \\
&\quad\subset \Gamma_\mathrm{orb}(\bX_1,\dots,\bX_n:
\Xi_1(N),\dots,\Xi_n(N);N,m,\delta)\cdot(V_1(N),\dots,V_n(N))
\end{align*} 
for all $N \geq N_0$. Hence we have the desired assertion thanks to the right
invariance of $\gamma_{\mathrm{U}(N)}$. 
\end{proof} 

\begin{remark}{\rm 
Jung's result \cite{Ju4} (or Lemma \ref{L-1.2}) says that the above definition
of\break $\chi_\mathrm{orb}(\bX_1,\dots,\bX_n)$ can work only when all the $\bX_i$'s
are hyperfinite since the microstates are not concentrated in a single (approximate)
unitary orbit for a general random multi-variable. Thus, to define
$\chi_\orb$ for general random multi-variables, we need an appropriate way to gather
together all unitary orbits
without ``overlap" in each matricial level. One potential way is to use the space of
unitarily equivalent classes of $*$-representations of $W^*(\bX_i)$ in $R^\omega$,
which plays a similar role of $\mathbb{R}_\geq^N$  for random variables. Note that
the restriction of diagonal matrices to $\mathbb{R}_\geq^N$ is needed in
the definition of $\chi_\mathrm{orb}$ to avoid ``overlap"; indeed, if $\mathbb{R}^N$
is used in place of $\mathbb{R}^N_\geq$, then the space of ``orbital microstates" has
the ``overlap" coming from the symmetry of $\mathfrak{S}_N$ acting on the eigenvalue
space $\mathbb{R}^N$. One more way we considered is to use a suitable fundamental
domain of the diagonal action of $\mathrm{U}(N)$ on
$\Gamma_R(X_{i1},\dots,X_{ir(i)};N,m,\delta)$ as a role of $\mathbb{R}^N_\geq$,
but we encountered some difficulty in this approach.
}\end{remark}

Except for the relation between $\chi_\mathrm{orb}$ and $\chi$ (Theorem \ref{T-2.6}),
all basic properties of $\chi_\mathrm{orb}$ can be extended to hyperfinite random
multi-variables in the same way, which are summarized in the next proposition. Note
that the assertion of Theorem \ref{T-2.6} is meaningless for hyperfinite random
multi-variables because both sides of the equality are $-\infty$ as long as at least
one of the $\bX_i$'s is not a single variable.

Only (4)--(8) of the proposition are somewhat non-trivial. Note that
(6) is the $\chi_\orb$ counterpart of \cite[Remark 9.2\,(e)]{V6} while
it is just a byproduct of (5). The proofs of (5), (7) and (8) are
essentially same as before in the case of $\chi_\mathrm{orb}$ for random variables;
for example, Lemma \ref{L-1.2} is used in place of Lemma \ref{L-1.1}. We will sketch
them and leave the full details to the reader. 

\begin{prop}\label{P-4.4}
$\chi_\mathrm{orb}$ for hyperfinite random multi-variables enjoys the following
properties: 
\begin{itemize}  
\item[(1)] $\chi_\mathrm{orb}(\bX) = 0$ for any single $\bX$.
\item[(2)] $\chi_\mathrm{orb}(\bX_1,\dots,\bX_n) \leq 0$.
\item[(3)] $\chi_\mathrm{orb}(\bX_1,\dots,\bX_n) \leq
\chi_\mathrm{orb}(\bX_1,\dots,\bX_k) + \chi_\mathrm{orb}(\bX_{k+1},\dots,\bX_n)$.   
\item[(4)] If $\bX_i^{(k)}=(X_{i1}^{(k)},\dots,X_{ir(i)}^{(k)})$ are
hyperfinite random multi-variables for $1\le i\le n$ and $k\in\bN$ such that
$\bX_1^{(k)}\sqcup\cdots\sqcup\bX_n^{(k)}\to
\bX_1\sqcup\cdots\sqcup\bX_n$ in the distribution sense as $k\to\infty$, then
\begin{equation*} 
\chi_\mathrm{orb}(\bX_1,\dots,\bX_n) \geq 
\limsup_{k\rightarrow\infty} \chi_\mathrm{orb}(\bX_1^{(k)},\dots,\bX_n^{(k)}). 
\end{equation*}  
\item[(5)] $\chi_\mathrm{orb}(\bX_1,\dots,\bX_n)$ depends only upon
$W^*(\bX_1),\dots,W^*(\bX_n)$; more precisely, 
$$
\chi_\orb(\bX_1,\dots,\bX_n)=\chi_\orb(\bX'_1,\dots,\bX'_n)
$$
for hyperfinite random multi-variables $\bX_i$ and $\bX'_i$ with
$W^*(\bX_i)=W^*(\bX'_i)$, $1\le i\le n$,
where the numbers of variables in $\bX_i$ and in $\bX'_i$ may be different.
\item[(6)] If $\bY_1,\dots,\bY_n$ are random multi-variables such that
$\bY_i\subset W^*(\bX_i)$ for $1\le i\le n$, then
$$
\chi_\orb(\bX_1,\dots,\bX_n)\le\chi_\orb(\bY_1,\dots,\bY_n).
$$
\item[(7)] If $\bX_1$ is freely independent of $\bX_2,\dots,\bX_n$, then
$$
\chi_\orb(\bX_1,\bX_2,\dots,\bX_n)=\chi_\orb(\bX_2,\dots,\bX_n).
$$
\item[(8)] $\chi_\mathrm{orb}(\bX_1,\dots,\bX_n) = 0$ if and only if
$\bX_1,\dots,\bX_n$ are freely independent.   
\end{itemize} 
\end{prop}

\begin{proof}
(4)\enspace For each $k\in\bN$ let $(\{\Xi_1^{(k)}(N)\},\dots,\{\Xi_n^{(k)}(N)\})$ be
an approximating $n$-tuple for $(\bX_1^{(k)},\allowbreak\dots,\bX_n^{(k)})$. For any
$\alpha<\limsup_k\chi_\orb(\bX_1^{(k)},\dots,\bX_n^{(k)})$ one can choose a sequence
$k_1<k_2<\dots$ such that $\chi_\orb(\bX_1^{(k_m)},\dots,\bX_n^{(k_m)})>\alpha$ and
\begin{equation}\label{F-4.1}
|\tau(X_{i_1j_1}^{(k_m)}\cdots X_{i_\ell j_\ell}^{(k_m)})
-\tau(X_{i_1j_1}\cdots X_{i_\ell j_\ell})|<{1\over m}
\end{equation}
for all $1\le i_t\le n$, $1\le j_t\le n(i_t)$, $1\le t\le\ell$ and $1\le\ell\le m$.
Furthermore, one can find a sequence $N_1<N_2<\dots$ such that for every $m\in\bN$
and $1\le i\le n$,
\begin{equation}\label{F-4.2}
\Xi_i^{(k_m)}(N)\in\Gamma(\bX_i^{(k_m)};N,m,1/m)
\quad\mbox{if $N\ge N_m$},
\end{equation}
and also
\begin{equation}\label{F-4.3}
{1\over N_m^2}\log\gamma_{\U(N)}^{\otimes n}
(\Gamma_\orb(\bX_1^{(k_m)},\dots,\bX_n^{(k_m)}:
\Xi_1^{(k_m)}(N_m),\dots,\Xi_n^{(k_m)}(N_m);N_m,m,1/m))>\alpha.
\end{equation}
For $1\le i\le n$ define
$$
\Xi_i(N):=\Xi_i^{(k_m)}(N)\quad
\mbox{if $N_m\le N<N_{m+1}$, $m\in\bN$}.
$$
By \eqref{F-4.1} and \eqref{F-4.2}, for $1\le i\le n$ we get
$\Xi_i(N)\in\Gamma(\bX_i;N,m,2/m)$ if $N_m\le N<N_{m+1}$, $m\in\bN$; hence
$(\{\Xi_1(N)\},\dots,\{\Xi_n(N)\})$ is an approximating $n$-tuple for
$(\bX_1,\dots,\bX_n)$. For each $m\in\bN$, if $(U_1,\dots,U_n)$ is in
$$
\Gamma_\orb(\bX_1^{(k_m)},\dots,\bX_n^{(k_m)}:
\Xi_1^{(k_m)}(N_m),\dots,\Xi_n^{(k_m)}(N_m);N_m,m,1/m),
$$
then $(U_i\Xi_i(N_m)U_i^*)_{i=1}^n=(U_i\Xi_i^{(k_m)}(N_m)U_i^*)_{i=1}^n$ is
in $\Gamma(\bX_1^{(k_m)}\sqcup\cdots\sqcup\bX_n^{(k_m)};N_m,m,1/m)$.
Since this set of microstates is included in
$\Gamma(\bX_1\sqcup\dots\sqcup\bX_n;N_m,m,2/m)$
thanks to \eqref{F-4.1}, it follows that
\begin{align*}
&\Gamma_\orb(\bX_1^{(k_m)},\dots,\bX_n^{(k_m)}:
\Xi_1^{(k_m)}(N_m),\dots,\Xi_n^{(k_m)}(N_m);N_m,m,1/m) \\
&\qquad\subset\Gamma_\orb(\bX_1,\dots,\bX_n:
\Xi_1(N_m),\dots,\Xi_n(N_m);N_m,m,2/m).
\end{align*}
Hence, by \eqref{F-4.3} we have
$$
{1\over N_m^2}\log\gamma_{\U(N)}^{\otimes n}
(\Gamma_\orb(\bX_1,\dots,\bX_n:\Xi_1(N_m),\dots,\Xi_n(N_m);N_m,m,2/m))>\alpha
$$
for all $m\in\bN$. This immediately implies that
$\chi_\orb(\bX_1,\dots,\bX_n)\ge\alpha$, and the result follows.

(5)\enspace Let $\bX_i=(X_{i1},\dots,X_{ir(i)})$ and
$\bX'_i=(X'_{i1},\dots,X'_{ir'(i)})$ be as stated in the proposition,
and choose their approximating $n$-tuples $(\{\Xi_1(N)\},\dots,\{\Xi_n(N)\})$ and
$(\{\Xi'_1(N)\},\dots,\allowbreak\{\Xi'_n(N)\})$, respectively, with
$\Xi_i(N) = (\xi_{ij}(N))_{j=1}^{r(i)}$ and
$\Xi'_i(N) = (\xi'_{ij}(N))_{j=1}^{r'(i)}$. We may assume that
$\|\xi_{ij}(N)\|_\infty,\|\xi'_{ij'}(N)\|_\infty\le
\max_{i,j,j'}\{\|X_{ij}\|_\infty,\|X'_{ij'}\|_\infty\}$ for all $i,j,j'$
and $N$. Now, it suffices to prove that for each $m\in\bN$ and $\delta>0$ there are
an $m'\in\bN$, a $\delta'>0$ and an $N_0\in\bN$ such that
\begin{align}\label{F-4.4}
&\gamma_{\U(N)}^{\otimes n}(\Gamma_\orb((\bX_i)_{i=1}^n:
(\Xi_i(N))_{i=1}^n;N,m',\delta')) \nonumber\\
&\qquad\le\gamma_{\U(N)}^{\otimes n}(\Gamma_\orb((\bX'_i)_{i=1}^n:
(\Xi'_i(N))_{i=1}^n;N,m,\delta))
\end{align}
for all $N\ge N_0$. The proof is essentially same as that of \eqref{F-2.5} but
more complicated since the right-hand side of \eqref{F-4.4} contains $\Xi'_i(N)$
differently from \eqref{F-2.5}. The Kaplansky density theorem enables us to choose
non-commutative self-adjoint polynomials $P_{ij}$ of $r(i)$ indeterminates for
$1 \leq j \leq r'(i)$, $1\leq i \leq n$ such that $\|P_{ij}(\bX_i)\|_\infty\le R$ and
$\Vert X'_{ij} - P_{ij}(\bX_i)\Vert_{1,\tau}$ is arbitrarily small; hence
$\bX_1'\sqcup\cdots\sqcup\bX_n'$ is arbitrarily approximated by
$(P_{1j}(\bX_1))_{j=1}^{r'(1)}\sqcup\cdots
\sqcup(P_{nj}(\bX_n))_{j=1}^{r'(n)}$ in distribution. Since $\Xi'_i(N)\to\bX'_i$ and
$(P_{ij}(\Xi_i(N)))_{j=1}^{r'(i)}\to(P_{ij}(\bX_i))_{j=1}^{r'(i)}$ in distribution
as $N\to\infty$, by Lemma \ref{L-1.2} one can find an $N_0\in\bN$ such that for every
$N\ge N_0$ and $1\le i\le n$ there exists a $V_i(N)\in\U(N)$ for which
$\|P_{ij}(\Xi_i(N))-V_i(N)\xi'_{ij}(N)V_i(N)^*\|_{2,\tr_N}$ is arbitrarily small for
$1\le j\le r'(i)$. Then one can choose an $m'\in\bN$ and a $\delta'>0$
such that
\begin{align*}
&\Gamma_\orb((\bX_i)_{i=1}^n:(\Xi_i(N))_{i=1}^n;N,m',\delta')
\cdot(V_i(N))_{i=1}^n \\
&\qquad\subset\Gamma_\orb((\bX'_i)_{i=1}^n:(\Xi'_i(N))_{i=1}^n;N,m,\delta)
\end{align*}
for all $N\ge N_0$, implying \eqref{F-4.4}.

(6)\enspace
Letting $\bX'_i:=\bX_i\sqcup\bY_i$ for $1\le i\le n$ we have by (5)
$$
\chi_\orb(\bX_1,\dots,\bX_n)=\chi_\orb(\bX'_1,\dots,\bX'_n)
\le\chi_\orb(\bY_1,\dots,\bY_n),
$$
since the latter inequality is obvious by definition.

(7)\enspace The proof is completely same as that of Proposition
\ref{P-2.9}; just replace $X_i$, $\xi_i(N)$ by $\bX_i$, $\Xi_i(N)$. See also
Proposition \ref{P-4.7} for its generalization. 

(8)\enspace The assertions (1) and (7) show that the freeness implies
$\chi_\orb=0$. The converse is proven by extending the transportation cost inequality
in Proposition \ref{P-3.5} to hyperfinite random multi-variables. The proof is same as
before, so only a few remarks are mentioned here. Set
$R := \max_{i,j} \Vert X_{ij}\Vert_\infty$ and let $\mathcal{A}_R$ be the universal
free product of $r(1)+\cdots+r(n)$ copies of $C[-R,R]$ with canonical generators
$Z_{ij}$ for $1 \le j \le r(i)$, $1\le i \le n$. By the $*$-homomorphism sending each
$Z_{ij}$ to $X_{ij}$ we obtain $\tau_{(\bX_1,\dots,\bX_n)} \in TS(\mathcal{A}_R)$ as
in \S3. Also, for every $\lambda \in \mathcal{P}(\mathrm{SU}(N)^n)$ we associate
$\widehat{\lambda} \in TS(\cA_R)$ in the same manner as in \S3 by the integral over
the unitary orbit $\{(U_1\Xi_1(N) U_1^*,\dots,U_n\Xi_n(N) U_n^*) :
(U_1,\dots,U_n) \in \mathrm{SU}(N)^n\}$ with respect to $\lambda$. Then, the
counterparts of Lemmas \ref{L-3.2} and \ref{L-3.3} are proven exactly in the same
way. Indeed, applying Lemma \ref{L-1.3} to 
\begin{equation*} 
\{(U_1,\dots,U_n) \in \mathrm{U}(N)^n : U_1\Xi_1(N)U_1^*,\dots,
U_n\Xi_n(N)U_n^*\ \text{are $(m,\delta)$-free}\}
\end{equation*}
one can show that $\lim_{N\to\infty}\widehat{\gamma_{\mathrm{SU}(N)}^{\otimes n}}=
\tau_{(\bX_1,\dots,\bX_n)}^\mathrm{free}$ weakly*, where
$\tau_{(\bX_1,\dots,\bX_n)}^\mathrm{free} \in TS(\mathcal{A}_R)$ is the free
product of the states $\tau_{\bX_i}$ on $C^*(Z_{ij},\dots,Z_{ir(i)})$ induced from the
original $\tau$ on $M$ via the $*$-homomorphism sending $Z_{ij}$ to $X_{ij}$ for
$1 \le j \le r(i)$. With these the same argument as before proves  
\begin{equation*} 
W_{2,\mathrm{free}}(\tau_{(\bX_1,\dots,\bX_n)},
\tau^\mathrm{free}_{(\bX_1,\dots,\bX_n)}) \le 4R
\sqrt{-\chi_\orb(\bX_1,\dots,\bX_n)},  
\end{equation*} 
from which we get the conclusion.
\end{proof}

Next, we introduce the $\chi_\orb(\cdots:\bv)$ in the presence of unitary random
variables, which will be necessary in the next section.

First, let us recall the $\Gamma$-set of microstates approximating
$\bX=(X_1,\dots,X_n)$ in the presence of unitary random variables.
In addition to $\bX$ let $\bv=(v_1,\dots,v_\ell)$ be an $\ell$-tuple of unitary random
variables in $(M,\tau)$. For $N,m\in\bN$ and $\delta>0$ we denote by
$\Gamma(\bX,\bv;N,m,\delta)$ the set of all $(A_1,\dots,A_n,V_1,\dots,V_\ell)$ in
$(M_N^{sa})^n\times\U(N)^\ell$ such that
$$
|\tr_N(h(A_1,\dots,A_n,V_1,\dots,V_\ell))-\tau(h(\bX,\bv))|<\delta
$$
for all $*$-monomials $h$ of $n+\ell$ indeterminates of degree not greater
than $m$, and by $\Gamma(\bX:\bv;N,m,\delta)$ the set of all
$(A_1,\dots,A_n)\in(M_N^{sa})^n$ such that
$(A_1,\dots,A_n,V_1,\dots,V_\ell)\in\Gamma(\bX,\bv;\allowbreak N,m,\delta)$
for some $(V_1,\dots,V_\ell)\in\U(N)^\ell$.

\begin{definition}\label{D-4.5}
For $1\le i\le n$ choose a microstate sequence
$\Xi_i(N)=(\xi_{i1}(N),\dots,\xi_{ir(i)}(N))$ in $(M_N^{sa})^{r(i)}$,
$N\in\bN$, such that $\Xi_i(N)$ converges to $\bX_i$ in the distribution sense as
$N\to\infty$. Moreover, let $\bv=(v_1,\dots,v_\ell)$ be unitary random variables in
$(M,\tau)$. For $N,m\in\bN$ and $\delta>0$ define
$\Gamma_\orb(\bX_1,\dots,\bX_n:\Xi_1(N),\dots,\Xi_n(N):\bv;N,m,\delta)$ to be the set
of all $(U_1,\dots,U_n)\in\U(N)^n$ such that $(U_i\Xi_i(N)U_i^*)_{i=1}^n$ is in
$\Gamma(\bX_1,\dots,\bX_n:\bv;N,m,\delta)$. Then we define
the {\it orbital free entropy} of $(\bX_1,\dots,\bX_n)$ in the presence of $\bv$ by
\begin{align*}
&\chi_\orb(\bX_1,\dots,\bX_n:\bv) \\
&:=\lim_{m\rightarrow\infty,\,\delta\searrow0}\limsup_{N\to\infty}
{1\over N^2}\log\gamma_{\U(N)}^{\otimes n}
\bigl(\Gamma_\orb(\bX_1,\dots,\bX_n:\Xi_1(N),\dots,\Xi_n(N):\bv;N,m,\delta)\bigr).
\end{align*}
Similarly to Lemma \ref{L-4.2} the above definition of
$\chi_\orb(\bX_1,\dots,\bX_n:\bv)$ is independent of the choice of an approximating
$n$-tuple $(\{\Xi_1(N)\},\dots,\{\Xi_n(N)\})$.
\end{definition}

The next proposition can be regarded as the $\chi_\orb$-counterpart of
\cite[Proposition 10.4]{V6}. In what follows, $\chi_u(\cdots)$ means the free
entropy of unitary random variables (see \cite[\S6.5]{HP}).   

\begin{prop}\label{P-4.6}
Let $\bv = (v_1,\dots,v_n)$ be a freely independent $n$-tuple of unitary random
variables with $\chi_u(v_i) > -\infty$ for all $1 \le i \le n$. If $\bX_1,\dots,\bX_n$
are freely independent of $\bv$, then  
\begin{align*} 
\chi_\orb(\bX_1,\dots,\bX_n) &\le 
\chi_\orb(v_1\bX_1v_1^*,\dots,v_n\bX_nv_n^*
:\bv) \\
&\le \chi_\orb(v_1\bX_1v_1^*,\dots,v_n\bX_nv_n^*).
\end{align*}
In particular, when the above $\bX_1,\dots,\bX_n$ are single self-adjoint random
variables $X_1,\dots,X_n$, one has
\begin{equation*} 
\chi(X_1,\dots,X_n) \le \chi(v_1X_1v_1^*,\dots,v_nX_nv_n^*). 
\end{equation*}    
\end{prop}

\begin{proof}
The latter assertion follows immediately from the first thanks to Theorem \ref{T-2.6}.
For the first assertion it is enough to prove only the first inequality. Choose
$\Xi_i(N)$ as in Definition \ref{D-4.5} with
$\|\Xi_i(N)\|_\infty\le\|\bX_i\|_\infty$ for $1\le i\le n$ and $N\in\bN$, where
$\|\Xi_i(N)\|_\infty:=\max_{1\le j\le r(i)}\|\xi_{ij}\|_\infty$ and
$\|\bX_i\|_\infty:=\max_{1\le j\le r(i)}\|X_{ij}\|_\infty$. For $N,m\in\bN$ and
$\delta,\rho>0$ we write for short
\begin{align*}
\Phi(N,m,\delta)
&:=\Gamma_\orb((\bX_i)_{i=1}^n:(\Xi_i(N))_{i=1}^n;N,m,\delta), \\
\widehat\Psi(N,m,\rho)
&:=\Gamma(v_1\bX_1v_1^*\sqcup\cdots\sqcup v_n\bX_nv_n^*,\bv;N,m,\rho), \\
\widetilde\Psi(N,m,\rho)
&:=\Gamma(v_1\bX_1v_1^*\sqcup\cdots\sqcup v_n\bX_nv_n^*:\bv;N,m,\rho), \\
\Psi(N,m,\rho)
&:=\Gamma_\orb((v_i\bX_iv_i^*)_{i=1}^n:(\Xi_i(N))_{i=1}^n
:\bv;N,m,\rho).
\end{align*}
We define two probability measures $\mu_N$ and $\nu_N$ on $\U(N)^n$ by
\begin{align*}
\mu_N&:={1\over\gamma_{\U(N)}^{\otimes n}(\Phi(N,m,\delta))}
\,\gamma_{\U(N)}^{\otimes n}\big|_{\Phi(N,m,\delta)}, \\
\nu_N&:={1\over\gamma_{\U(N)}^{\otimes n}(\Gamma(\bv;N,2m,\delta))}
\,\gamma_{\U(N)}^{\otimes n}\big|_{\Gamma(\bv;N,2m,\delta)},
\end{align*}
where $\Gamma(\bv;N,2m,\delta)$ is the $\Gamma$-set of unitary microstates
in $\U(N)^n$ approximating $\bv$ (see \cite[\S6.5]{HP}). Here we may and
do assume that $\chi_\orb(\bX_1,\dots,\bX_n) > -\infty$ so that $\mu_N$ is
well-defined for all sufficiently large $N \in \mathbb{N}$. Also, note that $\nu_N$
is well-defined for all sufficiently large $N \in \bN$ thanks to the assumption of
free independence for $\bv$. Furthermore, define
\begin{align*}
\Omega(N,3m,\delta)
&:=\bigl\{(U_1,\dots,U_n,V_1,\dots,V_n)\in\U(N)^n\times\U(N)^n: \\
&\hskip2cm\mbox{$(U_i\Xi_i(N)U_i^*)_{i=1}^n$ and $(V_i)_{i=1}^n$ are
$(3m,\delta)$-free}\bigr\}.
\end{align*}
For every $m\in\bN$ and $\rho>0$ one can choose a $\delta>0$ such that if
$\bU=(U_1,\dots,U_n)\in\Phi(N,m,\delta)$,
$\bV=(V_1,\dots,V_n)\in\Gamma(\bv;N,2m,\delta)$ and
$(\bU,\bV)\in\Omega(N,3m,\delta)$,
then $((V_iU_i\Xi_i(N)U_i^*V_i^*)_{i=1}^n,\allowbreak\bV)\in\widehat\Psi(N,m,\rho)$.
Since $\mu_N\otimes\nu_N$ is invariant under the $\U(N)$-action given by
$(\bU,\bV)\mapsto(\bU,W\bV W^*)$, $W\bV W^*:=(WV_iW^*)_{i=1}^n$, for $W\in\U(N)$, it
follows from Lemma \ref{L-1.3} (see the proof of \cite[Corollary 2.14]{V-IMRN}) that
$(\mu_N\otimes\nu_N)(\Omega(N,3m,\delta))\ge1/2$ whenever $N$ is large enough
(depending only on $m,\delta$). For each such $N$ one can choose a
$\bV\in\Gamma(\bv;N,2m,\delta)$ such that
\begin{equation}\label{F-4.5}
{1\over2}\le\mu_N(\Omega(N,3m,\delta:\bV))
={\gamma_{\U(N)}^{\otimes n}(\Phi(N,m,\delta)\cap\Omega(N,3m,\delta:\bV))
\over\gamma_{\U(N)}^{\otimes n}(\Phi(N,m,\delta))},
\end{equation}
where $\Omega(N,3m,\delta:\bV):=\{\bU\in\U(N)^n:
(\bU,\bV)\in\Omega(N,3m,\delta)\}$. From the above choice of $\delta$ we have
$$
\bigl\{(V_iU_i\Xi_i(N)U_i^*V_i^*)_{i=1}^n:\bU\in
\Phi(N,m,\delta)\cap\Omega(N,3m,\delta:\bV)\bigr\}
\subset\widetilde\Psi(N,m,\rho),
$$
that is,
$$
\bigl\{\bV\bU=(V_iU_i)_{i=1}^n:
\bU\in\Phi(N,m,\delta)\cap\Omega(N,3m,\delta:\bV)\bigr\}
\subset\Psi(N,m,\rho).
$$
Thanks to the left invariance of $\gamma_{\U(N)}$, this and \eqref{F-4.5} imply that
$$
{1\over2}\gamma_{\U(N)}^{\otimes n}(\Phi(N,m,\delta))
\le\gamma_{\U(N)}^{\otimes n}(\Psi(N,m,\rho)).
$$
Therefore,
$$
\chi_\orb(\bX_1,\dots,\bX_n)
\le\limsup_{N\to\infty}{1\over N^2}\log\gamma_{\U(N)}^{\otimes n}(\Psi(N,m,\rho)),
$$
implying the required inequality.
\end{proof}

The next proposition is exactly the $\chi_\orb$-counterpart of
\cite[Theorem 3.8]{V-IMRN}. 

\begin{prop}\label{P-4.7}
Let $\bv=(v_1,\dots,v_n)$ be unitary random variables. If $(\bX_1,v_1)$ is freely
independent of $\bX_2,\dots,\bX_n$ and $v_2,\dots,v_n$, then
$$
\chi_\orb(\bX_1,\dots,\bX_n:\bv)
=\chi_\orb(\bX_1:v_1)+\chi_\orb(\bX_2,\dots,\bX_n:v_2,\dots,v_n)
$$
whenever $\bX_1$ is regular in the presence of $v_1$, that is, replacing the $\limsup$
as $N\to\infty$
by $\liminf$ gives the same value in the definition $\chi_\orb(\bX_1:v_1)$.
\end{prop}

\begin{proof} 
Since the subadditivity
\begin{align}\label{F-4.6}
&\chi_\orb(\bX_1,\dots,\bX_n:v_1,\dots,v_n) \nonumber\\
&\quad\le\chi_\orb(\bX_1,\dots,\bX_k:v_1,\dots,v_k)
+\chi_\orb(\bX_{k+1},\dots,\bX_n:v_{k+1},\dots,v_n)
\end{align}
is obvious by definition, it suffices to show inequality $\ge$ for the required
equality. We can assume that $\chi_\orb(\bX_1:v_1)>-\infty$ and
$\chi_\orb(\bX_2,\dots,\bX_n:v_2,\dots,v_n)>-\infty$. We choose $\Xi_i(N)$ as in the
previous proof and for each $N,m\in\bN$ and $\delta,\rho>0$ write
\begin{align*}
\Phi(N,m,\delta)&:=\Gamma_\orb(\bX_1:\Xi_1(N):v_1;N,m,\delta) \\
&\qquad
\times\Gamma_\orb((\bX_i)_{i=2}^n:(\Xi_i(N))_{i=2}^n:(v_i)_{i=2}^n;N,m,\delta), \\
\Psi(N,m,\rho)&:=\Gamma_\orb((\bX_i)_{i=1}^n:(\Xi_i(N))_{i=1}^n:\bv;N,m,\rho).
\end{align*}
The assumption guarantees that $\Phi(N,m,\delta)$ is not of
$\gamma_{\U(N)}^{\otimes n}$-measure $0$ for all $N$ large enough. We will prove
that for each $m\in\bN$ and $\rho>0$ there is a $\delta>0$ such that
\begin{equation}\label{F-4.7}
{\gamma_{\U(N)}^{\otimes n}(\Psi(N,m,\rho)\cap\Phi(N,m,\delta))
\over\gamma_{\U(N)}^{\otimes n}(\Phi(N,m,\delta))}\ge{1\over2}
\end{equation}
for all sufficiently large $N$. The proof is similar to that of
\cite[Lemma 3.5]{V-IMRN}. First, note that
$\Gamma_\orb(\bX_1:\Xi_1(N):v_1;N,m,\delta)$ is invariant under the left action
$U_1\mapsto UU_1$ for $U\in\U(N)$. Hence the probability measure
$$
\mu_N:={1\over\gamma_{\U(N)}^{\otimes n}(\Phi(N,m,\delta))}
\,\gamma_{\U(N)}^{\otimes n}\big|_{\Phi(N,m,\delta)}
$$
is invariant under the same action of $\U(N)$ to the only first component. Next, for
any $m\in\bN$ and $\rho>0$, one can choose a $\delta>0$ so that if
$(\bA_1,V_1)\in\Gamma(\bX_1,v_1;N,m,\delta)$ with
$\|\bA_1\|_\infty\le\|\bX_1\|_\infty$ and $((\bA_i)_{i=2}^n,(V_i)_{i=2}^n)\in
\Gamma((\bX_i)_{i=2}^n,(v_i)_{i=2}^n;N,m,\delta)$ with
$\|\bA_i\|_\infty\le\|\bX_i\|_\infty$ and if $(\bA_1,V_1)$ and
$((\bA_i)_{i=2}^n,(V_i)_{i=2}^n)$ are $(m,\delta)$-free, then
$((\bA_i)_{i=1}^n,(V_i)_{i=1}^n)\in\Gamma((\bX_i)_{i=1}^n,(v_i)_{i=1}^n;\allowbreak
N,m,\rho))$. Lemma \ref{L-1.3} implies that
$$
\gamma_{\U(N)}(\{U\in\U(N):
\mbox{$(U\bA_1U^*,UV_1U^*)$ and $((\bA_i)_{i=2}^n,(V_i)_{i=2}^n)$ are
$(m,\delta)$-free}\})\ge{1\over2}
$$
for every $\bA_i$ and $V_i$ as above whenever $N$ is sufficiently large (depending
only on $m,\delta$). Then it follows that
$$
\gamma_{\U(N)}\bigl(\{U\in\U(N):(UU_1,(U_i)_{i=2}^n)\in\Psi(N,m,\rho)\}\bigr)
\ge{1\over2}
$$
for all $(U_1,(U_i)_{i=2}^n)\in\Phi(N,m,\delta)$ whenever $N$ is sufficiently large.
This implies that
$$
\mu_N(\Psi(N,m,\rho))
=\int_{\U(N)^n}\Biggl(\int_{\U(N)}\1_{\Psi(N,m,\rho)}(UU_1,(U_i)_{i=2}^n)
\,d\gamma_{\U(N)}(U)\Biggr)d\mu_N\ge{1\over2},
$$
implying \eqref{F-4.7}. Therefore, we obtain
\begin{align*}
&\limsup_{N\to\infty}{1\over N^2}\log
\gamma_{\U(N)}^{\otimes n}(\Psi(N,m,\rho)) \\
&\qquad\ge\limsup_{N\to\infty}{1\over N^2}\log
\gamma_{\U(N)}^{\otimes n}(\Phi(N,m,\delta)) \\
&\qquad\ge\liminf_{N\to\infty}{1\over N^2}\log
\Gamma_\orb(\bX_1:\Xi_1(N):v_1;N,m,\delta) \\
&\qquad\quad+\limsup_{N\to\infty}{1\over N^2}\log
\Gamma_\orb((\bX_i)_{i=2}^n:(\Xi_i(N))_{i=2}^n:(v_i)_{i=2}^n;N,m,\delta),
\end{align*}
and the desired inequality follows thanks to the regularity assumption
of $\bX_1$ (in the presence of $v_1$).
\end{proof}

\section{Orbital free entropy dimension}
\setcounter{equation}{0}

The microstate free entropy dimension $\delta$ and its modified one $\delta_0$
due to Voiculescu \cite{V2,V3} are defined for self-adjoint random variables based on
the microstate free entropy $\chi$ and the semicircular deformation  
\begin{equation}\label{F-5.1} 
(X_1+\varepsilon S_1,\dots,X_n + \varepsilon S_n), \qquad \varepsilon>0,
\end{equation}
where $(S_1,\dots,S_n)$ is a free semicircular system freely independent
of given self-adjoint random variables $X_1,\dots,X_n$. In this section we will
introduce the orbital version $\delta_{0,\mathrm{orb}}$ of $\delta_0$ (and also
$\delta_\mathrm{orb}$ of $\delta$), or in other words the dimension counterpart of
the orbital free entropy $\chi_\mathrm{orb}$ discussed in the previous sections.
Our essential idea to define $\delta_{0,\mathrm{orb}}$ is to replace $\chi$ by
$\chi_\mathrm{orb}$ and more importantly the semicircular deformation \eqref{F-5.1}
by the so-called liberation process
\begin{equation}\label{F-5.2} 
(v_1(t)X_1v_1(t)^*,\dots,v_n(t)X_nv_n(t)^*),\qquad t>0, 
\end{equation}
introduced by Voiculescu \cite{V6}, where $(v_1(t),\dots,v_n(t))$ is a
free $n$-tuple of multiplicative free unitary Brownian motions (see \cite{Bi})
freely independent of the $X_i$'s. The idea to use the liberation process
goes back to our attempt to define the dimension counterpart of
$\chi_\mathrm{proj}$; note that the space of projections with fixed traces is not
closed under the semicircular deformation \eqref{F-5.1} while it is under the
liberation process \eqref{F-5.2}.

Throughout the rest of this section, let $\bX_1,\dots,\bX_n$ be
hyperfinite random multi-variables in $(M,\tau)$ as treated in \S4.
           
\begin{definition}\label{D-5.1}{\rm
Let $\bv(t)=(v_1(t),\dots,v_n(t))$, $t\ge0$, be a freely independent
$n$-tuple of multiplicative free unitary Brownian motions (see \cite{Bi}) with
$v_i(0)=\1$ chosen to be freely independent of $\bX_1,\dots,\bX_n$.
(We may always assume that such extra variables exist in $(M,\tau)$.) Write
$v_i(t)\bX_iv_i(t)^*:=(v_i(t)X_{i1}v_i(t)^*,\dots,v_i(t)X_{ir(i)}v_i(t)^*)$ and
define the {\it modified orbital free entropy dimension} of $(\bX_1,\dots,\bX_n)$ by
$$
\delta_{0,\orb}(\bX_1,\dots,\bX_n)
:=\limsup_{\eps\searrow0}
{\chi_\orb(v_1(\eps)\bX_1v_1(\eps)^*,\dots,v_n(\eps)\bX_nv_n(\eps)^*:
\bv(\eps))\over|\log\eps^{1/2}|}.
$$
One may also define the orbital free entropy dimension $\delta_\orb(\bX_1,\dots,\bX_n)$
in the same manner by using
$\chi_\orb(v_1(\eps)\bX_1v_1(\eps)^*,\dots,v_n(\eps)\bX_nv_n(\eps)^*)$
without the presence of $\bv(\eps)$. However, we will deal with only
$\delta_{0,\orb}$ in this paper. 
}\end{definition}

\begin{remark}\label{R-5.2}{\rm
Let $(X_1,\dots,X_n)$ be an $n$-tuple of self-adjoint random variables and $\bv$ a
tuple of unitary random variables in $(M,\tau)$. The proof of Theorem \ref{T-2.6}
can be slightly modified to obtain
$$
\chi(X_1,\dots,X_n:\bv)=\chi_\orb(X_1,\dots,X_n:\bv)+\sum_{i=1}^n\chi(X_i).
$$
Applying this to $(v_1(\eps)X_1v_1(\eps)^*,\dots,v_n(\eps)X_nv_n(\eps)^*)$ and
$\bv(\eps)$ yields
\begin{align*}
&\chi(v_1(\eps)X_1v_1(\eps)^*,\dots,v_n(\eps)X_nv_n(\eps)^*:\bv(\eps)) \\
&\qquad=\chi_\orb(v_1(\eps)X_1v_1(\eps)^*,\dots,v_n(\eps)X_nv_n(\eps)^*:\bv(\eps))
+\sum_{i=1}^n\chi(X_i).
\end{align*}
Consequently, if $\chi(X_i)>-\infty$ for all $1\le i\le n$, then we have
$$
\delta_{0,\orb}(X_1,\dots,X_n)=\limsup_{\eps\searrow0}
{\chi(v_1(\eps)X_1v_1(\eps)^*,\dots,v_n(\eps)X_nv_n(\eps)^*:
\bv(\eps))\over|\log\eps^{1/2}|}.
$$
This formula might serve as the definition of $\delta_{0,\orb}$ for random variables
$(X_1,\dots,X_n)$ such that $\chi(X_i)>-\infty$ for $1\le i\le n$. However, it does
not make sense for hyperfinite random multi-variables
$(\bX_1,\dots,\bX_n)$ since $\chi(v_1(\eps)\bX_1v_1(\eps)^*\sqcup
\dots\sqcup v_n(\eps)\bX_nv_n(\eps)^*:\bv(\eps))=-\infty$ as long as at least one
of the $\bX_i$'s is not a single variable.
}\end{remark}

The next proposition summarizes properties of $\delta_{0,\orb}$;
(1)--(3) are rather obvious. The assertion (4) says that
$\delta_{0,\orb}(\bX_1,\dots,\bX_n)$ can be regarded as the (modified) orbital
free entropy dimension of the hyperfinite subalgebras $W^*(\bX_1),\dots,W^*(\bX_n)$.
Note that (6) is the orbital counterpart of \cite[Proposition 6.10]{V3}. Also, note
that (7) is the $\delta_{0,\orb}$-counterpart of Proposition \ref{P-2.9}, which 
slightly strengthens the second assertion of (6).

\begin{prop}\label{P-5.3}
$\delta_{0,\orb}$ for hyperfinite random multi-variables enjoys the following
properties{\rm:}
\begin{itemize}
\item[(1)] $\delta_{0,\orb}(\bX)=0$ for a single multi-variable $\bX$.
\item[(2)] $\delta_{0,\orb}(\bX_1,\dots,\bX_n) \le0$.
\item[(3)] $\delta_{0,\orb}(\bX_1,\dots,\bX_n)
\le\delta_{0,\orb}(\bX_1,\dots,\bX_k)
+\delta_{0,\orb}(\bX_{k+1},\dots,\bX_n)$ for every $1\le k<n$.
\item[(4)] $\delta_{0,\orb}(\bX_1,\dots,\bX_n)$ depends only upon
$W^*(\bX_1),\dots,W^*(\bX_n)$.
\item[(5)] If $\bY_1,\dots,\bY_n$ are random multi-variables such that
$\bY_i\subset W^*(\bX_i)$ for $1\le i\le n$, then
$$
\delta_{0,\orb}(\bX_1,\dots,\bX_n)\le\delta_{0,\orb}(\bY_1,\dots,\bY_n).
$$
\item[(6)] If $\chi_\orb(\bX_1,\dots,\bX_n)>-\infty$, then
$\delta_{0,\orb}(\bX_1,\dots,\bX_n)=0$. 
In particular, $\delta_{0,\mathrm{orb}}(\bX_1,\dots,\allowbreak\bX_n)=0$ if
$\bX_1,\dots,\bX_n$ are freely independent. 
\item[(7)] If $\bX_1$ is freely independent of $\bX_2,\dots,\bX_n$, then
$$
\delta_{0,\orb}(\bX_1,\bX_2\dots,\bX_n)
=\delta_{0,\orb}(\bX_2,\dots,\bX_n).
$$
\end{itemize}
\end{prop}

\begin{proof}
Since $\chi_\orb(\bX)=0$ for a single $\bX$, (1) is contained in (6). (2) is trivial
since $\chi_\orb(\bX_1,\dots,\allowbreak\bX_n:\bv)\le0$ for any $\bX_1,\dots,\bX_n$
and $\bv$. (3) follows from the subadditivity \eqref{F-4.6}.

(4)\enspace
For $1\le i\le n$ let $\bX'_i=(X'_{i1},\dots,X'_{ir'(i)})$ be another
random multi-variable with $W^*(\bX'_i)=W^*(\bX_i)$. To show the assertion, it
suffices to prove the equality of the modified orbital free entropies
$$
\chi_\orb(\bX_1,\dots,\bX_n:\bv)
=\chi_\orb(\bX'_1,\dots,\bX'_n:\bv)
$$
in the presence of unitary random variables $\bv$. But the proofs of Propositions
\ref{P-2.8} and \ref{P-4.4}\,(5) can be easily modified to prove this, so the details
are omitted.

(5) follows immediately from (4) as in the proof of Proposition \ref{P-4.4}\,(6).

(6)\enspace
Since $\chi_u(v_i(\eps)) > -\infty$ for every $\eps>0$
(see e.g., \cite[Proposition 10.10]{V6}), Proposition \ref{P-4.6} shows that  
$$
\chi_\orb(v_1(\eps)\bX_1v_1(\eps)^*,\dots,v_n(\eps)\bX_nv_n(\eps)^*:\bv(\eps))
\ge \chi_\orb(\bX_1,\dots,\bX_n)
$$
for every $\eps > 0$, from which the desired assertion immediately
follows. The latter assertion follows from Proposition \ref{P-4.4}\,(8).

(7)\enspace
The proof of Proposition \ref{P-4.6} shows that for every $m\in\bN$ and $\rho>0$
there is a $\delta>0$ such that
\begin{align*}
&{1\over2}\gamma_{\U(N)}\bigl(\Gamma_\orb(\bX_1:\Xi_1(N);N,m,\delta)\bigr) \\
&\qquad\le\gamma_{\U(N)}\bigl(\Gamma(v_1(\eps)\bX_1v_1(\eps)^*:
\Xi_1(N):v_1(\eps);N,m,\rho)\bigr)
\end{align*}
for all sufficiently large $N$. Since $\Gamma_\orb(\bX_1:\Xi_1(N);N,m,\delta)$ is
the whole $\U(N)$ whenever $N$ is large enough, $v_1(\eps)\bX_1v_1(\eps)^*$ is regular
in the presence of $v_1(\eps)$ as in Proposition \ref{P-4.7} and
$\chi_\orb(v_1(\eps)\bX_1v_1(\eps)^*:v_1(\eps))=0$ for every $\eps>0$. Therefore,
Proposition \ref{P-4.7} shows that 
\begin{align*} 
&\chi_\orb(v_1(\eps)\bX_1v(\eps)^*,v_2(\eps)\bX_2v_2(\eps)^*,\dots,
v_n(\eps)\bX_nv_n(\eps)^*:\bv(\eps)) \\
&\qquad= \chi_\orb(v_1(\eps)\bX_1v(\eps)^*,\dots,v_n(\eps)\bX_nv_n(\eps)^*:
v_2(\eps),\dots,v_n(\eps)) 
\end{align*} 
for every $\eps>0$, which immediately implies the required equality.
\end{proof}

Now, we examine how Jung's covering/packing approach \cite{Ju1,Ju2} to $\delta_0$
works for $\delta_{0,\mathrm{orb}}$ introduced above.
First, let us recall the notions of covering/packing numbers. Let $(\cX,d)$ be
a Polish space and $\Gamma\subset\cX$. Consider $\Gamma$ as a metric space with the
restriction of $d$ on $\Gamma$. For each $\eps>0$ we denote by $K_\eps(\Gamma)$
the minimum number of open $\eps$-balls covering $\Gamma$, and by $P_\eps(\Gamma)$
the maximum number of elements in a family of mutually disjoint open $\eps$-balls
in $\Gamma$, where $\eps$-balls in $\Gamma$ are taken as subsets of $\Gamma$. 
Those numbers will sometimes be denoted by $K_\eps(\Gamma,d)$ and
$P_\eps(\Gamma,d)$ to emphasize the metric $d$.
A subset $\{x_s:s\in S\}$ of $\Gamma$ is called an {\it $\eps$-net} of $\Gamma$ if
the open $\eps$-balls centered at $x_s$, $s\in S$, cover $\Gamma$, and also an
{\it $\eps$-separated set} of $\Gamma$ if the $\eps$-balls centered at $x_s$, $s\in S$,
are mutually disjoint. This definition is slightly different from that in
\cite{Sz2} but consistent with the definition of packing numbers used
here. Moreover, $\cN_\eps(\Gamma)$ stands for the open $\eps$-neighborhood of $\Gamma$.
Remark that $P_\eps(\Gamma) \ge K_{2\eps}(\Gamma) \ge P_{4\eps}(\Gamma)$ holds
in general, and thus if a lower/upper estimate for either $K_\eps(\Gamma)$ or
$P_\eps(\Gamma)$ was proven, then the essentially same estimate for the other would
immediately follow.

On the space $(M_N^{sa})^n$ ($\cong\bR^{nN^2}$) we consider the metric $d_2$ induced
from the Hilbert-Schmidt norm with respect to $\tr_N$.

\begin{definition}\label{D-5.4} 
Let $\bX_1,\dots,\bX_n$ and $\{\Xi_1(N)\},\dots,\{\Xi_n(N)\}$ be as in Definition
\ref{D-4.5}. Define the {\it orbital fractal free entropy dimension} of
$(\bX_1,\dots,\bX_n)$ by
$$
\delta_{1,\orb}(\bX_1,\dots,\bX_n)
:=\limsup_{\eps\searrow0}{\bK^\orb_\eps(\bX_1,\dots,\bX_n)\over|\log\eps|}-n
=\limsup_{\eps\searrow0}{\bP^\orb_\eps(\bX_1,\dots,\bX_n)\over|\log\eps|}-n,
$$
where
$$
\bK^\orb_\eps(\bX_1,\dots,\bX_n)
:=\lim_{m\to\infty,\delta\searrow 0}\limsup_{N\to\infty}
{1\over N^2}\log K_\eps\bigl(\Gamma_\orb((\bX_i)_{i=1}^n:
(\Xi_i(N))_{i=1}^n;N,m,\delta)\bigr)
$$
and $\bP^\orb_\eps(\bX_1,\dots,\bX_n)$ is similar with $P_\eps$ in place of $K_\eps$.
Indeed, it is seen from the proof of Lemma \ref{L-4.2} that the definitions of
$\bK^\orb_\eps(\bX_1,\dots,\bX_n)$, $\bP^\orb_\eps(\bX_1,\dots,\bX_n)$ and hence 
$\delta_{1,\orb}(\bX_1,\dots,\allowbreak\bX_n)$ are independent of the choice of
$(\{\Xi_1(N)\},\dots,\{\Xi_n(N))\}$.
\end{definition}

Let us then prove the equality $\delta_{0,\orb}=\delta_{1,\orb}$. Indeed,
the subtraction by $n$ in the above definition of $\delta_{1,\orb}$ is necessary to
get this equality. To do so we need a lemma, which says that $\bv(t)$ is regular;
namely, we have the same value if $\limsup$ is replaced by $\liminf$ in the definition
of $\chi_u(\bv(t))$ (see \cite[\S6.5]{HP}). Its proof is essentially same as in the
case of self-adjoint variables (the large deviation principle in \cite[5.4.10]{HP}
might be important).

\begin{lemma}\label{L-5.5}
Let $\bv(t)$, $t\ge0$, be as in Definition {\rm\ref{D-5.1}}. Then for every $t\ge0$,
$$
\lim_{m\to\infty,\delta\searrow 0}\liminf_{N\to\infty}
{1\over N^2}\log\gamma_{\U(N)}^{\otimes n}(\Gamma(\bv(t);N,m,\delta))
=\chi_u(\bv(t))=\sum_{i=1}^n\chi_u(v_i(t)).
$$
\end{lemma}

\begin{prop}\label{P-5.6}  
$$
\delta_{0,\orb}(\bX_1,\dots,\bX_n)
=\delta_{1,\orb}(\bX_1,\dots,\bX_n).
$$
\end{prop}

\begin{proof}
The idea of the proof is similar to that in \cite{Ju2}. First, by \cite[Lemma 8]{Bi}
there is a constant $K>0$ such that $\|v_i(t)-\1\|_\infty\le Kt^{1/2}$ for all
$0\le t\le1$. In what follows let $C:=K^2+2$, and let $N,m\in\bN$ and $\eps,\delta>0$
be arbitrary with restriction $\delta<\eps\le1$. Also let $\Xi_i(N)$ be
as in Definition \ref{D-4.5}.  

First let us prove the inequality $\ge$. One can choose a
$2(Cn\eps)^{1/2}$-separated subset $\{\bU_{Ns}=(U_{Nsi})_{i=1}^n : s\in S_N\}$ of
$\Gamma_\orb((\bX_i)_{i=1}^n:(\Xi_i(N))_{i=1}^n;N,m,\delta)$ with 
\begin{equation}\label{F-5.3}
|S_N|=P_{2(Cn\eps)^{1/2}}\bigl(\Gamma_\orb((\bX_i)_{i=1}^n:
(\Xi_i(N))_{i=1}^n;N,m,\delta)\bigr).
\end{equation}
(See the remark above Definition \ref{D-5.4} for the terminology of
``$\eps$-separated sets".)
Define two probability measures $\mu_N,\nu_N$ on $\U(N)^n$ by
\begin{align*}
\mu_N&:={1\over|S_N|}\sum_{s\in S_N}\delta_{\bU_{Ns}}
\quad\mbox{($\delta_{\bU_{Ns}}$ is the Dirac measure at $\bU_{Ns} \in \U(N)^n$)}, \\
\nu_N&:={1\over\gamma_{\U(N)}^{\otimes n}(\Gamma(\bv(\eps);N,m,\delta))}
\,\gamma_{\U(N)}^{\otimes n}\big|_{\Gamma(\bv(\eps);N,m,\delta)}.
\end{align*}
Write $\bU=(U_i)_{i=1}^n \in \U(N)^n$ etc., and set
\begin{align*}
\Omega(N,3m,\delta)&:=\bigl\{(\bU,\bV)\in\U(N)^n\times\U(N)^n:
\mbox{$(U_i\Xi_i(N)U_i^*)_{i=1}^n$ and $\bV$ are
$(3m,\delta)$-free}\bigr\}, \\
\Phi_s(N,3m,\delta)&:=\bigl\{\bV\in\Gamma(\bv(\eps);N,m,\delta):
\mbox{$(U_{Nsi}\Xi_i(N)U_{Nsi}^*)_{i=1}^n$ and $\bV$ are
$(3m,\delta)$-free}\bigr\}
\end{align*}
for $s\in S_N$. Since $\mu_N\otimes\nu_N$ is invariant under the $\U(N)$-action
$(\bU,\bV)\mapsto(\bU,W\bV W^*)$ for $W\in\U(N)$, by Lemma \ref{L-1.3} (as
\cite[Corollary 2.14]{V-IMRN}) we have
$$
{1\over2}\le(\mu_N\otimes\nu_N)(\Omega(N,3m,\delta))
={1\over|S_N|}\sum_{s\in S_N}\nu_N(\Phi_s(N,3m,\delta))
$$
so that
\begin{equation}\label{F-5.4}
\sum_{s\in S_N}\gamma_{\U(N)}^{\otimes n}(\Phi_s(N,3m,\delta))
\ge{1\over2}|S_N|\gamma_{\U(N)}^{\otimes n}(\Gamma(\bv(\eps);N,m,\delta))
\end{equation}
whenever $N$ is large enough. For every
$\bV=(V_1,\dots,V_n)\in\Gamma(\bv(\eps);N,m,\delta)$ we get
\begin{align}\label{F-5.5}
\|V_iU_{Nsi}-U_{Nsi}\|_{2,\tr_N}^2
&=\|V_i-I\|_{2,\tr_N}^2=\tr_N(2I-V_i-V_i^*) \nonumber\\
&\le\tau(2\1-v_i(\eps)-v_i(\eps)^*)+2\delta \nonumber\\
&<\|v_i(\eps)-\1\|_\infty^2+2\eps\le C\eps
\end{align}
so that $d_2(\bV\bU_{Ns},\bU_{Ns})<(Cn\eps)^{1/2}$. Hence it follows
that $\Phi_s(N,3m,\delta)\bU_{Ns}
:=\bigl\{\bV\bU_{Ns}:\bV\in\Phi_s(N,3m,\delta)\bigr\}$, $s\in S_N$,
are mutually disjoint. Furthermore, it is seen that for any $\rho>0$ we have
\begin{equation}\label{F-5.6}
\bigsqcup_{s\in S_N}\Phi_s(N,3m,\delta)\bU_{Ns}\subset
\Gamma_\orb((v_i(\eps)\bX_iv_i(\eps)^*)_{i=1}^n:
(\Xi_i(N))_{i=1}^n:\bv(\eps);N,m,\rho)
\end{equation}
if a sufficiently small $\delta\in(0,\eps)$ was chosen. By \eqref{F-5.6},
\eqref{F-5.4} and \eqref{F-5.3} we have
\begin{align*}
&\limsup_{N\to\infty}{1\over N^2}\log\gamma_{\U(N)}^{\otimes n}\bigl(
\Gamma_\orb((v_i(\eps)\bX_iv_i(\eps)^*)_{i=1}^n:
(\Xi_i(N))_{i=1}^n:\bv(\eps);N,m,\rho)\bigr) \\
&\qquad\ge\limsup_{N\to\infty}{1\over N^2}\log
\sum_{s\in S_N}\gamma_{\U(N)}^{\otimes n}(\Phi_s(N,3m,\delta)) \\
&\qquad\ge\limsup_{N\to\infty}{1\over N^2}\log\biggl({1\over2}|S_N|
\gamma_{\U(N)}^{\otimes n}(\Gamma(\bv(\eps);N,m,\delta))\biggr) \\
&\qquad\ge\mathbb{P}_{2(Cn\eps)^{1/2}}(\bX_1,\dots,\bX_n)
+\sum_{i=1}^n\chi_u(v_i(\eps))
\end{align*}
thanks to Lemma \ref{L-5.5}.
Note here that $\lim_{\eps\searrow0}\chi_u(v_i(\eps))/|\log\eps^{1/2}|=-1$ can be
easily derived from Voiculescu's computation \cite[Proposition 6.3]{V2} based on
\cite[Lemma 8]{Bi} and \cite[Proposition 1.6]{V6} since the spectrum of $v_i(\eps)$
is concentrated in a very small arc around $1$ for all sufficiently small $\eps >0$
(also see \cite[Proposition 6.1]{Sh1}). Hence the
above estimate implies the required inequality. 

Next let us prove the inequality $\le$. Let $\bU = (U_i)_{i=1}^n \in
\Gamma_\orb((v_i(\eps)\bX_i v_i(\eps)^*)_{i=1}^n:(\Xi_i(N))_{i=1}^n:\bv(\eps);
N,3m,\delta)$, which is accompanied by another
$\bV = (V_i)_{i=1}^n \in \Gamma(\bv(\eps);N,3m,\delta)$ by definition. One easily
observes that $\bV^*\bU = (V_i^* U_i)_{i=1}^n$ is in
$\Gamma_\orb((\bX_i)_{i=1}^n:(\Xi_i(N))_{i=1}^n;N,m,\delta)$. As similar to
\eqref{F-5.5} we have $d_2(\bU,\bV^*\bU) < (Cn\eps)^{1/2}$. Hence
$\Gamma_\orb((v_i(\eps)\bX_i v_i(\eps)^*)_{i=1}^n:(\Xi_i(N))_{i=1}^n:\bv(\eps);
N,3m,\delta)$ is included in
$\cN_{(Cn\eps)^{1/2}}(\Gamma_\orb((\bX_i)_{i=1}^n:(\Xi_i(N))_{i=1}^n;N,m,\delta))$.
Now choose an $\eps^{\1/2}$-net $\{\bU'_{Ns} : s\in S'_N\}$ of
$\Gamma_\orb((\bX_i)_{i=1}^n:(\Xi_i(N))_{i=1}^n;N,m,\delta)$
with
$$
|S'_N|=K_{\eps^{1/2}}(\Gamma_\orb((\bX_i)_{i=1}^n:(\Xi_i(N))_{i=1}^n;N,m,\delta)).
$$ 
Then $\Gamma_\orb(((v_i(\eps)\bX_i (v_i(\eps)^*)_{i=1}^n:(\Xi_i(N))_{i=1}^n:
\bv(\eps);N,3m,\delta)$ is clearly included in the union of the
$((Cn)^{1/2}+1)\eps^{1/2}$-balls $B_{((Cn)^{1/2}+1)\eps^{1/2}}(\bU'_{Ns})$ centered
at $\bU'_{Ns}$, $s \in S'_N$. By using the packing number estimate of $\U(N)$ due to
Szarek \cite{Sz1} one easily sees that there is a constant $C'>0$ independent of $N$
so that 
$$
\gamma_{\U(N)}^{\otimes n}(B_{((Cn)^{1/2}+1)\eps^{1/2}}(\bU'_{Ns}))
\le (C'((Cn)^{1/2}+1)\eps^{1/2})^{nN^2}
$$
as long as $\eps>0$ is small enough. Therefore we get 
\begin{align*}
&\gamma_{\U(N)}^{\otimes n}(\Gamma_\orb((v_i(\eps)\bX_i v_i(\eps)^*)_{i=1}^n:
(\Xi_i(N))_{i=1}^n:\bv(\eps);N,3m,\delta)) \\
&\qquad\le |S'_N|(C'((Cn)^{1/2}+1)\eps^{1/2})^{nN^2},   
\end{align*}
and hence 
\begin{align*}
&\frac{1}{N^2}\log\gamma_{\U(N)}^{\otimes n}
(\Gamma_\orb((v_i(\eps)\bX_i v_i(\eps)^*)_{i=1}^n:
(\Xi_i(N))_{i=1}^n:\bv(\eps);N,3m,\delta)) \\
&\qquad\le
\frac{1}{N^2}\log K_{\eps^{1/2}}(\Gamma_\orb((\bX_i)_{i=1}^n:
(\Xi_i(N))_{i=1}^n;N,m,\delta)) \\
&\qquad\qquad+ n\log\eps^{1/2} + n\log(C'((Cn)^{1/2}+1)). 
\end{align*}
Taking $\lim_{m\rightarrow\infty, \delta\searrow0}\limsup_{N\rightarrow\infty}$ of
both sides, we have 
\begin{align*}
&\chi_\orb(v_1(\eps)\bX_1 v_1(\eps)^*,\dots,v_n(\eps)\bX_n v_n(\eps)^*:\bv(\eps)) \\
&\qquad\le 
\bK^\orb_{\eps^{1/2}}(\bX_1,\dots,\bX_n) + n\log\eps^{1/2} + n\log(C'((Cn)^{1/2}+1)),  
\end{align*}
from which the desired inequality immediately follows.
\end{proof}

\begin{remark}\label{R-5.7}
The role of multiplicative free unitary Brownian motions $v_i(t)$ is not so essential
in the above proof of Proposition \ref{P-5.6}. In fact, besides the free independence
assumption, we used only the facts that $\|v_i(t)-\1\|_\infty\le Kt^{1/2}$ for small
$t\ge0$ and that $\lim_{\eps\searrow0}\chi_u(v_i(\eps))/|\log\eps^{1/2}|=-1$, while
Lemma \ref{L-5.5} is valid for general freely independent unitary random variables.
Consequently, we notice that the definition $\delta_{0,\orb}(\bX_1,\dots,\bX_n)$ in
Definition \ref{D-5.1} is equivalent when $\bv(t)$ is replaced by, for example,
$(e^{\sqrt{-1}\sqrt th_1},\dots,e^{\sqrt{-1}\sqrt th_n})$,
where $h_1,\dots,h_n$ are freely independent self-adjoint random variables with
$\chi(h_i)>-\infty$ for $1\le i\le n$ chosen to be freely independent of
$\bX_1,\dots,\bX_n$. The situation is similar to the case $\delta_0$ shown in
\cite{Ju2}.
\end{remark}

The main result of this section is the following exact relation between
$\delta_{0,\orb}$ and the usual $\delta_0$.

\begin{thm}\label{T-5.8} 
$$
\delta_0(\bX_1\sqcup\cdots\sqcup\bX_n) = \delta_{0,\orb}(\bX_1,\dots,\bX_n)
+ \sum_{i=1}^n \delta_0(\bX_i).
$$
\end{thm}

The rest of this section is devoted to the proof of the theorem. We will prove the
part ``$\le$" first and next ``$\ge$." The latter is more involved than the former. 

Let $\Xi_i(N)$ be chosen for $\bX_i$, $1\le i\le n$, as in Definition \ref{D-4.5}. For
$N,m\in\bN$ and $\delta>0$ define
\begin{align}\label{F-5.7}
&\Gamma(\bX_1\sqcup\cdots\sqcup\bX_n:(\Xi_i(N))_{i=1}^n;N,m,\delta) \nonumber\\
&\qquad:=\Gamma(\bX_1\sqcup\cdots\sqcup\bX_n;N,m,\delta)
\cap\{(U_i\Xi_i(N)U_i^*)_{i=1}^n:(U_1,\dots,U_n)\in\U(N)^n\},
\end{align}
where $U_i\Xi_i(N)U_i^*$ is as in Definition \ref{D-4.1}.
We need the following simple lemma. 

\begin{lemma}\label{L-5.9}
\begin{align*}
&\delta_0(\bX_1\sqcup\cdots\sqcup\bX_n) \\
&=\limsup_{\eps\searrow0}{1\over|\log\eps|}
\biggl(\lim_{m\to\infty,\delta\searrow 0}\limsup_{N\to\infty} 
{1\over N^2}\log K_\eps\bigl(
\Gamma(\bX_1\sqcup\cdots\sqcup\bX_n:(\Xi_i(N))_{i=1}^n;N,m,\delta)\bigr)\biggr).
\end{align*}
The same formula holds also when $K_\eps$ is replaced by $P_\eps$.
\end{lemma}

\begin{proof}
Thanks to Jung's covering/packing approach \cite{Ju2} to $\delta_0$ with additional
remarks \cite[p.\,455]{DH} and \cite[Lemma 2.2]{Ju3}, it suffices to show that for
every $m\in\bN$ and $\eps,\delta>0$ there are an $m_0\in\bN$ and a $\delta_0>0$
such that
$$
\Gamma(\bX_1\sqcup\cdots\sqcup\bX_n;N,m_0,\delta_0)
\subset\cN_\eps\bigl(\Gamma(\bX_1\sqcup\cdots\sqcup\bX_n:
(\Xi_i(N))_{i=1}^n;N,m,\delta)\bigr)
$$
for all sufficiently large $N$. But this can be easily verified by Lemma \ref{L-1.2}. 
\end{proof}

For a while fix an arbitrary $1 \le i \le n$. In what follows we assume that
$W^*(\bX_i)$ has both diffuse and atomic parts since this case is most involved and
needs all the ingredients of the proof. Let us decompose  
$$
W^*(\bX_i) = \bigoplus_{k=0}^{s(i)}M_{ik} 
$$
with $s(i) \in \bN\cup\{\infty\}$ such that $M_{i0}$ is diffuse and
$M_{ik}\cong M_{m_{ik}}(\bC)$ for $k\ge1$, and denote by $p_{ik}$ the central support
projection of $M_{ik}$, $k\ge0$. By Jung's result \cite{Ju1} one has 
\begin{equation}\label{F-5.8}
\delta_0(\bX_i) = 1 -\sum_{k=1}^{s(i)} \frac{\tau(p_{ik})^2}{m_{ik}^2}. 
\end{equation}
We choose and fix a matrix unit system
$\{e_{\alpha\beta}^{(ik)}:1\le\alpha,\beta\le m_{ik}\}$ of
$M_{ik}\cong M_{m_{ik}}(\bC)$ for each $k\ge1$. Let $\ell\in \bN$ be
arbitrary, and write $q_i^{(\ell)} := \1 - \sum_{k=0}^{s(i)\wedge\ell}p_{ik}$
with $s(i)\wedge\ell:=\min\{s(i),\ell\}$. We can choose,
for any sufficiently large $N\in\bN$, positive integers $n^{(\ell)}_{ik}(N)$ and
$m_{ik}^{(\ell)}(N)$ for $0\le k\le s(i)\wedge\ell$ such that
$n^{(\ell)}_{ik}(N) = m_{ik}m^{(\ell)}_{ik}(N)$, $k\ge1$, and
\begin{align}
&\sum_{k=0}^{s(i)\wedge\ell}n^{(\ell)}_{ik}(N) \le N, \nonumber\\
&\lim_{N\rightarrow\infty}\frac{n^{(\ell)}_{ik}(N)}{N}
= \tau(p_{ik}), \qquad 0 \le k\le s(i)\wedge\ell. \label{F-5.9}
\end{align} 
Moreover, choose orthogonal projections $P^{(\ell)}_{ik}(N)\in M_N(\bC)$
of rank $n^{(\ell)}_{ik}(N)$ for $0\le k\le s(i)\wedge\ell$ so that we can identify
$$
P^{(\ell)}_{ik}(N)M_N(\bC)P^{(\ell)}_{ik}(N)
= M_{m_{ik}}(\bC) \otimes M_{m^{(\ell)}_{ik}(N)}(\bC),\qquad 1\le k\le s(i)\wedge\ell.
$$
Under this identification, we set  
$$
\eta_{\alpha\beta}^{(ik\ell)}(N)
:= e_{\alpha\beta}^{(ik)}\otimes I_{m^{(\ell)}_{ik}(N)}
\in M_{m_{ik}}(\bC)\otimes M_{m^{(\ell)}_{ik}(N)}(\bC),
$$
for $1\le \alpha,\beta \le m_{ik}$ and $1\le k\le s(i)\wedge\ell$. Also we set
$$
Q_i^{(\ell)}(N) := I - \sum_{k=0}^{s(i)\wedge\ell}P^{(\ell)}_{ik}(N).
$$
 
\medskip
Now let us start the proof of the inequality $\le$. For each $\ell\in\bN$
fixed, we enlarge the given multi-variables $\bX_i = (X_{i1},\dots,X_{ir(i)})$,
$1 \le i\le n$, as follows: 
$$
\bX^{(\ell)}_i := \bX_i\sqcup\bigsqcup_{k=1}^{s(i)\wedge\ell}\Biggl(e_{11}^{(ik)},
\dots,e_{m_{ik}m_{ik}}^{(ik)},\sum_{\alpha,\beta=1}^{m_{ik}}e_{\alpha\beta}^{(ik)}
\Biggr),\qquad 1\le i\le n. 
$$
Since $W^*(\bX_i^{(\ell)}) = W^*(\bX_i)$, Proposition \ref{P-5.3}\,(4) gives
\begin{equation}\label{F-5.10}
\delta_{0,\orb}(\bX_1,\dots,\bX_n)
= \delta_{0,\orb}(\bX^{(\ell)}_1,\dots,\bX^{(\ell)}_n)
\end{equation} 
Moreover, since $\bX_i \subset \bX^{(\ell)}_i$, by \cite[Theorem 4.3]{V-IMRN} one has 
\begin{equation}\label{F-5.11}
\delta_0(\bX_1\sqcup\cdots\sqcup\bX_n)
\le \delta_0(\bX^{(\ell)}_1\sqcup\cdots\sqcup\bX^{(\ell)}_n).
\end{equation}  

The next lemma is plain to show by the use of Lemma \ref{L-1.2}.

\begin{lemma}\label{L-5.10} For each $1\le i\le n$ and for any sufficiently large
$N\in \bN$, one can find microstates $\xi^{(\ell)}_{ij}(N) \in M_N^{sa}$,
$1\le j\le r(i)$, in such a way that
$\|\xi^{(\ell)}_{ij}(N)\|_\infty \le \|\bX_i\|_\infty$, the $\xi_{ij}^{(\ell)}(N)$'s
are contained in
$$
P_{i0}^{(\ell)}(N)M_N(\bC)P_{i0}^{(\ell)}(N)\,\oplus\,\Biggl[
\bigoplus_{k=1}^{s(i)\wedge\ell}M_{m_{ik}}(\bC)\otimes\bC I_{m_{ik}^{(\ell)}(N)}
\Biggr]\,\oplus\,Q_i^{(\ell)}(N)M_N(\bC)Q_i^{(\ell)}(N), 
$$
and moreover  
\begin{align*}
\Xi_i^{(\ell)}(N) &:= 
\bigl(\xi_{i1}^{(\ell)}(N),\dots,\xi_{ir(i)}^{(\ell)}(N)\bigr) \\
&\qquad\sqcup\bigsqcup_{k=1}^{s(i)\wedge\ell}\Biggl(
\eta_{11}^{(ik\ell)}(N),\dots,\eta_{m_{ik}m_{ik}}^{(ik\ell)}(N),
\sum_{\alpha,\beta=1}^{m_{ik}}\eta_{\alpha\beta}^{(ik\ell)}(N)\Biggr)
\end{align*}
converges in the distribution sense to $\bX^{(\ell)}_i$ as
$N\rightarrow\infty$. 
\end{lemma}
\begin{proof} For each $1\le i\le n$, from hyperfiniteness one can choose
an approximating sequence of microstates
$$
(\xi_1(N),\dots,\xi_{r(i)}(N))\sqcup(P(N))\sqcup
\bigsqcup_{k=1}^{s(i)\wedge\ell}\Bigl(\zeta^{(k)}_{11}(N),\dots,
\zeta_{m_{ik}m_{ik}}^{(k)}(N),\zeta^{(k)}(N)\Bigr)\sqcup(Q(N))
$$
for
$$
(X_{i1},\dots,X_{ir(i)})\sqcup(p_{i0})\sqcup\bigsqcup_{k=1}^{s(i)\wedge\ell}
\Biggl(e_{11}^{(ik)},\dots,e_{m_{ik}m_{ik}}^{(ik)},
\sum_{\alpha,\beta=1}^{m_{ik}}e_{\alpha\beta}^{(ik)}\Biggr)\sqcup(q_i^{(\ell)})
$$
with the norm condition $\|\xi_j(N)\|_\infty \le \|X_{ij}\|_\infty$, $1\le j\le r(i)$.
Compare these $P(N)$, $Q(N)$, $\zeta^{(k)}_{\alpha\alpha}(N)$'s and
$\zeta^{(k)}(N)$'s with $P_{i0}^{(\ell)}(N)$, $Q^{(\ell)}_i(N)$,
$\eta_{\alpha\alpha}^{(ik\ell)}(N)$'s and
$\sum_{\alpha,\beta=1}^{m_{ik}}\eta_{\alpha\beta}^{(ik\ell)}(N)$'s in mixed moments.
By Lemma \ref{L-1.2}, for sufficiently large $N$ we then get unitaries $U_N$ which
intertwine the two families approximately in the sense that
$\|U_NP(N)U_N^*-P_{i0}^{(\ell)}(N)\|_{m,\tr_N}$ etc.\ go to $0$ as $N\to\infty$ for
all $m\in\bN$. Then one can get a new approximating sequence of microstates for
$X_{ij}p_{i0}$, $X_{ij}q_i^{(\ell)}$ ($1\le j \le r(i)$) and
$e_{11}^{(ik)},\dots,e_{m_{ik}m_{ik}}^{(ik)}$,
$\sum_{\alpha,\beta=1}^{m_{ik}}e_{\alpha\beta}^{(ik)}$ ($1\le k\le s(i)\wedge\ell$)
by sending the $\xi_j(N)$'s via $\mathrm{Ad}\,U_N$ and cutting with
$P_{i0}^{(\ell)}(N)$ or $Q_i^{(\ell)}(N)$ such that the part of those corresponding to
$e_{11}^{(ik)},\dots,e_{m_{ik}m_{ik}}^{(ik)}$,
$\sum_{\alpha,\beta=1}^{m_{ik}}e_{\alpha\beta}^{(ik)}$ ($1\le k\le s(i)\wedge\ell$)
are exactly $\eta_{11}^{(ik\ell)}(N),\dots,\eta_{m_{ik}m_{ik}}^{(ik\ell)}(N)$,
$\sum_{\alpha,\beta=1}^{m_{ik}}\eta_{\alpha\beta}^{(ik\ell)}(N)$
($1\le k\le s(i)\wedge\ell$). Then the desired microstates can easily
be made from those.
\end{proof}

Remark that the commutant of $\Xi^{(\ell)}_i(N)$ includes
$$
\bC P_{i0}^{(\ell)}(N)\,\oplus\,\Bigg[\bigoplus_{k=1}^{s(i)\wedge\ell}
\bC I_{m_{ik}}\otimes M_{m^{(\ell)}_{ik}(N)}(\bC)\Bigg]\,\oplus\,\bC Q_i^{(\ell)}(N).
$$
We denote by $\U_i^{(\ell)}(N)$ the unitary group of this algebra, i.e.,
\begin{equation}\label{F-5.12}
\U_i^{(\ell)}(N) := \bT P_{i0}^{(\ell)}(N)\,\oplus\,\Bigg[
\bigoplus_{k=1}^{s(i)\wedge\ell}I_{m_{ik}}\otimes\U(m_{ik}^{(\ell)}(N))
\Bigg]\,\oplus\,\bT Q_i^{(\ell)}(N).
\end{equation}
We then have
\begin{equation}\label{F-5.13}
\lim_{N\rightarrow\infty}\frac{1}{N^2}\dim_{\bR}\U_i^{(\ell)}(N)
= \lim_{N\rightarrow\infty}
\frac{2+\sum_{k=1}^{s(i)\wedge\ell}m_{ik}^{(\ell)}(N)^2}{N^2}
= \sum_{k=1}^{s(i)\wedge\ell}\frac{\tau(p_{ik})^2}{m_{ik}^2}. 
\end{equation}

Consider the embedding 
$$
\Psi_N^{(\ell)} : ([U_i])_{i=1}^n \in \prod_{i=1}^n \U(N)/\U_i^{(\ell)}(N) 
\mapsto (U_i\Xi_i^{(\ell)}(N)U_i^*)_{i=1}^n \in (M_N^{sa})^{n(\ell)}, 
$$
where $n(\ell)$ is the sum of the numbers of variables in $\bX_i^{(\ell)}$ for
$1\le i\le n$ and $[U_i]$ denotes the coset determined by $U_i\in\U(N)$. We introduce
the ``embedding" metric $d_{2,E}$ on the homogeneous space
$\prod_{i=1}^n \U(N)/\U_i^{(\ell)}(N)$ by 
$$
d_{2,E}(([U_i])_{i=1}^n,([V_i])_{i=1}^n)
:= \|\Psi_N^{(\ell)}(([U_i])_{i=1}^n) - \Psi_N^{(\ell)}(([V_i])_{i=1}^n)\|_{2,\tr_N}
$$
for $([U_i])_{i=1}^n,([V_i])_{i=1}^n \in \prod_{i=1}^n \U(N)/\U_i^{(\ell)}(N)$ with
$U_i,V_i \in \U(N)$, $1\le i \le n$. Another natural metric on
$\prod_{i=1}^n \U(N)/\U_i^{(\ell)}(N)$ is the quotient metric $d_{2,Q}$ induced from
$d_2$ on $\U(N)$ by 
$$
d_{2,Q}([U_i]_{i=1}^n,[V_i]_{i=1}^n)
:= \inf\Biggl\{d_2((U_i)_{i=1}^n,(V_iW_i)_{i=1}^n) :
(W_i)_{i=1}^n \in \prod_{i=1}^n \U_i^{(\ell)}(N) \Biggr\}.
$$
It is plain to see that 
\begin{equation}\label{F-5.14}
d_{2,E}(([U_i])_{i=1}^n,([V_i])_{i=1}^n)
\le C_\ell\,d_{2,Q}(([U_i])_{i=1}^n,([V_i])_{i=1}^n) 
\end{equation}
for all $([U_i])_{i=1}^n,([V_i])_{i=1}^n \in \prod_{i=1}^n \U(N)/\U_i^{(\ell)}(N)$,
where 
$$
C_\ell := 2n(\ell)\max\{\|\bX_i\|_\infty,\,m_{ik}:
1\le k\le s(i)\wedge\ell,\,1\le i\le n\}.
$$

Viewing $\prod_{i=1}^n\U(N)/\U_i^{(\ell)}(N)$ as
$\U(N)^n/\prod_{i=1}^n\U_i^{(\ell)}(N)$ we have the canonical quotient map 
$$
\Phi^{(\ell)}_N : \U(N)^n \rightarrow \prod_{i=1}^n\U(N)/\U_i^{(\ell)}(N).   
$$
Now let $\eps>0$ be arbitrary, and note that the $\eps$-covering number of
$\Gamma(\bX_1^{(\ell)}\sqcup\cdots\sqcup\bX_n^{(\ell)}:
(\Xi_i^{(\ell)}(N))_{i=1}^n;N,m,\delta)$ with respect to $d_2$ is equal to that of
$\Phi_N^{(\ell)}(\Gamma_\orb((\bX_i^{(\ell)})_{i=1}^n:
(\Xi^{(\ell)}_i(N))_{i=1}^n;\allowbreak N,m,\delta))$ with respect to $d_{2,E}$
since $\Psi_N^{(\ell)}$ isometrically maps the latter set to the former.
By \eqref{F-5.14} we hence get 
\begin{align}\label{F-5.15}
&K_\eps(\Gamma(\bX_1^{(\ell)}\sqcup\cdots\sqcup\bX_n^{(\ell)}:
(\Xi_i^{(\ell)}(N))_{i=1}^n; N,m,\delta)) \notag\\
&\qquad\le 
K_{C_\ell\eps}(\Phi_N^{(\ell)}(\Gamma_\orb((\bX_i^{(\ell)})_{i=1}^n:
(\Xi^{(\ell)}_i(N))_{i=1}^n;N,m,\delta))),   
\end{align}  
where the above right-hand side is counted with respect to $d_{2,Q}$. Here we need
the following simple (probably known) fact on the packing/covering numbers in
homogeneous spaces. For the convenience of the reader we give it with a proof. 

\begin{lemma}\label{L-5.11}{\rm (cf.~\cite[Lemma 6]{Sz1})} Let $G$ be a compact group
with a bi-invariant metric $d$, and $H$ be its closed subgroup. Let
$\pi : G \rightarrow G/H$ be the canonical quotient map sending $g \in G$ to the coset
$gH$, and equip $G/H$ with the quotient metric
$d_Q(g_1 H, g_2 H) := \min\{ d(g_1,g_2 h) : h \in H\}$. Then, for any
$\Gamma \subset G$ with $\pi^{-1}(\pi(\Gamma)) = \Gamma$ and for every $\eps>0$, one
has
$$
K_\eps(\Gamma) \ge 
K_\eps(H)\cdot P_{2\eps}(\pi(\Gamma)) \ge K_\eps(H)\cdot K_{4\eps}(\pi(\Gamma)).
$$
\end{lemma} 
\begin{proof}
The ball centered at $x$ of radius $r$ in a metric space is denoted by $B_r(x)$. One
can choose an $\varepsilon$-net $\{ g_i : i\in I \}$ of $\Gamma$ with cardinality
$|I| = K_{\varepsilon}(\Gamma)$, and a $2\varepsilon$-separated set
$\{ g'_j : j\in J\}$ of $\pi(\Gamma)$ with cardinality $|J| = P_{2\varepsilon}$. Let
$I_j := \{ i \in I : B_\varepsilon(g_i)\cap g'_j H \neq \emptyset\}$ for $j \in J$;
then it is clear that $I_{j_1} \cap I_{j_2} = \emptyset$ if $j_1 \neq j_2$. On the
other hand, $\{g'_j{}^{-1} g_i : i\in I_j\}$ gives an $\varepsilon$-net of $H$ so that
$|I_j| \ge K_\varepsilon(H)$ for all $j \in J$. Hence
$K_\varepsilon(\Gamma) = |I| \ge \sum_{j\in J} |I_j|
\ge K_\varepsilon(H)\cdot P_{2\varepsilon}(\pi(\Gamma))$. Then the assertion follows
thanks to the obvious relation between covering and packing numbers.   
\end{proof} 

By the above lemma and \eqref{F-5.15} we have 
\begin{align}
&K_\eps(\Gamma(\bX_1^{(\ell)}\sqcup\cdots\sqcup\bX_n^{(\ell)}:
(\Xi_i^{(\ell)}(N))_{i=1}^n;N,m,\delta)) \notag \\
&\quad\le 
K_{(C_\ell/4)\eps}(\Gamma_\orb((\bX_i^{(\ell)})_{i=1}^n:
(\Xi^{(\ell)}_i(N))_{i=1}^n;N,m,\delta))\cdot
K_{(C_\ell/4)\eps}\Biggl(\prod_{i=1}^n\U_i^{(\ell)}(N)\Biggr)^{-1}. \label{F-5.16}
\end{align} 
Identify   
\begin{align*}
\prod_{i=1}^n\U_i^{(\ell)}(N) 
&= \prod_{i=1}^n\Biggl(\bT\,\times\,\Biggl(\prod_{k=1}^{s(i)\wedge\ell}
\U(m_{ik}^{(\ell)}(N))\Biggr)\,\times\,\bT\Biggr) 
\end{align*}
in the obvious way, and consider the $\ell^\infty$-product metric induced from
$\|\cdot\|_{2,\tr_{m_{ik}^{(\ell)}(N)}}$ on $\U(m_{ik}^{(\ell)}(N))$.
Since the $m_{ik}$'s as well as $\ell$ are independent of $N$ in
\eqref{F-5.12}, this metric is clearly equivalent to the original metric $d_2$
(restricted on $\prod_{i=1}^n\U_i^{(\ell)}(N)$) uniformly in all sufficiently large
$N$. Then, by \cite[Theorem 7]{Sz2} with the help of \cite[Lemma 5]{Sz1}, there is a
constant $C'_\ell>0$ independent of $N$ such that
\begin{equation}\label{F-5.17}
K_{(C_\ell/4)\eps}\Biggl(\prod_{i=1}^n\U_i^{(\ell)}(N)\Biggr) 
\ge \biggl(\frac{C'_\ell}{(C_\ell/4)\eps}
\biggr)^{\sum_{i=1}^n\dim_\bR\U_i^{(\ell)}(N)}
\end{equation}
as long as $\eps > 0$ is small enough. By using \eqref{F-5.16}, \eqref{F-5.17} and
\eqref{F-5.13} we thus get 
\begin{align*} 
&\lim_{m\rightarrow\infty,\delta\searrow0}\limsup_{N\rightarrow\infty}
\frac{1}{N^2}\log K_\eps(\Gamma(\bX_1^{(\ell)}\sqcup\cdots\sqcup\bX_n^{(\ell)}:
(\Xi_i^{(\ell)}(N))_{i=1}^n;N,m,\delta)) \\
&\qquad\le 
\bK^\orb_{C_\ell\eps}(\bX_1^{(\ell)},\dots,\bX_n^{(\ell)})
+ \biggl(\log\biggl({C_\ell\over4C'_\ell}\biggr) +\log\eps\biggr)
\sum_{i=1}^n\sum_{k=1}^{s(i)\wedge\ell}\frac{\tau(p_{ik})^2}{m_{ik}^2}. 
\end{align*}
By Lemma \ref{L-5.9}, \eqref{F-5.11}, Proposition \ref{P-5.6} and
\eqref{F-5.10}, this implies that 
$$
\delta_0(\bX_1\sqcup\cdots\sqcup\bX_n) \le 
\delta_{0,\orb}(\bX_1,\dots,\bX_n)
+ n - \sum_{i=1}^n\sum_{k=1}^{s(i)\wedge\ell}\frac{\tau(p_{ik})^2}{m_{ik}^2}
$$
Since $\ell$ is arbitrary, we get the inequality $\le$ in Theorem \ref{T-5.8} thanks
to \eqref{F-5.8}. 

\medskip
Let us turn to the proof of the inequality $\ge$. Keep an arbitrary $\ell\in\bN$.
Since $M_{i0}$ is diffuse, one can choose
$p_{i01}^{(\ell)},\dots,p_{i0\ell}^{(\ell)} \in M_{i0}$ in such a way that
$p_{i01}^{(\ell)}+\cdots+p_{i0\ell}^{(\ell)} = p_{i0}$ and
$\tau(p_{i01}^{(\ell)}) = \cdots = \tau(p_{i0\ell}^{(\ell)})$ $(=\tau(p_{i0})/\ell)$.
For $1\le i\le n$ consider a new hyperfinite random multi-variables $\bY_i^{(\ell)}$
in $W^*(\bX_i)$ given by
$$
\bY_i^{(\ell)} := (p_{i01}^{(\ell)},\dots,p_{i0\ell}^{(\ell)})
\sqcup\bigsqcup_{k=1}^{s(i)\wedge\ell}
\Biggl(e_{11}^{(ik)},\dots,e_{m_{ik}m_{ik}}^{(ik)},
\sum_{\alpha,\beta=1}^{m_{ik}}e_{\alpha\beta}^{(ik)}\Biggr)\sqcup(q_i^{(\ell)}).
$$
Then, similarly to \eqref{F-5.10} one has
\begin{equation}\label{F-5.18}
\delta_{0,\orb}(\bX_1,\dots,\bX_n)
= \delta_{0,\orb}(\bX_1\sqcup\bY_1^{(\ell)},\dots,\bX_n\sqcup\bY_n^{(\ell)}).
\end{equation}
On the other hand, one can apply \cite[Corollary 4.2]{Ju3}, a corollary of
the so-called hyperfinite inequality due to Jung, to obtain the equality
\begin{equation}\label{F-5.19}
\delta_0(\bX_1\sqcup\dots\sqcup\bX_n)
= \delta_0(\bX_1\sqcup\bY_1^{(\ell)}\sqcup\dots\sqcup\bX_n\sqcup\bY_n^{(\ell)}) 
\end{equation}
unlike the previous \eqref{F-5.11}.
Let $n_{ik}^{(\ell)}(N)$, $P_{ik}^{(\ell)}(N)$ ($0\le k\le s(i)\wedge\ell$),
$m_{ik}^{(\ell)}(N)$ ($1\le k\le s(i)\wedge\ell$) and $Q_i^{(\ell)}(N)$ be as in the
proof of the part ``$\le$." Moreover, for any sufficiently large $N$, one can choose
$n_{i01}^{(\ell)}(N),\dots,n_{i0\ell}^{(\ell)}(N) \in \bN$ so that 
\begin{align}
&n_{i01}^{(\ell)}(N)+\cdots+n_{i0\ell}^{(\ell)}(N)
= n_{i0}^{(\ell)}(N), \nonumber\\
&\lim_{N\rightarrow\infty}\frac{n_{i0k}^{(\ell)}(N)}{N}
= {\tau(p_{i0})\over\ell}, \qquad 1\le k\le\ell. \label{F-5.20}
\end{align}
Then, choose orthogonal projections $P_{i0k}^{(\ell)}(N)\in M_N(\bC)$ of rank
$n_{i0k}^{(\ell)}(N)$ for $0\le k\le\ell$ such that
$\sum_{k=1}^\ell P_{i0k}^{(\ell)}(N) = P_{i0}^{(\ell)}(N)$. A special approximating
microstates for $\bY_i^{(\ell)}$ is given by
\begin{align*}
\widehat\Xi_i^{(\ell)}&:=(P_{i01}^{(\ell)}(N),\dots,P_{i0\ell}^{(\ell)}(N)) \\
&\qquad\sqcup\bigsqcup_{k=1}^{s(i)\wedge\ell}\Biggl(
\eta_{11}^{(ik\ell)}(N),\dots,\eta_{m_{ik},m_{ik}}^{(ik\ell)}(N),
\sum_{\alpha,\beta=1}^{m_{ik}}\eta_{\alpha\beta}^{(ik\ell)}(N)\Biggr)
\sqcup(Q_i^{(\ell)}(N)).
\end{align*}
The next lemma is proven in the same way as Lemma \ref{L-5.10}, so the
details are left to the reader.  

\begin{lemma}\label{L-5.12} For $1\le i\le n$ and for any sufficiently large $N$,
one can find microstates $\Xi_i^{(\ell)}(N)$ for $\bX_i$ in such a way that
$\|\Xi_i^{(\ell)}(N)\|_\infty\le\|\bX_i\|_\infty$ and
$\Xi_i^{(\ell)}(N)\sqcup\widehat\Xi_i^{(\ell)}(N)$ converges to
$\bX_i\sqcup\bY_i^{(\ell)}$ in the distribution sense as $N\to\infty$.
\end{lemma}

We denote by $\U_i^{(\ell)}(N)$ the unitary group of the commutant of
$\Xi_i^{(\ell)}(N)\sqcup\widehat\Xi_i^{(\ell)}(N)$ and by $\widehat\U_i^{(\ell)}(N)$
that of the commutant of $\widehat\Xi_i^{(\ell)}$. Since the commutant of
$\widehat\Xi_i^{(\ell)}$ is
\begin{align*}
&\bigoplus_{k=1}^\ell P_{i0k}^{(\ell)}(N)M_N(\bC)P_{i0k}^{(\ell)}(N) \\
&\qquad\quad\oplus\Biggl[\bigoplus_{k=1}^{s(i)\wedge\ell}\bC I_{m_{ik}}
\otimes M_{m_{ik}^{(\ell)}(N)}(\bC)\Biggr]
\oplus Q_i^{(\ell)}(N)M_N(\bC)Q_i^{(\ell)}(N),
\end{align*}
the real dimension of $\widehat\U_i^{(\ell)}(N)$ is
$$
\dim_\bR\widehat\U_i^{(\ell)}(N)=\sum_{k=1}^\ell n_{i0k}^{(\ell)}(N)^2
+\sum_{k=1}^{s(i)\wedge\ell}m_{ik}^{(\ell)}(N)^2
+\Biggl(N-n_{i0}^{(\ell)}(N)-\sum_{k=1}^{s(i)\wedge\ell}n_{ik}^{(\ell)}(N)\Biggr)^2
$$
so that by \eqref{F-5.9} and \eqref{F-5.20} we have
\begin{equation}\label{F-5.21}
\lim_{N\rightarrow\infty}\frac{N^2-\dim_\bR\widehat\U_i^{(t\ell)}(N)}{N^2}
= 1-{\tau(p_{i0})^2\over\ell}
-\sum_{k=1}^{s(i)\wedge\ell}{\tau(p_{ik})^2\over m_{ik}^2}-\tau(q_i^{(\ell)})^2.
\end{equation}

Introduce the embeddings
\begin{align*}
\Psi_N^{(\ell)}&:([U_i])_{i=1}^n\in\prod_{i=1}^n\U(N)/\U_i^{(\ell)}(N)
\mapsto\bigl(U_i(\Xi_i^{(\ell)}(N)\sqcup\widehat\Xi_i^{(\ell)}(N))U_i^*\bigr)_{i=1}^n
\in(M_N^{sa})^{n(\ell)}, \\
\widehat\Psi_N^{(\ell)}&:(\<U_i\>)_{i=1}^n
\in\prod_{i=1}^n\U(N)/\widehat\U_i^{(\ell)}(N)
\mapsto\bigl(U_i(\widehat\Xi_i^{(\ell)}(N))U_i^*\bigr)_{i=1}^n
\in(M_N^{sa})^{\widehat n(\ell)},
\end{align*}
where $n(\ell)$ is the sum of the numbers of variables in $\bX_i\sqcup\bY_i^{(\ell)}$,
$1\le i\le n$, and $\widehat n(\ell)$ is that of variables in $\bY_i^{(\ell)}$,
$1\le i\le n$. Moreover, we introduce the ``embedding" metric $\widehat d_{2,E}$ in
terms of $\widehat\Psi_N^{(\ell)}$ and the quotient metric $\widehat d_{2,Q}$ on the
homogeneous space $\prod_{i=1}^n \U(N)/\widehat\U_i^{(\ell)}(N)$ in the same way as
in the proof of the part ``$\le$." The ``embedding" metric $d_{2,E}$ in terms of
$\Psi_N^{(\ell)}$ is also introduced on $\prod_{i=1}^n\U(N)/\U_i^{(\ell)}(N)$. From
the trivial inclusions $\U_i^{(\ell)}(N)\subset\widehat\U_i^{(\ell)}(N)$ and
$\widehat\Xi_i^{(\ell)}(N)\subset\Xi_i^{(\ell)}(N)\sqcup\widehat\Xi_i^{(\ell)}(N)$,
we have the well-defined surjective map
$$
([U_i])_{i=1}^n\in\prod_{i=1}^n\U(N)/\U_i^{(\ell)}(N)\mapsto
(\<U_i\>)_{i=1}^n\in\prod_{i=1}^n\U(N)/\widehat\U_i^{(\ell)}(N)
$$
so that
\begin{equation}\label{F-5.22}
\widehat d_{2,E}((\<U_i\>)_{i=1}^n,(\<V_i\>)_{i=1}^n)
\le d_{2,E}(([U_i])_{i=1}^n,([V_i])_{i=1}^n).
\end{equation}
The next lemma is essentially \cite[Lemma 5.4]{Ju1}. In fact, the argument there works
well when $\bigoplus_{k=1}^\ell\bC p_{i0k}^{(\ell)}\oplus
\bigoplus_{k=1}^{s(i)\wedge\ell}M_{ik}\oplus\bC q_i^{(\ell)}$ and $\Xi_i^{(\ell)}(N)$
here play the roles of $M$ and $\pi$ there. Thus a chosen constant $C_\ell$ depends
only on $m_{ik}$, $1\le k\le s(i)\wedge\ell$, $1\le i\le n$, as well as $\ell$.

\begin{lemma}\label{L-5.13} There is a constant $C_\ell>0$ independent of $N$ such that 
$$
\widehat d_{2,Q}((\<U_i\>)_{i=1}^n,(\<V_i\>)_{i=1}^n)
\le C_\ell\,\widehat d_{2,E}((\<U_i\>)_{i=1}^n,(\<V_i\>)_{i=1}^n)
$$
for all $(\<U_i\>)_{i=1}^n, (\<V_i\>)_{i=1}^n \in
\prod_{i=1}^n \U(N)/\widehat\U_i^{(\ell)}(N)$.
\end{lemma}

Viewing $\prod_{i=1}^n\U(N)/\U_i^{(\ell)}(N)$ as
$\U(N)^n/\prod_{i=1}^n\U_i^{(\ell)}(N)$ we have the canonical quotient map 
$$
\Phi_N^{(\ell)} : \U(N)^n \rightarrow \prod_{i=1}^n\U(N)/\U_i^{(\ell)}(N),    
$$
and denote by $\mu_N^{(\ell)}$ the left-invariant probability measure on
$\prod_{i=1}^n\U(N)/\U_i^{(\ell)}(N)$ induced from $\gamma_{\U(N)}^{\otimes n}$.
In what follows, let $\eps,\delta>0$ be arbitrary with restriction $\delta <\eps$,
and for $N,m\in\bN$ we write for short
\begin{align*}
&\Gamma_\orb(\eps,N,3m,\delta)
:=\Gamma_\orb((v_i(\eps)(\bX_i\sqcup\bY_i^{(\ell)})v_i(\eps)^*)_{i=1}^n: \\
&\hskip5cm (\Xi_i^{(\ell)}(N)\sqcup\widehat\Xi_i^{(\ell)}(N))_{i=1}^n:
\bv(\eps);N,3m,\delta), \\
&\Gamma_\orb(N,m,\delta)
:=\Gamma_\orb((\bX_i\sqcup\bY_i^{(\ell)})_{i=1}^n:
(\Xi_i^{(\ell)}(N)\sqcup\widehat\Xi_i^{(\ell)}(N))_{i=1}^n;N,m,\delta), \\
&\Gamma(N,m,\delta)
:=\Gamma(\bX_1\sqcup\bY_1^{(\ell)}\sqcup\cdots\sqcup\bX_n\sqcup\bY_n^{(\ell)}:
(\Xi_i^{(\ell)}(N)\sqcup\widehat\Xi_n^{(\ell)}(N))_{i=1}^n:N,m,\delta)
\end{align*}
(see \eqref{F-5.7}). The following inequality is trivial:
\begin{equation}\label{F-5.23} 
\gamma_{\U(N)}^{\otimes n}(\Gamma_\orb(\eps,N,3m,\delta))
\le\mu_N^{(\ell)}(\Phi_N^{(\ell)}(\Gamma_\orb(\eps,N,3m,\delta))). 
\end{equation}

Assume $(U_i)_{i=1}^n \in \Gamma_\orb(\eps,N,3m,\delta)$; then there is a
$(V_i)_{i=1}^n \in \Gamma(\bv(\eps);N,3m,\delta)$ such that
$(V_i^*U_i)_{i=1}^n\in\Gamma_\orb(N,m,\delta)$ and hence for any
$\xi\in\Xi_i^{(\ell)}(N)\sqcup\widehat\Xi_i^{(\ell)}(N)$ we have
$$
\|U_i\xi U_i^* - V_i^* U_i\xi U_i^* V_i\|_{2,\tr_N} \le 2R(C\eps)^{1/2}
$$
so that
$$
d_{2,E}(([U_i])_{i=1}^n,([V_i^*U_i])_{i=1}^n) \le 2R(n(\ell)C\eps)^{1/2},
$$
where $R := \max\bigl\{\|\bX_i\sqcup\bY_i^{(\ell)}\|_\infty,\,m_{ik}:
1\le k\le s(i)\wedge\ell,\,1\le i\le n\bigr\}$
and $C>0$ is the same constant as in the proof of Proposition \ref{P-5.6}. Therefore
we get 
$$
\Phi_N^{(\ell)}(\Gamma_\orb(\eps,N,3m,\delta))
\subset\cN_{L\eps^{1/2}}(\Phi_N^{(\ell)}(\Gamma_\orb(N,m,\delta))), 
$$
where the right-hand side is the open $L\eps^{1/2}$-neighborhood of
$\Gamma_\orb(N,m,\delta)$ with respect to the metric $d_{2,E}$ and
$L:=2R(n(\ell)C)^{1/2}+1$. Here note that the $\eps^{1/2}$-covering number of
$\Phi_N^{(\ell)}(\Gamma_\orb(N,m,\delta))$ with respect to $d_{2,E}$ is equal to
$K_{\eps^{1/2}}(\Gamma(N,m,\delta))$ with respect to $d_2$ (as noted just above
\eqref{F-5.16}). Hence the above inclusion immediately implies that
\begin{equation}\label{F-5.24}
\mu_N^{(\ell)}(\Phi_N^{(\ell)}(\Gamma_\orb(\eps,N,3m,\delta)))
\le K_{\eps^{1/2}}(\Gamma(N,m,\delta))
\cdot\mu_N^{(\ell)}(\mathrm{Ball}((L+1)\eps^{1/2},d_{2,E})),
\end{equation}
where $\mathrm{Ball}((L+1)\eps^{1/2},d_{2,E})$ stands for the $(L+1)\eps^{1/2}$-ball
in $\prod_{i=1}^n\U(N)/\U_i^{(\ell)}(N)$ with respect to $d_{2,E}$. The measure of
this ball can be estimated from above by packing numbers as follows:
\begin{align*}
\mu_N^{(\ell)}(\mathrm{Ball}((L+1)\eps^{1/2},d_{2,E}))^{-1}
&\ge P_{(L+1)\eps^{1/2}}\Biggl(\prod_{i=1}^n\U(N)/\U_i^{(\ell)}(N),
d_{2,E}\Biggr) \\
&\ge P_{(L+1)\eps^{1/2}}\Biggl(\prod_{i=1}^n\U(N)/\widehat\U_i^{(\ell)}(N),
\widehat d_{2,E}\Biggr) \\
&\ge P_{C_\ell(L+1)\eps^{1/2}}\Biggl(\prod_{i=1}^n\U(N)/\widehat\U_i^{(\ell)}(N),
\widehat d_{2,Q}\Biggr).
\end{align*}
The second and the third inequalities in the above follow from \eqref{F-5.22} and
Lemma \ref{L-5.13}, respectively. Furthermore, the packing estimate due to
Jung \cite[Lemma 5.2 and \S8]{Ju1} (based on \cite{Sz2}) says that there is a constant
$C'_\ell>0$ independent of $N$ such that
$$
P_{C_\ell(L+1)\eps^{1/2}}\Biggl(\prod_{i=1}^n\U(N)/\widehat\U_i^{(\ell)}(N),
\widehat d_{2,Q}\Biggr)
\ge\biggl({C'_\ell\over C_\ell(L+1)\eps^{1/2}}
\biggr)^{nN^2-\sum_{i=1}^n\dim_\bR\widehat U_i^{(\ell)}(N)}
$$
as long as $\eps>0$ is small enough. Therefore we get
\begin{equation}\label{F-5.25}
\mu_N^{(\ell)}(\mathrm{Ball}((L+1)\eps^{1/2},d_{2,E}))
\le\biggl({C_\ell(L+1)\eps^{1/2}\over C'_\ell}
\biggr)^{nN^2-\sum_{i=1}^n\dim_\bR\widehat U_i^{(\ell)}(N)}
\end{equation}
for all sufficiently small $\eps>0$.

Combining \eqref{F-5.23}--\eqref{F-5.25} and
\eqref{F-5.21} altogether implies that
\begin{align*}
&\limsup_{N\to\infty}{1\over N^2}\log
\gamma_{\U(N)^n}(\Gamma_\orb(\eps,N,3m,\delta)) \\
&\quad\le
\limsup_{N\to\infty}{1\over N^2}\log
K_{\eps^{1/2}}(\Gamma(N,m,\delta)) \\
&\qquad\qquad +\sum_{i=1}^n\Biggl(1-{\tau(p_{i0})^2\over\ell}
-\sum_{k=1}^{s(i)\wedge\ell}{\tau(p_{ik})^2\over m_{ik}^2}-\tau(q_i^{(\ell)})^2\Biggr)
\cdot\biggl(\log\eps^{1/2}+\log{C_\ell(L+1)\over C'_\ell}\biggr)
\end{align*}
whenever $\eps>0$ is sufficiently small. Take
$\lim_{m\to\infty,\delta\searrow0}$ and then $\limsup_{\eps\searrow0}$ after dividing
by $|\log\eps^{1/2}|$; then by Lemma \ref{L-5.9}, \eqref{F-5.18},
Proposition \ref{P-5.6} and \eqref{F-5.19} we have
$$
\delta_{0,\orb}(\bX_1,\dots,\bX_n)
\le\delta_0(\bX_1\sqcup\cdots\sqcup\bX_n)
-\sum_{i=1}^n\Biggl(1-{\tau(p_{i0})^2\over\ell}
-\sum_{k=1}^{s(i)\wedge\ell}{\tau(p_{ik})^2\over m_{ik}^2}-\tau(q_i^{(\ell)})^2\Biggr).
$$
Hence the inequality $\ge$ in Theorem \ref{T-5.8} follows by taking 
$\ell\rightarrow\infty$ thanks to \eqref{F-5.8}.

\section{Applications}
\setcounter{equation}{0}

\subsection{Immediate corollaries} The next corollary is immediate from
Theorem \ref{T-5.8} and (5)--(7) of Proposition \ref{P-5.3}.

\begin{cor}\label{C-6.1}
Let $\bX_1,\dots,\bX_n$ be hyperfinite random multi-variables in $(M,\tau)$.
\begin{itemize}
\item[(1)] If $\bY_1,\dots,\bY_n$ are random multi-variables such that
$\bY_i\subset W^*(\bX_i)$ for $1\le i\le n$, then
$$
\delta_0(\bX_1\sqcup\cdots\sqcup\bX_n)-\sum_{i=1}^n\delta_0(\bX_i)
\le\delta_0(\bY_1\sqcup\cdots\sqcup\bY_n)-\sum_{i=1}^n\delta_0(\bY_i).
$$
In addition, if $W^*(\bY_i)$'s are all diffuse, then
$\delta_0(\bX_1\sqcup\cdots\sqcup\bX_n)\le\delta_0(\bY_1\sqcup\cdots\sqcup\bY_n)$.
\item[(2)] If $\chi_\orb(\bX_1,\dots,\bX_n)>-\infty$, then
$\delta_0(\bX_1\sqcup\cdots\sqcup\bX_n)=\sum_{i=1}^n\delta_0(\bX_i)$. In particular,
this is the case if $\bX_1,\dots,\bX_n$ are freely independent.
\item[(3)] Let $\bX$ be a hyperfinite random multi-variable and $\bY$ a general random
multi-variable. If $\bX$ is freely independent of $\bY$, then
$\delta_0(\bX\sqcup\bY)=\delta_0(\bX)+\delta_0(\bY)$.
\end{itemize}
\end{cor}

The following is an immediate corollary of Theorem \ref{T-5.8} too. But we note that
it can also be shown by a direct method of estimating the covering numbers of orbital
microstate spaces.     

\begin{cor}\label{C-6.2} {\rm (General upper bound of $\delta_{0,\mathrm{orb}}$)}
\enspace For any $n$-tuple of hyperfinite self-adjoint multi-variables
$\mathbf{X}_1,\dots,\mathbf{X}_n$,
\begin{equation*}
\delta_{0,\mathrm{orb}}(\mathbf{X}_1,\dots,\mathbf{X}_n)
\le -(n-1)\,\delta_0\big(W^*(\mathbf{X}_1)\cap\cdots\cap W^*(\mathbf{X}_n)\big), 
\end{equation*}
and equality holds either when the $\mathbf{X}_k$'s are the same or when the
$\mathbf{X}_k$'s are freely independent with amalgamation over their common subalgebra
$W^*(\mathbf{X}_1)\cap\cdots\cap W^*(\mathbf{X}_n)$. Here,
$\delta_0\big(W^*(\mathbf{X}_1)\cap\cdots\cap W^*(\mathbf{X}_n)\big)$ denotes
the unique value of $\delta_0(\mathbf{X})$ with
$W^*(\mathbf{X}) = W^*(\mathbf{X}_1)\cap\cdots\cap W^*(\mathbf{X}_n)$ due to \cite{Ju1}.  
\end{cor}
\begin{proof} By Proposition \ref{P-5.3}\,(5) we have 
$$
\delta_{0,\mathrm{orb}}(\mathbf{X}_1,\dots,\mathbf{X}_n) \le 
\delta_{0,\mathrm{orb}}(\mathbf{X},\dots,\mathbf{X}), 
$$ 
and the right-hand side is $-(n-1)\,\delta_0(\mathbf{X})$ by Theorem \ref{T-5.8}. The
second equality condition follows from \cite{BDJ}. 
\end{proof}

\subsection{Liberation process vs $\delta_0$} Let $X_1,\dots,X_n$ be an $n$-tuple of
single self-adjoint random variables and $S_1,\dots,S_n$ a standard semicircular system
freely independent of $X_1,\dots,X_n$. Concerning the condition of $X_1,\dots,X_n$
having f.d.a.\ (i.e., finite-dimensional approximants), it is known (see
\cite[7.3.9]{HP} and \cite{Ju2}) that the following are all equivalent:
\begin{itemize}
\item $X_1,\dots,X_n$ has f.d.a.\ (i.e., finite-dimensional approximants);
\item $\chi(X_1+\eps S_1,\dots,X_n+\eps S_n)>-\infty$ for all $\eps>0$;
\item $\delta(X_1,\dots,X_n)\ge0$ (also $\delta_0(X_1,\dots,X_n)\ge0$);
\item $\delta(X_1+\eps S_1,\dots,X_n+\eps S_n)=n$ (also
$\delta_0(X_1+\eps S_1,\dots,X_n+\eps S_n)=n$) for all $\eps>0$.
\end{itemize}
The next proposition gives similar equivalent conditions in terms of the orbital theory.

\begin{prop}\label{P-6.3}
Let $\bX_1,\dots,\bX_n$ be hyperfinite random multi-variables, and
$\bv(t) = (v_i(t))_{i=1}^n$, $t \ge 0$, be a freely independent $n$-tuple of free
unitary Brownian motions that is freely independent of $\bX_1\sqcup\cdots\sqcup\bX_n$.
Then the following conditions are equivalent:
\begin{itemize}
\item[(i)] $\bX_1\sqcup\cdots\sqcup\bX_n$ has f.d.a.
\item[(ii)] $\chi_\orb(v_1(\eps)\bX_1v_1(\eps)^*,\dots,v_n(\eps)\bX_nv_n(\eps)^*)
>-\infty$ for all $\eps>0$.
\item[(iii)] $\delta_{0,\orb}(\bX_1,\dots,\bX_n)>-\infty$.
\item[(iv)] $\delta_{0,\orb}(\bX_1,\dots,\bX_n)\ge-n$.
\item[(v)] $\delta_{0,\orb}(v_1(\eps)\bX_1v_1(\eps)^*,\dots,v_n(\eps)\bX_n v_n(\eps)^*)
=0$ for all $\eps>0$.
\end{itemize}
\end{prop}

\begin{proof}
(i) $\Rightarrow$ (iv).\enspace
Note that (i) implies $\delta_0(\bX_1\sqcup\cdots\sqcup\bX_n)\ge0$ as mentioned above.
Hence $\delta_{0,\orb}(\bX_1,\dots,\bX_n)\ge-\sum_{i=1}^n\delta_0(\bX_i)\ge-n$ by
Theorem \ref{T-5.8}.

(iv) $\Rightarrow$ (iii) is trivial.

(iii) $\Rightarrow$ (ii).\enspace
Condition (iii) implies that there is a sequence $\eps_k\searrow0$ such that
\begin{align*}
&\chi_\orb(v_1(\eps_k)\bX_1v_1(\eps_k)^*,\dots,v_n(\eps_k)\bX_nv_n(\eps_k)^*) \\
&\qquad
\ge \chi_\orb(v_1(\eps_k)\bX_1v_1(\eps_k)^*,\dots,v_n(\eps_k)\bX_nv_n(\eps_k)^*:
\bv(\eps_k)) > -\infty.
\end{align*}
Hence (ii) follows if we see that
$\chi_\orb(v_1(t)\bX_1 v_1(t)^*,\dots,\allowbreak v_n(t)\bX_n v_n(t)^*)$ is increasing
in $t\ge0$. To show this, let $t>s\ge0$ and set $w_i(t,s):=v_i(t)v_i(s)^*$,
$\bY_i:=v_i(s)\bX_iv_i(s)^*$. Note \cite{Bi} that $w_i(t,s)$ has the same distribution
as $v_i(t-s)$ and is freely independent of $\bY_1\sqcup\cdots\sqcup\bY_n$.
Hence we have
\begin{align*}
&\chi_\orb(v_1(t)\bX_1v_1(t)^*,\dots,v_n(t)\bX_nv_n(t)^*) \\
&\qquad=\chi_\orb(w_1(t,s)\bY_1w_1(t,s)^*,\dots,w_n(t,s)\bY_nw_n(t,s)^*) \\
&\qquad\ge\chi_\orb(\bY_1,\dots,\bY_n)
\end{align*}
thanks to Proposition \ref{P-4.6}.

(ii) $\Rightarrow$ (v) is Proposition \ref{P-5.3}\,(6).

(v) $\Rightarrow$ (i).\enspace
From (v) there is a sequence $\eps_k\searrow0$ such that
$$
\chi_\orb(v_1(\eps_k)\bX_1v_1(\eps_k)^*,\dots,v_n(\eps_k)\bX_nv_n(\eps_k)^*)>-\infty.
$$
Hence, for every $m\in\bN$, $\delta>0$ and $k\in\bN$ we have
$$
\Gamma_\orb((v_i(\eps_k)\bX_iv_i(\eps_k)^*)_{i=1}^n:
(\Xi_i(N))_{i=1}^n;N,m,\delta/2)\ne\emptyset
$$
for sufficiently large $N$, where $\Xi_i(N)=(\xi_{i1}(N),\dots,\xi_{ir(i)}(N))$,
$1\le i \le n$, are chosen as in Definition \ref{D-4.5} with
$\|\xi_{ij}(N)\|_\infty \le \|X_{ij}\|_\infty$. Since
$v_i(\eps_k)X_{ij}v_i(\eps_k)\to X_{ij}$ strongly as $k\to\infty$ for $1\le j\le r(i)$
with $\bX_i = (X_{i1},\dots,X_{ir(i)})$, $1\le i \le n$, this implies
that $\Gamma_R(\bX_1\sqcup\cdots\sqcup\bX_n;N,m,\delta)\ne\emptyset$ for every
$R\ge\max\{\|X_{ij}\|_\infty:1\le j \le r(i), 1\le i\le n\}$ if $N$ is sufficiently
large, so (i) follows.
\end{proof}

The next corollary is immediate from Theorem \ref{T-5.8} and the above proposition. 

\begin{cor}\label{C-6.4}
If $\bX_1\sqcup\cdots\sqcup\bX_n$ has f.d.a., then 
$$
\delta_0(v_1(\eps)\bX_1v_1(\eps)^*\sqcup\cdots\sqcup v_n(\eps)\bX_nv_n(\eps)^*)
=\sum_{i=1}^n\delta_0(\bX_i)
$$
for every $\eps>0$. In particular, $\delta_0$ is discontinuous at $0$ in the liberation
process. 
\end{cor}

It might be worth mentioning that the above corollary provides another route toward
Brown's observation \cite{Br} based on the liberation process instead of the
semicircular deformation.   

\subsection{Lower semicontinuity for $\delta_0$} The next lemma partially
strengthens \cite[Lemma 7.3]{Ju1}. 

\begin{lemma}\label{L-6.5} Let $\bX = (X_1,\dots,X_r)$ be hyperfinite self-adjoint
multi-variables. If $\bX^{(k)} = (X_1^{(k)},\dots,X_r^{(k)})$ is a sequence
of {\rm (}not necessarily hyperfinite{\rm)} random multi-variables converging to $\bX$
in mixed moments as $k\to\infty$, then one has
$$
\liminf_{k\to\infty}\delta_0(\bX^{(k)})\ge\delta_0(\bX).
$$ 
\end{lemma}
\begin{proof}
As same as in the proof of Theorem \ref{T-5.8} we assume that $W^*(\bX)$ has both
diffuse and atomic parts, and decompose
$$
W^*(\mathbf{X})=\bigoplus_{j=0}^s M_j
$$
with possibly $s=\infty$, where $M_0$ is diffuse and $M_j \cong M_{m_j}(\bC)$ for
$j\ge1$. Let $p_j$ be the central support projection of $M_j$. Fix an arbitrary
$\ell\in\bN$, and one can choose projections
$p_{01}^{(\ell)},\dots,p_{0\ell}^{(\ell)} \in M_0$ in such a way that
$p_{01}^{(\ell)}+\dots+p_{0\ell}^{(\ell)}=p_0$ and
$\tau(p_{01}^{(\ell)})=\cdots=\tau(p_{0\ell}^{(\ell)})$. Define a finite-dimensional
subalgebra $\cP_\ell \subset W^*(\mathbf{X})$ by
$$
\cP_\ell:=\bigoplus_{i=1}^\ell \bC p_{0i}^{(\ell)}
\,\oplus\,\bigoplus_{j=1}^{s\wedge\ell} M_j\,\oplus\bC q^{(\ell)},
$$
with $q^{(\ell)} := \1-\sum_{j=0}^{s\wedge\ell}p_j$. Here we choose and fix a matrix
unit system $\{e_{\alpha\beta}^{(j)}\}_{1\le\alpha,\beta\le m_j}$ of
$M_j \cong M_{m_j}(\bC)$ for $1\le j\le s\wedge\ell$. For any element $e$ chosen from 
$p_{0i}^{(\ell)}$ ($1\le i\le \ell$), $e_{\alpha\beta}^{(j)}$
($1\le\alpha,\beta\le m_j$, $1\le j\le s\wedge\ell$) and $q_\ell$ (altogether
forming a matrix unit system of $\cP_\ell$) and for any $\delta>0$, one can choose a
non-commutative polynomial $P$ of $r$ indeterminates such that
$\|e-P(X_1,\dots,X_r)\|_{\tau,2}<\delta$. Then
$\|e-P(X_1^{(k)},\dots,X_r^{(k)})\|_{\tau,2}<\delta$ for all sufficiently large $k$.
This shows that, for every $\eps>0$, one has $\cP_\ell \subset_\delta W^*(\bX^{(k)})$
for all $k$ large enough. (See \cite{CS} for the notation ``$\subset_\delta$".) Hence,
by \cite[Lemma 4]{CS} we see that, for every $\ell\in\bN$ and every $\eps>0$, there is
a $k_0\in\bN$ such that for every $k\ge k_0$ one can find a new matrix unit system 
(depending on $k$) $\hat{p}_{0i}$ ($1\le i\le\ell$), $\hat{e}_{\alpha\beta}^{(j)}$
($1\le\alpha,\beta\le m_j$, $1\le j\le s\wedge\ell$) and $\hat{q}$ inside
$W^*(\bX^{(k)})$ satisfying
$$
\sum_{i=1}^\ell \hat{p}_{0i}
+ \sum_{j=1}^{s\wedge\ell} \sum_{\alpha=1}^{m_j}\hat{e}_{\alpha\alpha}^{(j)}
+ \hat{q} = \1
$$
and
\begin{align*}
&\|\hat{p}_{0i}-p_{0i}^{(\ell)}\|_{2,\tau}<\eps,\quad 1\le i\le \ell, \\
&\|\hat{e}_{\alpha1}^{(j)}-e_{\alpha1}^{(j)}\|_{2,\tau} < \eps,
\quad 1\le\alpha\le m_j,\ 1\le j\le s\wedge\ell, \\
&\|\hat{q} - q^{(\ell)}\|_{2,\tau}< \eps.
\end{align*}
For every $k\ge k_0$, since $W^*(\bX^{(k)})$ has a subalgebra (depending on $k$)
$$
\bigoplus_{i=1}^\ell \bC\hat{p}_{0i}\,
\oplus\,\bigoplus_{j=1}^{s\wedge\ell}
\mathrm{Alg}\bigl(\{\hat{e}_{\alpha\beta}^{(j)} : 1\le\alpha,\beta\le m_j\}\bigr)\,
\oplus\bC\hat{q},
$$
it follows from \cite[Corollary 7.2]{Ju1} that
\begin{align}\label{F-6.1}
\delta_0(\bX^{(k)})
&\ge
1-\sum_{i=1}^\ell\tau(\hat{p}_{0i})^2
-\sum_{j=1}^{s\wedge\ell}\frac{\tau\big(\sum_{\alpha=1}^{m_j}
\hat{e}_{\alpha\alpha}^{(j)}\big)^2}{m_j^2}-\tau(\hat{q})^2 \notag \\
&\ge 1-\ell\biggl(\frac{\tau(p_0)}{\ell}+\eps\biggr)^2
-\sum_{j=0}^{s\wedge\ell}\frac{(\tau(p_j)+2m_j\eps)^2}{m_j^2} 
-\Biggl(\tau\Biggl(\1-\sum_{j=0}^{s\wedge\ell}p_j\Biggr)+\varepsilon\Biggr)^2,
\end{align}
since 
\begin{align*}
\Bigg\|\sum_{\alpha=1}^{m_j}\hat{e}_{\alpha\alpha}^{(j)}-p_j\Bigg\|_{2,\tau}
&\le \sum_{\alpha=1}^{m_j}\|\hat{e}_{\alpha1}^{(j)}\hat{e}_{1\alpha}^{(j)}
-e_{\alpha1}^{(j)}e_{1\alpha}^{(j)}\|_{2,\tau} \\
&\le \sum_{\alpha=1}^{m_j}\Big\{\|(\hat{e}_{\alpha1}^{(j)}
-e_{\alpha1}^{(j)})\hat{e}_{1\alpha}^{(j)}\|_{2,\tau}
+\|e_{\alpha1}^{(j)}(\hat{e}_{1\alpha}^{(j)}-e_{1\alpha}^{(j)})\|_{2,\tau}\Big\} \\
&\le 2m_j\eps.
\end{align*}
Since the value of \eqref{F-6.1} converges to
$\delta_0(\bX)=1-\sum_{j=1}^s \tau(p_j)^2/m_j^2$ as $\eps\searrow0$ and then
$\ell\to\infty$, we get the desired assertion.  
\end{proof}

\begin{remark}
If $\bX^{(k)} \subset W^*(\bX)$ for all $k$ is further assumed in the above lemma,
the consequence becomes $\lim_{k\rightarrow\infty}\delta_0(\bX^{(k)}) = \delta_0(\bX)$
due to \cite[Corollary 7.2]{Ju1}. 
\end{remark} 

\begin{prop}\label{P-6.7} {\rm (Lower semicontinuity of $\delta_0$ under convergence
inside hyperfinite subalgebras)}\enspace Let
$\bX_i=(X_{i1},\dots,X_{ir(i)})$, $1\le i\le n$, be hyperfinite random
multi-variables. For each $1\le i \le n$ assume that
$\bX_i^{(k)}=(X_{i1}^{(k)},\dots,X_{ir(i)}^{(k)}) \subset W^*(\bX_i)$
is a sequence of hyperfinite random multi-variables converging to $\bX_i$ in mixed
moments as $k\to\infty$. Then one has
\begin{equation*}
\liminf_{k\rightarrow\infty}\delta_0(\bX^{(k)}_1\sqcup\cdots\sqcup\bX^{(k)}_n)
\ge \delta_0(\bX_1\sqcup\cdots\sqcup\bX_n). 
\end{equation*}
\end{prop}   
\begin{proof} 
We have   
\begin{align*} 
\liminf_{k\rightarrow\infty}
\delta_0(\mathbf{X}_1^{(t)}\sqcup\cdots\sqcup\mathbf{X}_n^{(k)}) 
&= \liminf_{k\rightarrow\infty}\Biggl(\delta_{\orb,0}(\bX_1^{(k)},\dots,\bX_n^{(k)})
+ \sum_{k=1}^n\delta_0(\bX_k^{(k)})\Biggr) \\
&\ge \liminf_{k\rightarrow\infty}\Biggl(\delta_{\orb,0}(\bX_1,\dots,\bX_n)
+ \sum_{i=1}^n\delta_0(\bX_i^{(k)})\Biggr) \\
&\ge \delta_{\orb,0}(\bX_1,\dots,\bX_n)
+ \sum_{i=1}^n \liminf_{k\rightarrow\infty}\delta_0(\bX_i^{(k)}) \\
&\ge \delta_{\orb,0}(\bX_1,\dots,\bX_n) + \sum_{i=1}^n \delta_0(\mathbf{X}_i) \\
&= \delta_0(\bX_1\sqcup\cdots\sqcup\bX_n), 
\end{align*}
where the first and the last equalities are due to Theorem \ref{T-5.8}, the second
inequality is Proposition \ref{P-5.3}\,(5) and the last inequality is Lemma
\ref{L-6.5}.
\end{proof}

\begin{remark} It is worth pointing out that the invariance of $\delta_0$
for hyperfinite von Neumann algebras, a corollary of Jung's result \cite{Ju1}, can
also be shown from Lemma \ref{L-6.5} and the argument in \cite[Remark 6.13]{V2}
(together with \cite[Theorem 4.3 and Corollary 4.5]{V-IMRN}). Similarly, Proposition
\ref{P-6.7} and the same pattern of argument show that if
$\bX_1,\dots,\bX_n$ and $\bX_1',\dots,\bX_n'$ are hyperfinite random multi-variables
and $W^*(\bX_i) = W^*(\bX_i')$ for $1\le i \le n$, then
$\delta_0(\bX_1\sqcup\cdots\sqcup\bX_n) = \delta_0(\bX_1'\sqcup\cdots\sqcup\bX_n')$
holds. In fact, this is a ``part" of \cite[Corollary 4.2]{Ju3} and also a consequence
of Theorem \ref{T-5.8} and Proposition \ref{P-5.3}\,(4). However, those two results
of Jung \cite{Ju1,Ju3} were crucially used in the proofs of Lemma \ref{L-6.5} and
Theorem \ref{T-5.8} while it is desirable to prove Theorem \ref{T-5.8}
without the use of the hyperfinite inequality in \cite{Ju3}. Our discussion on lower
semicontinuity for $\delta_0$ should be regarded as a kind of ``converse argument" of
\cite[Remark 6.13]{V2}, thus suggesting that a kind of lower semicontinuity for
$\delta_0$ is essentially equivalent to the affirmative solution of ``the entropy
dimension problem" in \cite[\S2.6]{V-Survey}.
\end{remark}

\subsection{A few remarks on definitions of $\delta_0$ and
$\delta_{0,\orb}$} In the course of proving
Theorem \ref{T-5.8} we examined several ideas and observed some small facts concerning
the free entropy dimension $\delta_0$ and Voiculescu's liberation process, which may
be of independent interest. Here we give a brief summary of them.

Let $X_1,\dots,X_n$ be non-commutative self-adjoint random variables in a tracial
$W^*$-probability space, and we choose a freely independent $n$-tuple of free
multiplicative unitary Brownian motions $\mathbf{v}(t) = (v_k(t))_{k=1}^n$, $t \ge 0$.
The next proposition can be shown directly by an argument as in the part
``$\le$" of Proposition \ref{P-5.6} applied to $\chi$ rather than $\chi_\mathrm{orb}$
while we will derive it simply from Theorem \ref{T-5.8}.

\begin{prop}\label{P-6.9} With the above assumption, 
$$
n+\limsup_{\eps\searrow0}
\frac{\chi(v_1(\eps)X_1 v_1(\eps)^*,\dots,v_n(\eps) X_n v_n(\eps)^* :
\mathbf{v}(\eps))}{|\log\eps^{1/2}|} \le \delta_0(X_1,\dots,X_n).  
$$
Furthermore, if $\chi(X_i) > -\infty$ for all $1\le i\le n$, then equality holds and
$\delta_0(X_1,\dots,X_n) = \delta_{0,\mathrm{orb}}(X_1,\dots,X_n) + n$.  
\end{prop}
\begin{proof} One has nothing to do when $\chi(X_i) = -\infty$ for some $i$ since the
$\chi$ in the left-hand side always becomes $-\infty$. Thus we may and do assume
$\chi(X_i) > -\infty$ for all $1\le i\le n$. By Remark \ref{R-5.2} the left-hand side
is nothing but $n+\delta_{0,\mathrm{orb}}(X_1,\dots,X_n)$, and hence the assertion
follows from Theorem \ref{T-5.8} thanks to $\delta_0(X_i)=1$ (see
\cite[Proposition 6.3]{V2} and \cite[Corollary 6.7]{V3}).  
\end{proof}

Hence, no difference occurs even if the semicircular deformation is replaced by the
liberation process in the definition of $\delta_0$ when all $\chi(X_i)$ are finite.
Furthermore, we can prove the following:  

\begin{prop}\label{P-6.10} If $\mathbf{S}=(S_i)_{i=1}^n$ is a free semicircular system
freely independent of the other random variables, then
\begin{align*}
&\delta_0(X_1,\dots,X_n) \\
&\quad= n\, +\limsup_{\eps\searrow0}
\frac{\chi(v_1(\eps)(X_1+\eps^{1/2}S_1)v_1(\eps)^*,\dots,
v_n(\eps)(X_n+\eps^{1/2}S_n)v_n(\eps)^* :
\mathbf{S}\sqcup\mathbf{v}(\eps))}{|\log\eps^{1/2}|}.
\end{align*}
\end{prop}
\begin{proof} (Sketch)\enspace Fix $0<\delta<\eps\le1$ arbitrarily and
let $R\ge\max_{1\le i\le n}\Vert X_i\Vert_\infty +2$. We write
\begin{align*}
\Gamma(\eps,N,3m,\delta/2^m) 
&:= 
\Gamma_R((v_i(\eps)(X_i+\eps^{1/2}S_i)v_i(\eps)^*)_{i=1}^n:
\mathbf{S}\sqcup\mathbf{v}(\eps);N,3m,\delta/2^m), \\
\Gamma(N,m,\delta) 
&:= 
\Gamma_R((X_i)_{i=1}^n;N,m,\delta).
\end{align*}
By an essentially same argument as in the part ``$\le$" of Proposition \ref{P-5.6}
one can prove that $\Gamma(\eps,N,3m,\delta/2^m)$ is included in
$\cN_{L\eps^{1/2}}(\Gamma(N,m,\delta))$ with respect to $d_2$ (the metric induced from
$\tr_N$) with some constant $L>0$. Then one has 
\begin{align*}
\Lambda_N^{\otimes n}(\Gamma(t,N,3m,\gamma/2^m)) \le 
K_{\eps^{1/2}}(\Gamma(N,m,\gamma))\cdot
\Lambda_N^{\otimes n}(\text{Ball}((L+1)\eps^{1/2},d_2).
\end{align*}
From this one can derive the inequality $\ge$ for the desired equality similarly to
the part ``$\le$" of Proposition \ref{P-5.6}. On the other hand, we can modify the
proof of Proposition \ref{P-4.6} to show the inequality
\begin{align*}
&\chi(v_1(\eps)(X_1+\eps^{1/2}S_1)v_1(\eps)^*,\dots,
v_n(\eps)(X_n+\eps^{1/2}S_n)v_n(\eps)^*:
\mathbf{S}\sqcup\mathbf{v}(t)) \\
&\qquad\ge\chi(X_1+\eps^{1/2}S_1,\dots,X_n+\eps^{1/2}S_n:\mathbf{S}),
\end{align*}
which gives the reverse inequality.
\end{proof}

\begin{prop}\label{P-6.11}
With the same assumption as Proposition \ref{P-6.10},
\begin{align*}
&\delta_{0,\mathrm{orb}}(X_1,\dots,X_n) \\
&\qquad=\limsup_{\eps\searrow0}
{\chi_{\mathrm{orb}}(v_1(\eps)(X_1+\eps^{1/2}S_1)v_1(\eps)^*,\dots,
v_n(\eps)(X_n+\eps^{1/2}S_n)v_n(t)^*:\mathbf{S}\sqcup\mathbf{v}(\eps))
\over|\log\eps^{1/2}|}.
\end{align*}
\end{prop}

\begin{proof} (Sketch)\enspace The proofs of Lemma \ref{L-2.4} and Theorem \ref{T-2.6}
can slightly be modified to show the equality
\begin{align*}
&\chi(v_1(\eps)(X_1+\eps^{1/2}S_1)v_1(\eps)^*,\dots,
v_n(\eps)(X_n+\eps^{1/2}S_n)v_n(\eps)^*:
\mathbf{S}\sqcup\mathbf{v}(\eps)) \\
&\qquad=\chi_{\mathrm{orb}}(v_1(\eps)(X_1+\eps^{1/2}S_1)v_1(\eps)^*,\dots,
v_n(\eps)(X_n+\eps^{1/2}S_n)v_n(\eps)^*:\mathbf{S}\sqcup\mathbf{v}(\eps)) \\
&\qquad\qquad+\sum_{i=1}^n\chi(X_i+\eps^{1/2}S_i).
\end{align*}
Hence the desired formula follows from Theorem \ref{T-5.8} and Proposition \ref{P-6.10}
thanks to $\delta_0(X_i)=1+\lim_{\eps\searrow0}\chi(X_i+\eps^{1/2}S_i)/|\log\eps^{1/2}|$
by \cite[Proposition 6.3]{V2}.
\end{proof}

We were interested in the formulas in Propositions \ref{P-6.10} and \ref{P-6.11}
because those together immediately imply Theorem \ref{T-5.8} for self-adjoint
variables $X_1,\dots,X_n$. But the direct proof of Proposition \ref{P-6.11}
seems difficult since the orbital theory does not fit well in additive operations like
the semicircular deformation.

\end{document}